\documentclass[11pt]{article}

\usepackage{amsfonts}
\usepackage{amssymb}
\usepackage{bm}
\usepackage{amsmath}
\usepackage{amsthm}
\usepackage{stmaryrd}
\usepackage{color}
\usepackage{xcolor}
\usepackage[subnum]{cases}
\usepackage{rotating}
\usepackage{enumitem}
\usepackage{accents}
\usepackage[title]{appendix}
\usepackage{mathtools}
\usepackage{caption}
\usepackage{url}
\usepackage{framed}
\usepackage{lscape}
\usepackage{multirow}
\usepackage{latexsym}
\usepackage{mathrsfs}
\usepackage{array}
\usepackage{makecell}
\usepackage{float}
\usepackage{tabularray}% for long table

\usepackage{tablefootnote}
\usepackage{tikz}

\usepackage[usestackEOL]{stackengine}
\usepackage{pdflscape}
\usepackage{geometry}
\newtheorem{proposition}{Proposition}[section]

\newtheorem{definition}{Definition}[section]
\newtheorem{lemma}{Lemma}[section]
\newtheorem{corollary}{Corollary}[section]
\newtheorem{theorem}{Theorem}[section]
\newtheorem{remark}{Remark}[section]

\newtheorem{example}{Example}[section]

\definecolor{lightgreen}{rgb}{0.1,0.5,0.1}
\usepackage[colorlinks=true,
breaklinks=true,
bookmarks=true,
urlcolor=blue,
citecolor=lightgreen,
linkcolor=lightgreen,
bookmarksopen=false,
draft=false]{hyperref}

\usepackage{endnotes}
\let\footnote=\endnote

\usepackage{graphicx}
\usepackage{multirow}
\usepackage{hhline}

\usepackage{appendix}

\usepackage{color}

\usepackage{amssymb}
\usepackage{mathtools}
\usepackage{stmaryrd}
\expandafter\def\expandafter\normalsize\expandafter{%
    \normalsize%
    \setlength\abovedisplayskip{5pt}%
    \setlength\belowdisplayskip{3pt}%
}

\usepackage{multirow}
\usepackage{array}
\usepackage{makecell}
\usepackage{float}
\usepackage{tabularray}
\usepackage{tablefootnote}
\usepackage{tikz}

\usepackage{enumitem}

\usepackage{url}

\usepackage{natbib}
\bibpunct[, ]{(}{)}{,}{a}{}{,}%
\def\BIBand{and}%

\begin{document}

\title{Wasserstein distributionally robust optimization\protect\\ and its tractable regularization formulations}

\author{Hong T. M. Chu\footnotemark[1],\quad Meixia Lin\footnotemark[2]\,\,\,\footnotemark[4], \quad Kim-Chuan Toh\footnotemark[3]}

\maketitle

\renewcommand{\thefootnote}{\fnsymbol{footnote}}

\footnotetext[1]{Department of Mathematics, National University of Singapore, Singapore, \texttt{hongtmchu@u.nus.edu}}
\footnotetext[2]{Engineering Systems and Design, Singapore University of Technology and Design, Singapore, \texttt{meixia\_lin@sutd.edu.sg}}
\footnotetext[3]{Department of Mathematics and Institute of Operations Research and Analytics, National University of Singapore, Singapore, \texttt{mattohkc@nus.edu.sg}}
\footnotetext[4]{Corresponding author}
\renewcommand{\thefootnote}{\arabic{footnote}}

\begin{abstract}
	We study a variety of Wasserstein distributionally robust optimization (WDRO) problems where the distributions in the ambiguity set are chosen by constraining their Wasserstein discrepancies to the empirical distribution. Using the notion of weak Lipschitz property, we derive lower and upper bounds of the corresponding worst-case loss quantity and propose sufficient conditions under which this quantity coincides with its regularization scheme counterpart. Our constructive methodology and elementary analysis also directly characterize the closed-form of the approximate worst-case distribution. Extensive applications show that our theoretical results are applicable to various problems, including regression, classification and risk measure problems.
\end{abstract}

\medskip
\noindent
{\bf AMS subject classification:} {90C15, 90C17, 90C47}
\\[5pt]

{\bf Keywords:} {Wasserstein discrepancy, Wasserstein distributionally robust optimization, regularized optimization, worst-case loss quantity, data-driven decision making.}

\section{Introduction}
A central question of interest in many machine learning and operations research applications is the selection of an appropriate decision variable \(\beta\) from a decision space \(\mathcal{B}\). This often involves minimizing the expected risk of prediction errors, that is, 
\begin{equation*}
\inf_{\beta\in\mathcal{B}} \ \mathrm{E}_{\mathbb{P}_{\mathrm{true}}} [\ell (Z;\beta) ],
\end{equation*}
where \(Z\) is a random variable in a given space \(\mathcal{Z}\), with the probability distribution \(\mathbb{P}_{\mathrm{true}}\), and \(\ell \colon \mathcal{Z}\times \mathcal{B} \rightarrow \mathbb{R} \) is a loss function. In practice, the ground-truth distribution \(\mathbb{P}_{\mathrm{true}} \) is usually unknown. Instead, one only has access to an empirical distribution \(\mathbb{P}_N\coloneqq\sum_{i=1}^{N}\mu_i {\boldsymbol \chi}_{\{Z^{(i)}\}} \), where $\mathcal{Z}_N\coloneqq \{Z^{(1)},\dots,Z^{(N)}\} \subset \mathcal{Z}$ is a training dataset, \(\{\mu_i\}_{i=1}^{N} \) are nonnegative weights satisfying \(\sum_{i=1}^{N}\mu_i=1\), and \( {\boldsymbol \chi}_{\{Z^{(i)}\}}\) is the point mass at  \(Z^{(i)}\). The associated optimization problem
\begin{align}
\inf_{\beta\in\mathcal{B}} \left\{ \mathrm{E}_{\mathbb{P}_N} [\ell (Z;\beta) ] = \sum_{i=1}^{N}\mu_i\ell(Z^{(i)};\beta) \right\}  \label{eq: erm}
\end{align}
is often known as the empirical risk minimization (ERM) problem \citep{vapnik2015uniform}.  

\textbf{Robustness approach.} One major criticism of the ERM problem is that the empirical distribution \(\mathbb{P}_N\) might differ from the ground-truth distribution \(\mathbb{P}_{\mathrm{true}} \) considerably, and \(\mathcal{Z}_N\) might be unreliable due to errors in the data collection. This motivated the study of the corresponding distributionally robust optimization problem. Rather than relying on one single distribution \(\mathbb{P}_N\), it hedges against a set of distributions \({\mathfrak{M}}\) in the space of all probabilities on $\mathcal{Z}$, denoted as \(\mathcal{P}(\mathcal{Z}) \). Formally, the distributionally robust optimization (DRO) problem takes the form of solving a minimax problem
\begin{equation*}
\inf_{\beta\in\mathcal{B}} \ \sup_{\mathbb{P}\in{\mathfrak{M}} }  \ \mathrm{E}_{\mathbb{P}} [\ell (Z;\beta) ].
\end{equation*}
Here, the set \({\mathfrak{M}}\subset\mathcal{P}(\mathcal{Z}) \) is referred to as the ambiguity set or uncertainty set, which is often constructed by imposing some statistical conditions on the set of probability distributions under consideration. For example, the moment-based ambiguity set can be defined via certain moment constraints as in \citet{delage2010distributionally,goh2010distributionally,wiesemann2014distributionally}, or as the confidence region of a goodness-of-fit hypothesis test \citep{bertsimas2018robust}. Alternatively, the ambiguity set can be constructed as \( {\mathfrak{M}} = \left\{\mathbb{P}\in \mathcal{P}(\mathcal{Z}) \mid \mathcal{D}(\mathbb{P},\mathbb{P}_N) \leq \delta  \right\} \), where \(\delta\) is a nonnegative tuning parameter, and \(\mathcal{D}(\cdot,\cdot)\) defines a certain discrepancy on \(\mathcal{P}(\mathcal{Z}) \), such as the Prohorov and total variation metric \citep{gibbs2002choosing,erdougan2006ambiguous}, Kullback-Leibler and \(\chi^2 \) divergence \citep{hu2013kullback,jiang2016data}, or Wasserstein metric \citep{shafieezadeh2015distributionally,mohajerin2018data,blanchet2019quantifying,gao2022finite}. These choices of \({\mathfrak{M}}\) are often named as discrepancy-based ambiguity sets. 

In this work, we focus on the ambiguity set based on the Wasserstein discrepancy. Given two probability distributions \(\mathbb{P},\mathbb{Q}\in\mathcal{P}(\mathcal{Z}) \) and an extended nonnegative-valued function \(d\colon \mathcal{Z}\times\mathcal{Z}\rightarrow [0,\infty] \), the Wasserstein discrepancy\footnotemark\footnotetext{In this work, we use the term ``Wasserstein discrepancy'' instead of the commonly-used ``Wasserstein metric'' in the literature, since \(d(\cdot,\cdot)\) here is not required to be a metric.} with respect to \(d(\cdot,\cdot)\) and an exponent  \(r\in[1,\infty)\) is defined via the optimal transport problem \citep{villani2009optimal,peyre2017computational} as
\begin{align}
\mathcal{W}_{d,r}(\mathbb{P},\mathbb{Q}) \coloneqq \left(\inf_{\pi \in \Pi(\mathbb{P},\mathbb{Q})} \int_{\mathcal{Z}\times\mathcal{Z}} d^{r}(z',z)\mathrm{d}\pi(z',z)\right)^{\frac{1}{r}},\label{eq:defWd}
\end{align}
where \(\Pi(\mathbb{P},\mathbb{Q})\) \cite[Definition 1.1]{villani2009optimal} denotes the set of all joint probability distributions between \(\mathbb{P}  \) and \(\mathbb{Q}\), that is, by letting \(\sigma(\mathcal{Z})\) denote the set of all measurable sets in \(\mathcal{Z}\),
\begin{equation*}
\Pi(\mathbb{P},\mathbb{Q})  =  \left\{
\begin{matrix}
\pi \in\mathcal{P}(\mathcal{Z}\times\mathcal{Z}) \,\text{such that}\,\forall\, A,B\in \sigma(\mathcal{Z})\\
\pi (A\times \mathcal{Z}) = \mathbb{P}(A),  \pi(\mathcal{Z}\times B) = \mathbb{Q}(B)
\end{matrix} \right\}.
\end{equation*} 
Intuitively, the Wasserstein discrepancy \eqref{eq:defWd} can be understood as finding the minimum cost to move the mass of \(\mathbb{P}\) to that of \(\mathbb{Q}\). Accordingly, the Wasserstein distributionally robust optimization (WDRO) problem considers the ambiguity set \({\mathfrak{M}}= \left\{ \mathbb{P}\in\mathcal{P}(\mathcal{Z}) \mid \mathcal{W}_{d,r}(\mathbb{P},\mathbb{P}_N) \leq \delta \right\}\), where \(\delta\) is a nonnegative scalar. Formally, the WDRO problem takes the form as
\begin{align} 
\inf_{\beta\in\mathcal{B}} \ \sup_{\mathbb{P} : \mathcal{W}_{d,r}(\mathbb{P},\mathbb{P}_N) \leq \delta }  \ \mathrm{E}_{\mathbb{P}} [\ell (Z;\beta) ].\label{eq:WDRO}
\end{align}

\textbf{Regularization approach.} Another criticism of the ERM problem is that the resulting estimator \(\hat{\beta}\) might exhibit unsatisfactory out-of-sample performance or overfitting phenomena \citep{plan2012robust,feng2014robust}. To overcome this deficiency, a common approach is to modify the objective function in the ERM problem \eqref{eq: erm} by adding a regularization term. Specifically, the regularization scheme takes the form as
\begin{equation*}
	\inf_{\beta\in\mathcal{B}} \ \mathrm{E}_{\mathbb{P}_N} [\ell (Z;\beta) ] + \delta \varpi(\beta),
\end{equation*}
where \(\delta\) is a nonnegative tuning parameter, $\varpi: \mathcal{B}\rightarrow (-\infty,+\infty]$ is a regularization function. When \(\mathcal{B}\subset\mathbb{R}^n\), there are some popular options for $\varpi(\cdot)$, such as the ridge penalty \(\|\cdot\|^2_2 \) \citep{horel1962application,hoerl1970ridge}, the Lasso penalty \(\|\cdot\|_1 \) or its variant combining with a group or a fused penalty \citep{belloni2011square,bunea2013group,stucky2017sharp,jiang2021simultaneous}. In addition, it is also a topic of interest to study the appropriate value for the tuning parameter \(\delta\). In certain scenarios such as the Lasso and group Lasso problems \citep{bickel2009simultaneous,lounici2011oracle}, it has been shown that {appropriate values}  of \(\delta\) {depend} on the noise level of the dataset. In response, several works have shown that the optimal \(\delta\) can be made independent of the noise level if the original expectation term \( \mathrm{E}_{\mathbb{P}_N}[\ell(Z;\beta)]\) is replaced by \(\left( \mathrm{E}_{\mathbb{P}_N}[\ell(Z;\beta)]\right)^{\frac{1}{2}}\) when \(\ell(\cdot)\) is the squared loss function  \citep{belloni2011square,bunea2013group,stucky2017sharp}.

\textbf{Equivalence interpretation.} In both the robustness and regularization approaches, one of the key ingredients is the tuning parameter \(\delta \), which essentially controls how conservative the new scheme is compared to the original ERM problem. While the regularization scheme is more tractable and preferable in terms of computational consideration, the WDRO scheme \eqref{eq:WDRO} is more favorable in accommodating the geometric structure of the data space via the cost function \(d(\cdot,\cdot)\) and the
intuitive understanding of the level of robustness via the radius \(\delta\). In order to draw the connections and take into account the advantages of both of them, the equivalence between the WDRO problem and the regularization scheme has received increasing attention 
over the last few years. Specifically, let
\begin{align}
\mathcal{S}\coloneqq \sup_{\mathbb{P}\colon  \mathcal{W}_{d,r}(\mathbb{P},\mathbb{P}_N) \leq \delta }  \mathrm{E}_{\mathbb{P}} [\ell (Z;\beta) ], \label{eq:defS}
\end{align} 
then it aims to study sufficient conditions under which the following equivalence holds:
\begin{align} 
\mathcal{S} = \left[ \left( \mathrm{E}_{\mathbb{P}_N}[\ell(Z;\beta)]\right)^{\frac{1}{r}} + {\boldsymbol L}(\beta) \delta\right]^{r}, \label{eq:main_equation}
\end{align}
where \(r\in[1,\infty) \), and ${\boldsymbol L}(\beta)$ acts as a penalty on \(\beta\) which might also depend on other factors such as \(\mathbb{P}_N\) and \(d(\cdot,\cdot)\). In particular, the equivalence \eqref{eq:main_equation} gives a probabilistic explanation of the penalty parameter in the regularization model based on the WDRO interpretation. {Compared to \eqref{eq:defS} which involves a minimax problem, solving problem \eqref{eq:main_equation} appears to be more tractable and computationally favorable in certain circumstances, thanks to
{the} efficient algorithms studied in \cite{li2018efficiently,li2018highly,luo2019solving,zhang2020efficient,tang2020sparse,Chu2021OnRS}.}

There are typically two main streams of works to tackle \eqref{eq:main_equation} in the literature. First, the worst-case loss quantity \(\mathcal{S}\) defined in \eqref{eq:defS} can be viewed as an optimization problem with a single inequality constraint, thus its dual counterpart can be inspected as a univariate optimization problem. In addition, the complications in verifying the interchangeability condition for \(\sup\) and \(\inf\) also suggest that it might be beneficial to replace \eqref{eq:defWd} with its dual problem. Consequently, to guarantee the equality \eqref{eq:main_equation} instead of just the inequality derived from weak duality, many existing works impose relatively strong assumptions on the underlying problem in order to prove/use the strong duality and/or guarantee the existence of the dual optimal variables \citep{mohajerin2018data,blanchet2019quantifying,blanchet2019robust,Chu2021OnRS,zhang2022simple,gao2023distributionally,Zhen2023}. Second, when \(r=1\), \eqref{eq:main_equation} can be rewritten as 
%\({\boldsymbol L}(\beta) \delta =  \sup\limits_{\mathbb{P}\colon  \mathcal{W}_{d,r}(\mathbb{P},\mathbb{P}_N) \leq \delta }\mathrm{E}_{(Z',Z)\sim \pi\in\Pi(\mathbb{P},\mathbb{P}_N)} [\ell (Z';\beta) - \ell (Z;\beta) ]\)
\({\boldsymbol L}(\beta) \delta =  \sup\left\{ \mathrm{E}_{(Z',Z)\sim \pi} [\ell (Z';\beta) - \ell (Z;\beta) ] \mid \pi\in\Pi(\mathbb{P},\mathbb{P}_N) \mbox{ with } \mathbb{P}\in\mathcal{P}(\mathcal{Z}),  \mathcal{W}_{d,r}(\mathbb{P},\mathbb{P}_N) \leq \delta  \right\} \).
This observation suggests that \({\boldsymbol L}(\beta) \) is related to certain Lipschitz-type properties of the loss function \(\ell(\cdot,\beta)\) \citep{shafieezadeh2015distributionally,shafieezadeh2019regularization,an2021generalization,gao2022finite}. However, the usual Lipschitz condition is not enough to guarantee \eqref{eq:main_equation} since the Lipschitz constant can always be chosen arbitrarily large. Thus, one often needs some other conditions such as the tightness at certain points \citep{shafieezadeh2019regularization,gao2022finite}, the convexity of \(\ell\) \citep{wu2022generalization} or the differentibility of $\ell$ almost-everywhere with nonexpansive gradients \citep{an2021generalization}. 

\vspace{0.3cm}
In this paper, we study  sufficient conditions to establish the equivalence between the worst-case loss quantity in the WDRO problem and its associated regularization scheme. Our proposed sufficient conditions generalize the existing results from various perspectives, particularly by relaxing the required assumptions on the loss function and cost function. Moreover, our constructive approaches and elementary proofs directly characterize the closed forms of the approximate worst-case distributions. The generality of our theoretical results are demonstrated through their applications to various problems, including regression, classification and risk measure problems.

The remaining part of this paper is organized as follows. In Section \ref{sec:contribution}, we summarize our main contributions and compare them with the existing results in the literature. We derive our main theoretical results in Section \ref{sec:main-results} and present their applications in Section \ref{sec:app} and Section \ref{sec: generalrisk}. The conclusion is given in Section \ref{sec: conclusion}.

\vspace{0.2cm}
\noindent\textbf{Notations.} Throughout this paper, \((\mathcal{Z},\mathcal{A}) \) or \(\mathcal{Z}\) denotes a measurable space, where \(\mathcal{A}\) is a given \(\sigma\)-algebra on \(\mathcal{Z}\) such that \(\{z\}\in\mathcal{A} \) for any \(z\in\mathcal{Z} \)  \cite[Section 1.2]{cohn2013measure}. In particular, when \( \mathcal{Z} \in \mathfrak{B}(\mathbb{R}^n) \), where \(\mathfrak{B}(\mathbb{R}^n)\) is the Borel \(\sigma\)-algebra on \(\mathbb{R}^n\), then \(\mathcal{A}\) is always understood as \(\mathcal{A}\coloneqq \left\{ A \mid A \subset \mathcal{Z}, A\in \mathfrak{B}(\mathbb{R}^n)  \right\} \). A set \(A \subset \mathcal{Z} \) is called measurable if \(A\in\mathcal{A}\);   a function \(f\colon \mathcal{Z}\rightarrow  [-\infty,\infty]  \) is called measurable if \( \{z \in\mathcal{Z} \mid f(z) \leq t\} \in\mathcal{A} \) for any \(t\in\mathbb{R} \) \cite[Proposition 2.1.1]{cohn2013measure}; and \(\mathcal{Z}\times\mathcal{Z}\) denotes the Cartesian product measurable space with \(\sigma\)-algebra \(\mathcal{A}\times\mathcal{A}\) \cite[Section 5.1]{cohn2013measure}. 

In this work, the function \(\ell \colon \mathcal{Z}\times\mathcal{B}\rightarrow \mathbb{R} \) is assumed that \(\ell(\cdot,\beta) \colon \mathcal{Z} \rightarrow \mathbb{R} \) is measurable for any \(\beta\in\mathcal{B}\). A function \(\mathbb{P} \colon \mathcal{A}\rightarrow [0,\infty] \) is a probability  if it is countably additive, \(\mathbb{P}(\emptyset)=0 \) and \(\mathbb{P}(\mathcal{Z})=1 \); the space of all probabilities on \(\mathcal{Z}\) is denoted by \(\mathcal{P}(\mathcal{Z}) \); and the expectation of a measurable function \(f\) of a real-valued random variable \(Z\) on \((\mathcal{Z},\mathcal{A},\mathbb{P}) \) is denoted by \( \mathrm{E}_{\mathbb{P}}[f(Z)] = \int_{\mathcal{Z}} f(z)\mathrm{d}\mathbb{P}(z) \) \cite[Section 10.1]{cohn2013measure}. 

Define the indicator function \({\boldsymbol\delta}_{S} \colon  \mathcal{Z}  \rightarrow \mathbb{R}\) of a set $S\subset\mathcal{Z} $ as ${\boldsymbol\delta}_{S}(z)=0$ if $z\in S$, and $\infty$ otherwise. Define the point mass function (Dirac measure) ${\boldsymbol\chi}_{\{\hat{z}\}}\in\mathcal{P}(\mathcal{Z}) $ at point $\hat{z}\in \mathcal{Z}$ as ${\boldsymbol\chi}_{\{\hat{z}\} }(A) = 1$ if $\hat{z} \in A$, and $0$ otherwise, for any measurable set \(A\subset \mathcal{Z}\). Given two functions \(f\colon \mathcal{X}\rightarrow \mathbb{R}\) and \(g\colon\mathcal{Y}\rightarrow\mathbb{R}\), we define the function \(f\otimes g\colon \mathcal{X}\times\mathcal{Y}\rightarrow\mathbb{R} \) as \((x,y) \rightarrow f(x)\cdot g(y) \) for any $(x,y) \in \mathcal{X}\times \mathcal{Y}$. We adopt the convention of extended arithmetic such that \(0\cdot\infty=0 \). Denote the inner product on \(\mathbb{R}^n\) by \(\langle x,y\rangle = \sum_{i=1}^{n}x_i y_i  \) for any \(x,y\in\mathbb{R}^n\). Let \(\|\cdot\|_{\mathbb{R}^{n}}\) be an arbitrary norm on \(\mathbb{R}^{n}\) and \(\|\cdot\|_{{\mathbb{R}^{n}},*} \) be its dual norm defined as \(\|x\|_{\mathbb{R}^n,*}\coloneqq  \max_{y\in\mathbb{R}^n} \left\{ \langle x,y\rangle \mid \|y\|_{\mathbb{R}^n} =1 \right\}  \). Given matrices \(A\in\mathbb{R}^{n_1\times n_2}, B \in\mathbb{R}^{n_1\times n_3}, C\in\mathbb{R}^{n_3\times n_2} \), we denote the horizontal concatenation of \(A\) and \(B\) by \([A,B]\in\mathbb{R}^{n_1\times(n_2+n_3)}\), and the vertical concatenation of \(A\) and \(C\) by \([A;C]\in\mathbb{R}^{(n_1+n_3)\times n_2}\). Let  \(\mathbb{R}_+ \coloneqq [0,\infty)  \). For any real number $t$, the sign function is defined as ${\rm sgn}(t)=-1$ if $t<0$, and ${\rm sgn}(t)=1$ otherwise.

\section{Main contributions} 
\label{sec:contribution}
In this section, we shall summarize our main contributions and compare them with the existing results in the literature. We first state some notations which will be used. Let \(\mathcal{Z}_N\coloneqq\{Z^{(1)},\dots,Z^{(N)}\}\subset\mathcal{Z} \) be a given dataset and \(\mathbb{P}_N\coloneqq\sum_{i=1}^{N}\mu_i {\boldsymbol \chi}_{\{Z^{(i)}\}}\in\mathcal{P}(\mathcal{Z})\) be the corresponding empirical distribution. In addition, let $r\in [1,\infty)$ be a scalar and \(d\colon \mathcal{Z}\times\mathcal{Z}\rightarrow [0,\infty] \) be a measurable function on \(\mathcal{Z}\times\mathcal{Z}\). Suppose the loss function $\ell:\mathcal{Z}\times \mathcal{B}\rightarrow \mathbb{R}$ takes the form as
\begin{equation*}
\ell\colon (z;\beta) \mapsto \psi_{\beta}^{r}(z), \text{ with } 
\begin{cases}
\psi_{\beta}\colon\mathcal{Z}\rightarrow \mathbb{R} \ \ \text{ if } r =  1, \\
\psi_{\beta}\colon\mathcal{Z}\rightarrow \mathbb{R}_+ \text{ if } r > 1.
\end{cases}
\end{equation*}
For notational simplicity, let \(\mathcal{I} \) and \(\mathcal{U}\) (depending on some scalar \(L_{\beta}^{\mathcal{Z}_N}\), which will be discussed in detail later in Section \ref{sec:main-results}) be defined as:
\begin{equation*}
\begin{array}{ll}
\mathcal{I}& \coloneqq \inf\limits_{\rho\geq 0} \left\{ \rho\delta^r +\mathrm{E}_{\mathbb{P}_N} \left[ \sup\limits_{z'\in\mathcal{Z}} \left\{ \ell(z';\beta) - \rho d^r(z',Z) \right\} \right]\right\},\\
\mathcal{U}  &\coloneqq \left(\left( \mathrm{E}_{\mathbb{P}_N}[\ell(Z;\beta)]\right)^{\frac{1}{r}} + L_{\beta}^{\mathcal{Z}_N}\delta\right)^{r}.
\end{array}
\end{equation*}
{In particular, \(\left( \mathrm{E}_{\mathbb{P}_N}[\ell(Z;\beta)]\right)^{\frac{1}{r}}\) is well-defined as  \(\psi_{\beta}\) is assumed to be nonnegative {when $r>1$}. } It can be seen that \(\mathcal{S}\) defined in \eqref{eq:defS} satisfies \(\mathcal{S} \leq \mathcal{I} \), whose proof can be found in Appendix~\ref{append:upper-bounds}. Since \(\mathcal{S}\) is a supremum quantity over a feasible set, there exists a sequence of feasible distributions \(\{\mathbb{P}_k\}_{k=1}^{\infty}\) whose expectations \(\{\mathrm{E}_{\mathbb{P}_k} [\ell (Z;\beta) ]\}_{k=1}^{\infty}\) converge to \(\mathcal{S}\) as $k\rightarrow \infty$. More precisely, for any \(\epsilon>0\), we are interested in characterizing \(\tilde{\mathbb{P}}_{\epsilon}\in \mathcal{P}(\mathcal{Z}) \) which satisfies \(\mathcal{W}_{d,r}(\tilde{\mathbb{P}}_{\epsilon} ,\mathbb{P}_N) \leq \delta\) and
$\mathcal{S} - \mathrm{E}_{\tilde{\mathbb{P}}_{\epsilon} } [\ell (Z;\beta) ]   \leq  \epsilon$. In particular, when \(\mathcal{S}\) is attainable, one might even characterize \(\tilde{\mathbb{P}}_0 \) satisfying \(\mathcal{W}_{d,r}(\tilde{\mathbb{P}}_{0} ,\mathbb{P}_N) \leq \delta\) and \(\mathcal{S}  = \mathrm{E}_{\tilde{\mathbb{P}}_{0} } [\ell (Z;\beta) ]   \). Note that, in general, it is not guaranteed that \(\tilde{\mathbb{P}}_0 \) exists, see for example in \citet[Example 2]{mohajerin2018data}. In this work, we call \(\tilde{\mathbb{P}}_{\epsilon} \) as an approximate worst-case distribution \citep{gao2023distributionally}. 

\begin{table}[h]
\fontsize{8pt}{6pt}\selectfont 
\caption{(Informal) comparison of our results with the existing results.}
\label{table:comp-result}
\renewcommand{\arraystretch}{3} 
\begin{tabular}{|c|c|c|c|c|c|l|} 
		\hline 
		\multicolumn{4}{|c|}{Assumptions}  &\multicolumn{2}{c|}{Conclusions}  & \multirow{2}{*}{References} \\  \cline{1-6}
		\(d(\cdot,\cdot)\) & \(\psi_{\beta}\) & \(r\) & Others & \(\mathcal{S},\mathcal{U},\mathcal{I} \) & \(\tilde{\mathbb{P}}_{\epsilon}, \tilde{\mathbb{P}}_{0} \) & \\ \hline  
		norm & max-of-concave & 1& convex \(\mathcal{Z}\subset \mathbb{R}^{n}\)  & \(\mathcal{S}=\mathcal{I}\) & \(\tilde{\mathbb{P}}_{\epsilon}\) & \citet[Thm~4.2, 4.4]{mohajerin2018data}  \\ \hline
		\makecell[c]{extended \\norm} & absolute/logistic & 2 or 1 & \(\mathcal{Z}=\mathbb{R}^{n}\)  & \(\mathcal{S}=\mathcal{U}\) & \(\tilde{\mathbb{P}}_{0} {}^{(\natural)}  \)& \citet[Thm~1,2]{blanchet2019robust}  \\ \hline
		\makecell[c]{extended \\semi-norm} & absolute & 2 & \(\mathcal{Z}=\mathbb{R}^{n}\)  & \(\mathcal{S}=\mathcal{U}\) & \(\tilde{\mathbb{P}}_{0} {}^{(\natural)}\) & \citet[Thm~4]{Chu2021OnRS}  \\ \hline
		metric & globally Lipschitz& 1 & \makecell[c]{tightness\\ at infinity\({}^{(\Join)}\)} & \(\mathcal{S}=\mathcal{U}\) & \(\tilde{\mathbb{P}}_{\epsilon}  \) & \citet[Col~2]{gao2022wasserstein} \\ \hline 		
		\makecell[c]{lower \\ semicontinuous} & \makecell[c]{upper\\ semicontinuous} & \([1,\infty)\) & \makecell[c]{\(d(\cdot,\cdot)\)  positive \\ definite\({}^{(\diamondsuit)} \) }  & \(\mathcal{S}=\mathcal{I}\) & \(\tilde{\mathbb{P}}_0 \) & \citet[Thm~1]{blanchet2019quantifying}  \\ \hline
		cost function & \makecell[c]{interchangeability \\ principle} & \([1,\infty) \) &\makecell[c]{\(d(\cdot,\cdot) \)  positive \\ definite\({}^{(\diamondsuit)} \)}    & \(\mathcal{S}=\mathcal{I}\)  &\makecell[c]{--\({}^{(\flat)}\)} &  \citet[Thm~1]{zhang2022simple} \\ \hline 
		cost function\({}^{(\sharp)}\) & measurable & \([1,\infty) \) & \makecell[c]{\((\mathcal{Z},d)\) locally \\ compact\({}^{(\sharp)}\)} & \(\mathcal{S}=\mathcal{I}\)  &\(\tilde{\mathbb{P}}_{0}\) &  \citet[Thm~1, Col 1]{gao2023distributionally} \\ \hline 
		{cost function}\({}^{(\dagger)}\) & { smooth} &{\([1,\infty) \)}& {\(\mathcal{Z}=\mathbb{R}^{n} \)} & {\(\mathcal{S}=\mathcal{U} \)} & {\makecell[c]{--\({}^{(\flat)}\)} }& \cite[Thm 3.2]{shafieezadeh2023new}  \\\hline
		extended norm & \makecell[c]{convex, \\
			piecewise linear} &  \([1,\infty) \) & \(\mathcal{Z}=\mathbb{R}^{n}\) & \(\mathcal{S}= \mathcal{U}\) & \(\tilde{\mathbb{P}}_{\epsilon}{}^{(\natural)}\)  & \citet[Thm~5,6,7]{wu2022generalization} \\ \hline
		cost function & \((L_{\beta}^{\mathcal{Z}_N},d)\)-Lipschitz &  \([1,\infty) \) & \makecell[c]{tightness\\ conditionally\({}^{(\circledast)}\)} & \(\mathcal{S}= \mathcal{U}\) & \(\tilde{\mathbb{P}}_{\epsilon}{}^{(\natural)}\)  & this work \\ \hline 
\multicolumn{7}{l}{
\footnotesize \makecell[l]{\( {}^{(\sharp)} \) As commented in \citet[Section 1]{blanchet2019quantifying}, the proof given in \citet[Lemma 2]{gao2023distributionally} implicitly \\
assumes that \((\mathcal{Z},d) \) is locally compact. We also notice from \citet[Remark 2, Remark 5]{gao2023distributionally} that \\
\citet[Lemma 2, Corollary 2]{gao2023distributionally} requires \(d(\cdot,\cdot)\) to be a metric.\\
\( {}^{(\natural)} \) The analytical formula of \(\tilde{\mathbb{P}}_{\epsilon}\) (or \(\tilde{\mathbb{P}}_{0}\)) is not given explicitly in the result, but can be constructed directly from the proof.\\
\( {}^{(\flat)} \) The analytical formula of \(\tilde{\mathbb{P}}_{\epsilon}\) (or \(\tilde{\mathbb{P}}_{0}\)) is not a trivial implication from the corresponding result, as far as we understand.\\
\({}^{(\diamondsuit)} \) Positive definiteness/point-separating: \(d(z',z)=0\) \textit{if and only if} \(z'=z\).\\
{\({}^{(\dagger)}\) \(d\) is required to be lower bounded by a metric with compact sublevel set (Assumption 2.1(ii)), hence it must be positive definite.}\\
{\({}^{(\Join)}\) See the discussion after Remark~\ref{remark}.}\\
{\({}^{(\circledast)}\) See Theorem~\ref{thm:main_r1} and Theorem~\ref{thm:main_r2}.}
}} 
\end{tabular} 
\end{table}

In recent years, it has been an emerging topic to study the relationships among $\mathcal{S}$, $\mathcal{I}$ and $\mathcal{U}$. Table~\ref{table:comp-result} shows an informal comparison of our results with the existing  results in the literature; for detailed discussions, see Section~\ref{sec:main-results}. We should note that while \(\mathcal{U}\) is a computationally tractable quantity, \(\mathcal{I}\)  is less computationally friendly as its evaluation involves solving a one-dimensional minimization problem with a complicated objective function. Thus the equivalence \(\mathcal{S}=\mathcal{U}\) is much more desirable than the equivalence \(\mathcal{S}=\mathcal{I}\). Correspondingly, the conditions needed for the equivalence \(\mathcal{S}=\mathcal{U}\) to hold are also stronger. 

We summarize our main contributions in this paper as follows.
\begin{itemize}
	\item We first propose a lower bound \(\mathcal{L}\) for \(\mathcal{S}\) (see Theorem \ref{thm:main}(a)). Then, we characterize a certain property named as the weak \((L_{\beta}^{\mathcal{Z}_N},d)\)-Lipschitz property and prove that under this condition, we have  \(\mathcal{S}\leq \mathcal{U}\) (see Theorem~\ref{thm:main}(b)). The bounds \(\mathcal{L}\) and \(\mathcal{U}\) are demonstrated to exhibit some characteristics that are consistent with the existing literature (see the discussions after Theorem~\ref{thm:main}). 
	\item We propose sufficient conditions for \(\mathcal{S}=\mathcal{U}\) in the cases where \(r=1\) (see Theorem~\ref{thm:main_r1}) and \(r>1 \) (see Theorem~\ref{thm:main_r2}). Our result is a generalization of many existing results in the literature, which is discussed in detail in Section \ref{sec:main-results}. It is worth noting that our proofs do not involve verifying the validity of interchangeability of \(\inf\) and \(\sup\). In particular, we do not use the strong duality result and/or the existence of primal-dual optimizer of the Wasserstein problem \(\mathcal{W}_{d,r}\), which relaxes the assumptions needed for \(d(\cdot,\cdot)\). As a byproduct, our constructive approach directly characterizes the analytic formulation of \(\tilde{\mathbb{P}}_{\epsilon} \).
	\item Although studying sufficient conditions for  \(\mathcal{S}=\mathcal{U}\) is a widely explored topic, we state certain scenarios in which the existing results do not apply. But our results still can cover theses circumstances under suitable conditions,	for instance, even when the globally Lipschitz condition fails (see Example~\ref{exam:non-Lipschitz}); \(\ell\) is not convex (see Example~\ref{exam:remove_delta} and Example~\ref{exam:nonlin-class}); \(d(\cdot,\cdot)\) is not a metric, not positive definite, not convex (see Example~\ref{ex:linear-square}); and \((\mathcal{Z},d)\) is not a locally compact metric space (see Example~\ref{eg:functional-LR}).
	\item We demonstrate the versatility of our theoretical results by applying them to various applications, including regression, classification and risk measure problems (see Section~\ref{sec:app} and Section~\ref{sec: generalrisk}). Table~\ref{table:supper-long-table} shows an informal summary of some applications.
\end{itemize}

{ \fontsize{9pt}{6pt}\selectfont
	\begin{longtblr}
		[caption = {\small (Informal) summary of our applications.},
		label = {table:supper-long-table},
		]{
			colspec = {|l|l|l|},
		} 
		\hline
		Application & Formulation & Remark\\
		\hline
		\SetCell[c=3]{c}{Regression loss functions} \\ \hline
		\makecell[l]{Higher-order\\regression} & 
		$\mathrm{E}_{\mathbb{P}} [\left|Y-\langle \beta,X\rangle\right|^r ]$, $r\geq 1$ 
		& \SetCell[r=3]{l}{ 
			\makecell[l]{
				(1) $(X,Y,\beta)\in \mathbb{R}^n\times \mathbb{R}\times \mathbb{R}^n$ \\[0.5em]
				(2) $d((x',y'),(x,y))$ takes one of \\
				$\boldsymbol{\cdot}$ $\|[x'-x;y'-y]\|_{\mathbb{R}^{n+1}}$\\
				$\boldsymbol{\cdot}$ $\|x'-x\|_{\mathbb{R}^n}+{\boldsymbol \delta}_{\{0\}}(y'-y)$\\
				$\boldsymbol{\cdot}$ ${\boldsymbol \delta}_{\{{\bf 0}_{\mathbb{R}^{|\mathcal{I}^c| +1}}\}}([x'_{\mathcal{I}^c}-x_{\mathcal{I}^c};y'-y])$ \\
				\quad $+\|x'_{\mathcal{I}}-x_{\mathcal{I}}\|_{\mathbb{R}^{|\mathcal{I}|}}$, where \(\mathcal{I}\subset \{1,2,\dots,n\} \)\\
				$\boldsymbol{\cdot}$ $\inf_{\bar{x}\in\mathbb{R}^s}\left\{\|\bar{x}\|_{\mathbb{R}^s} \mid B^T\bar{x}=x'-x \right\}$\\
				\quad $+{\boldsymbol \delta}_{\{0\}}(y'-y)$, where \(B\in\mathbb{R}^{s\times n} \)
			}
		} \\ 
		\cline{1-2}
		\makecell[l]{Lower partial\\moments}  & 
		\makecell[l]{
			$\mathrm{E}_{\mathbb{P}} [ (Y-\langle \beta,X\rangle-\tau)_{+}^r ]$, \\ $r\geq 1$, $\tau\in \mathbb{R}$ }
		& \\ 
		\cline{1-2}
		\makecell[l]{Higher-order \\$\tau$-insensitive \\ regression} & 
		\makecell[l]{
			$\mathrm{E}_{\mathbb{P}} [ (\left|Y-\langle \beta,X\rangle\right|-\tau)_{+}^r ]$, \\ $r\geq 1$, $\tau\in \mathbb{R}$} & \\ 
		\hline
		\makecell[l]{Nonparametric \\
			scalar-on-function \\
			linear regression} & \makecell[l]{$\mathrm{E}_{\mathbb{P}} \left[ h^r( Y - \int\limits_{0}^{1} \left(X(t) \beta(t)) \mathrm{d}t
			\right)\right]$,\\  \(\beta\in\mathfrak{L}^2[0,1]\),  $r\geq 1$} 
		&\SetCell[r=2]{l}{
			\makecell[l]{
				(1) $(X,Y)\in \mathfrak{L}^2[0,1]\times\mathbb{R}$\\[0.5em]
				(2) $d((x',y'),(x,y)) = {\boldsymbol \delta}_{\{0\}}(y'-y)$\\
				\quad \quad $+ \left(\int_{0}^{1}\left|x'(t)-x(t)\right|^2\mathrm{d}t\right)^{1/2}  $\\[0.5em]
				(3) For $s\in \mathbb{R}$ ,$h(s)$ takes one of\\
				$\left|s\right|$, $(s-\tau)_{+}$ or $(\left|s\right|-\tau)_{+}$, $\tau\in \mathbb{R}$
		}}\\
		\cline{1-2}
		\makecell[l]{Parametric \\
			scalar-on-function \\
			linear regression} & \makecell[l]{$\mathrm{E}_{\mathbb{P}} 
			\left[  h^r( Y - \int\limits_{0}^{1} (X(t) \sum\limits_{j=1}^{n} \beta_j {\boldsymbol g}_j(t)) \mathrm{d}t
			)\right]$,\\[0.5em]
			\(\beta\in\mathbb{R}^n\), 
			\(\{{\boldsymbol g}_j\}_{j=1}^n \subset \mathfrak{L}^2[0,1]\), $r\geq 1$} &\\
		\hline
		\makecell[l]{ Log-cosh\\loss regression} & 
		$\mathrm{E}_{\mathbb{P}} [ \log(\cosh(Y-\langle \beta,X\rangle)) ]$ 
		& \SetCell[r=4]{l}{
			\makecell[l]{
				(1) $(X,Y,\beta)\in \mathbb{R}^n\times \mathbb{R}\times \mathbb{R}^n$ \\[0.5em]
				(2) $d((x',y'),(x,y))$ takes one of \\
				$\boldsymbol{\cdot}$ $\|[x'-x;y'-y]\|_{\mathbb{R}^{n+1}}$\\
				$\boldsymbol{\cdot}$ $\|x'-x\|_{\mathbb{R}^n}+{\boldsymbol \delta}_{\{0\}}(y'-y)$
		}} 
		\\
		\cline{1-2}
		\makecell[l]{ Huber\\loss regression} & \makecell[l]{\(\mathrm{E}_{\mathbb{P}} [ h(Y-\langle \beta,X\rangle) ]\) \\ where 
			\(h(t)= \begin{cases}
				t^2/2 & \text{if } \left|t\right| \leq 1,\\
				\left|t\right| -\frac{1}{2} & \text{otherwise}
			\end{cases}\)}  & \\
		\\
		\cline{1-2}
		\makecell[l]{ Quantile\\loss regression} & \makecell[l]{\(\mathrm{E}_{\mathbb{P}} [ h(Y-\langle \beta,X\rangle) ]\) \\ where 
			\(h(t)= \begin{cases}
				\gamma t & \text{if } t \geq 0,\\
				-t & \text{otherwise},
			\end{cases}\) \quad $\gamma\in (0,1)$ }  & \\
		\hline
		
		\makecell[l]{ Ridge linear\\ ordinary regression} &
		$\mathrm{E}_{\mathbb{P}}[(Y + \langle\beta,X\rangle)^2]$ &
		\makecell[l]{
			(1) $(X,Y,\beta)\in \mathbb{R}^n\times \mathbb{R}\times \mathbb{R}^n$, $Z=(X,Y)$ \\[0.5em]
			(2) $d(z',z) = \|z'-z\|_2 \|z'+z\|_2$
		}\\
		\hline
		\makecell[l]{ Hard sigmoid \\ /HardTanh } &
		$\mathrm{E}_{\mathbb{P}} \left[\max \left\{ 0, \min \left\{1, \frac{\langle \beta,Z\rangle +1}{2} \right\} \right\} \right]$ & 
		\makecell[l]{ (1) $(Z,\beta)\in \mathbb{R}^n\times \mathbb{R}^n$\\
			(2) \(d(z',z) = \|z'-z\|_{\mathbb{R}^n} \)\\
			(3) Equivalence holds conditionally
		}
		\\
		\hline %\pagebreak
		\SetCell[c=3]{c}{Classification loss functions} \\ \hline
		\makecell[l]{Higher-order\\hinge loss\\binary classification} & 
		$\mathrm{E}_{\mathbb{P}} [\left( 1-Y\cdot\langle\beta,X\rangle \right)_{+}^r ]$, $r\geq 1$ 
		& \SetCell[r=5]{l}{
			\makecell[l]{
				(1) $(X,Y,\beta)\in \mathbb{R}^n\times \{-1,1\}\times \mathbb{R}^n$ \\[0.5em]
				(2) $d((x',y'),(x,y))=\|x'-x\|_{\mathbb{R}^n}$\\
				\quad \quad  $+ {\boldsymbol \delta}_{\{0\}}(y'-y)$
		}} \\ 
		\cline{1-2}
		\makecell[l]{Higher-order\\support vector\\machine classification} & 
		$\mathrm{E}_{\mathbb{P}} [\left| 1-Y\cdot\langle\beta,X\rangle\right|^r ]$, $r\geq 1$ 
		& \\
		\cline{1-2}
		Log-exponential loss & $\mathrm{E}_{\mathbb{P}} [
		\log (1+\exp(-Y\cdot \langle\beta,X\rangle))]$ & \\
		\cline{1-2}
		Smooth hinge loss & \makecell[l]{ $\mathrm{E}_{\mathbb{P}} [h(Y\cdot\langle \beta,X \rangle )]$ with \\
			\(h(t) =\left\{\begin{array}{ll}
				0 & \text{if } t \geq 1,\\
				(1-t)^2/2 & \text{if } 0<t<1,\\
				1/2-t & \text{otherwise}
			\end{array}\right.\)}
		& \\
		\cline{1-2}
		\makecell[l]{Truncated \\
			pinball loss }
		& \makecell[l]{$\mathrm{E}_{\mathbb{P}} [h(Y\cdot\langle \beta,X \rangle)]$ with \\ \(h(t)= \begin{cases}
				1-t & \text{if } t \leq 1,\\
				\tau_1 (t-1) & \text{if } 1 < t < \tau_2+1,\\
				\tau_1 \tau_2 & \text{otherwise},
			\end{cases} \) \\
			\(\tau_1 \in [0,1],\tau_2 \geq 0 \)
		}& \\
		\hline
		\makecell[l]{ Binary \\ cross-entropy loss } &
		$\mathrm{E}_{\mathbb{P}} [\beta Z\log(\beta Z) + (1-\beta Z)\log(1-\beta Z)]$ & 
		\makecell[l]{ (1) $(Z,\beta)\in (0,1)\times (0,1)$\\
			(2) \(d(z',z) = \left|z'-z\right| \)\\
			(3) Equivalence holds conditionally
		}
		\\
		\hline
		\SetCell[c=3]{c}{Generalization to risk measure} \\ \hline
		\makecell[l]{$\nu$-support \\ vector \\regression} & 
		$\mathrm{CVaR}_{\alpha}^{\mathbb{P}}(\left|Y - \langle \beta,X\rangle \right|)$
		& \makecell[l]{
			(1) $(X,Y,\beta)\in \mathbb{R}^n\times\mathbb{R}\times \mathbb{R}^n$\\
			(2) $d((x',y'),(x,y))$\\
			\quad \quad $= \|(x',y')-(x,y)\|_{\mathbb{R}^{n+1}}$\\
			(3) $\alpha\in (0,1)$} \\ 
		\hline
		\makecell[l]{$\nu$-support \\ vector \\
			machine} & 
		$\mathrm{CVaR}_{\alpha}^{\mathbb{P}}(-Y\cdot\langle\beta,X\rangle)$
		& \makecell[l]{
			(1) $(X,Y,\beta)\in \mathbb{R}^n\times\mathbb{R}\times \mathbb{R}^n$\\
			(2) $d((x',y'),(x,y))=\|x'-x\|_{\mathbb{R}^n}$\\
			\quad \quad $+ {\boldsymbol \delta}_{\{0\}}(y'-y)$\\
			(3) $\alpha\in (0,1)$} \\ 
		\hline
		\makecell[l]{Higher \\
			moment\\ coherent \\risk measures} & 
		\makecell[l]{ $\inf\limits_{t\in \mathbb{R}} \Big\{ t+\frac{1}{1-\alpha}
			\left(\mathrm{E}_{\mathbb{P}} [\left( \langle\beta,Z\rangle -t\right)_{+}^r ]\right)^{\frac{1}{r}}\Big\}$\\
			$r\geq 1$}
		& \makecell[l]{
			(1) $(Z,\beta)\in \mathbb{R}^n\times \mathbb{R}^n$\\
			(2) $d(z',z)=\|z'-z\|_{\mathbb{R}^n}$\\
			(3) $\alpha\in (0,1)$} \\ 
		\hline
	\end{longtblr}
}

\section{Theoretical analysis of the equivalence} \label{sec:main-results}
In this section, we will establish our main results for deriving the equivalence between the worst-case loss quantity in the WDRO problem and its regularization scheme counterpart. We first give the following lemma, with the proof in Appendix \ref{append:proof_lemma1}, to quantify the Wasserstein discrepancy between any distribution and a singleton, which will be used in the subsequent analysis.

\begin{lemma}\label{lemma: singleton}
	Given any distribution \(\mathbb{P}\in\mathcal{P}(\mathcal{Z}) \) and any point $\hat{z}\in \mathcal{Z}$, for any scalar $r\geq 1$ and any extended nonnegative-valued measurable function \(d\colon\mathcal{Z}\times\mathcal{Z} \rightarrow [0,\infty] \), we have
	\begin{equation*}
		\mathcal{W}_{d,r}(\mathbb{P},{\boldsymbol \chi}_{\{\hat{z}\}}) =  \left(\int_{\mathcal{Z}}d^{r}(z,\hat{z})\mathrm{d}\mathbb{P}(z)\right)^{\frac{1}{r}}. 
	\end{equation*}
\end{lemma}

Before stating our main results, we need the following definition of the cost function. 
\begin{definition}
	The function $d(\cdot,\cdot)$ defined on $\mathcal{Z}\times\mathcal{Z}$ is called a cost function if it is extended nonnegative-valued, measurable, and vanishes whenever two arguments are the same, that is, for any $z',z\in \mathcal{Z}$, $d(z',z)\in [0,\infty]$ and \(d(z,z)=0\). 
\end{definition}

Next, we introduce a weak Lipschitz property for functions on $\mathcal{Z}$ with respect to a given cost function $d(\cdot,\cdot)$, where the weak Lipschitz constant depends on the second input of $d(\cdot,\cdot)$. Note that it is different from the Lipschitz property used in \citet{shafieezadeh2015distributionally,an2021generalization,gao2022finite}, as the latter does not depend on the input variables. 

\begin{definition}[Weak Lipschitz property] \label{def:Lip}
	Given a function \(f\colon\mathcal{Z}\rightarrow\mathbb{R}\), a cost function \(d(\cdot,\cdot)\) on \(\mathcal{Z}\times\mathcal{Z}\) and a subset \(\mathcal{S}\subset\mathcal{Z} \), \(f\) is called \((L_{f}^{\mathcal{S}},d)\)-Lipschitz at \(\mathcal{S}\) if for any \(z\in\mathcal{S}, z'\in\mathcal{Z}\), one has 
	\begin{equation*} 
		\left|f(z') -f(z)\right| \leq L_{f}^{\mathcal{S}}d(z',z),
	\end{equation*} 
	where $L_{f}^{\mathcal{S}}\in [0,\infty)$ is a constant depending on $f$ and $\mathcal{S}$. 
\end{definition}

The classical Lipschitz property can be seen as a special case of Definition~\ref{def:Lip} when \((\mathcal{Z},d) \) is a metric space and \(\mathcal{S}=\mathcal{Z}\). In other words, any Lipschitz function is weak Lipschitz, while the reverse is not always true, for example, see Example~\ref{exam:non-Lipschitz}.  

\subsection{Lower and Upper bounds of the worst-case loss quantity}
We now derive a lower bound and an upper bound of the worst-case loss quantity for a certain class of loss functions in the following theorem. 
\begin{theorem}\label{thm:main}
	Let \(\mathcal{Z}_N\coloneqq\{Z^{(1)},\dots,Z^{(N)}\}\subset\mathcal{Z} \) be a given dataset and \(\mathbb{P}_N\coloneqq\sum_{i=1}^{N}\mu_i {\boldsymbol \chi}_{\{Z^{(i)}\}}\in\mathcal{P}(\mathcal{Z})\) be the corresponding empirical distribution. In addition, let $r\in [1,\infty)$ be a scalar and \(d(\cdot,\cdot)\) be a cost function on \(\mathcal{Z}\times\mathcal{Z}\). Suppose the loss function $\ell:\mathcal{Z}\times \mathcal{B}\rightarrow \mathbb{R}$ takes the form as
	\begin{equation*}
		\ell\colon (z;\beta) \mapsto \psi_{\beta}^{r}(z), \text{ with } \begin{cases}
			\psi_{\beta}\colon\mathcal{Z}\rightarrow \mathbb{R} \ \ \text{ if } r =  1, \\
			\psi_{\beta}\colon\mathcal{Z}\rightarrow \mathbb{R}_+ \text{ if } r > 1.
		\end{cases}
	\end{equation*}
	Let \(\mathcal{S}\) be defined as in \eqref{eq:defS}. Then the following statements hold for any \(\delta\geq 0\). 	
	\begin{enumerate}[label=(\alph*)]		
		\item Let \(\mathcal{L}_i \coloneqq  \sup_{\mathbb{P}\in \mathcal{P}(\mathcal{Z})} \left\{\mathrm{E}_{\mathbb{P}}[\ell(Z;\beta)] \middle\vert  \mathcal{W}_{d,r}(\mathbb{P},{\boldsymbol \chi}_{\{Z^{(i)}\}}) \leq \delta \right\}\) for \(i=1,\dots,N\). Then
		\begin{equation*}
			\mathcal{S} \geq \mathcal{L} \coloneqq  \sum_{i=1}^{N}\mu_i \mathcal{L}_i \geq \mathrm{E}_{\mathbb{P}_N}[\ell(Z;\beta)] .
		\end{equation*}
		
		\item Suppose \(\psi_{\beta}\) is \((L_{\beta}^{\mathcal{Z}_N},d)\)-Lipschitz at \(\mathcal{Z}_N\) with \(L_{\beta}^{\mathcal{Z}_N}\in (0,\infty)\), then 
		\begin{equation*}
			\mathcal{S}\leq  \mathcal{U}=  \left(\left(\mathrm{E}_{\mathbb{P}_N}[\ell(Z;\beta)] \right)^{\frac{1}{r}}+ L_{\beta}^{\mathcal{Z}_N}\delta\right)^{r}.
		\end{equation*}
		\item Suppose \(\psi_{\beta}\) is \((0,d)\)-Lipschitz at \(\mathcal{Z}_N\), then 
		\begin{equation*}
			\mathcal{S} = \mathrm{E}_{\mathbb{P}_N}[\ell(Z;\beta)].
		\end{equation*}
	\end{enumerate}	
\end{theorem}
\textbf{Proof.}
(a) For any collection \(\left\{\tilde{\mathbb{P}}^{(i)}\right\}_{i=1}^N \subseteq \mathcal{P}(\mathcal{Z}) \) such that \(\mathcal{W}_{d,r}(\tilde{\mathbb{P}}^{(i)},{\boldsymbol \chi}_{\{Z^{(i)}\}})\leq\delta\) for all $i=1,\cdots,N$. It follows from Lemma \ref{lemma: singleton} that for any \(i = 1,\cdots,N,\)
\begin{equation*}
	\mathcal{W}_{d,r}(\tilde{\mathbb{P}}^{(i)},{\boldsymbol \chi}_{\{Z^{(i)}\}})= \left(\int_{\mathcal{Z}} d^{r}(z,Z^{(i)})\mathrm{d}\tilde{\mathbb{P}}^{(i)}(z)\right)^{\frac{1}{r}}  \leq \delta. 
\end{equation*}
Hence we can construct
\begin{align} 
	\tilde{\mathbb{P}}\coloneqq \sum_{i=1}^{N}\mu_i\tilde{\mathbb{P}}^{(i)} \quad \text{and} \quad \tilde{\pi}\coloneqq \sum_{i=1}^{N}\left(\mu_i \tilde{\mathbb{P}}^{(i)}\otimes{\boldsymbol \chi}_{\{Z^{(i)}\}}\right).\label{eq:construct_P}
\end{align}  
Then we have \(\tilde{\mathbb{P}}\in\mathcal{P}(\mathcal{Z}),\tilde{\pi}\in\mathcal{P}(\mathcal{Z}\times\mathcal{Z})  \), and \(\tilde{\pi}\in\Pi(\tilde{\mathbb{P}},\mathbb{P}_N) \), since for any measurable sets \(A,B\subset\mathcal{Z}\),
\begin{equation*} 
	\begin{array}{llll}
		\tilde{\pi}(\mathcal{Z}\times B) &= \sum_{i=1}^{N}  \mu_i  \tilde{\mathbb{P}}^{(i)} (\mathcal{Z}){\boldsymbol \chi}_{\{Z^{(i)}\}}(B)   
		&= \sum_{i=1}^{N}\mu_i{\boldsymbol \chi}_{\{Z^{(i)}\}}(B) = \mathbb{P}_N (B),\\
		\tilde{\pi}(A\times\mathcal{Z})  &= \sum_{i=1}^{N}  \mu_i \tilde{\mathbb{P}}^{(i)}(A)  {\boldsymbol \chi}_{\{Z^{(i)}\}}(\mathcal{Z})
		&= \sum_{i=1}^{N}\mu_i\tilde{\mathbb{P}}^{(i)}(A)  = \tilde{\mathbb{P}}(A). 
	\end{array}
\end{equation*} 
Therefore, we can see that
\begin{equation*}
	\begin{array}{ll}
		\mathcal{W}_{d,r}(\tilde{\mathbb{P}},\mathbb{P}_N) &\leq \left(\int_{\mathcal{Z}\times\mathcal{Z}} d^{r}(\tilde{z},z)\mathrm{d}\tilde{\pi}(\tilde{z},z)\right)^{\frac{1}{r}}
	 = \left(\sum_{i=1}^{N}\mu_i \int_{\mathcal{Z}}{d}^{r}(\tilde{z},Z^{(i)})\mathrm{d}\tilde{\mathbb{P}}^{(i)}(\tilde{z}) \right)^{\frac{1}{r}} \\
	&=\left(\sum_{i=1}^{N}\mu_i 
	\left( \mathcal{W}_{d,r}(\tilde{\mathbb{P}}^{(i)},{\boldsymbol \chi}_{\{Z^{(i)}\}})\right)^r \right)^{\frac{1}{r}} \leq \delta.
	\end{array}
\end{equation*}
Moreover, according to \eqref{eq:construct_P}, we have 
\begin{equation*}
	\mathrm{E}_{\tilde{\mathbb{P}}} [\ell (Z;\beta) ] = \sum_{i=1}^N \mu_i \mathrm{E}_{\tilde{\mathbb{P}}^{(i)}}[\ell(Z;\beta)].
\end{equation*}
By taking supremum on all possible \(\left\{\tilde{\mathbb{P}}^{(i)}\right\}_{i=1}^N \) such that $\mathcal{W}_{d,r}(\tilde{\mathbb{P}}^{(i)},{\boldsymbol \chi}_{\{Z^{(i)}\}}) \leq \delta$ for all $i=1,\cdots,N$, we have
\begin{equation*}
	\begin{array}{ll}
		\mathcal{S} = \sup_{\mathbb{P}\colon  \mathcal{W}_{d,r}(\mathbb{P},\mathbb{P}_N) 
			\leq \delta  }  \mathrm{E}_{\mathbb{P}} [\ell (Z;\beta) ]
			\geq \sum\limits_{i=1}^{N}\mu_i   \sup\limits_{\mathbb{P}\in\mathcal{P}(\mathcal{Z})} \left\{\mathrm{E}_{\mathbb{P}}[\ell(Z;\beta)] \middle\vert \mathcal{W}_{d,r}\left(\mathbb{P},{\boldsymbol \chi}_{\{Z^{(i)}\}}\right) \leq \delta \right\}
			=\sum_{i=1}^{N}\mu_i\mathcal{L}_i = \mathcal{L}.
	\end{array}
\end{equation*}
Besides, since \(\mathcal{W}_{d,r}({\boldsymbol \chi}_{\{Z^{(i)}\}},{\boldsymbol \chi}_{\{Z^{(i)}\}}) = 0 \leq \delta \) by Lemma~\ref{lemma: singleton}, we have that 
\begin{equation*}
	\begin{array}{ll}
		\mathcal{L}_i = \sup\limits_{\mathbb{P}\in\mathcal{P}(\mathcal{Z})} \left\{\mathrm{E}_{\mathbb{P}}[\ell(Z;\beta)] \middle\vert \mathcal{W}_{d,r}\left(\mathbb{P},{\boldsymbol \chi}_{\{Z^{(i)}\}}\right) \leq \delta \right\}
		\geq \mathrm{E}_{{\boldsymbol \chi}_{\{Z^{(i)}\}}}[\ell(Z;\beta)] =  \ell(Z^{(i)};\beta),
	\end{array}
\end{equation*}
and hence \(\mathcal{L}=\sum_{i=1}^{N}\mu_i\mathcal{L}_i \geq \sum_{i=1}^{N}\mu_i \ell(Z^{(i)};\beta) = \mathrm{E}_{\mathbb{P}_N}[\ell(Z;\beta)] \).

(b) Let \(\epsilon>0\) be an arbitrary scalar. Fix any \(  \tilde{\mathbb{P}} \in \mathcal{P}(\mathcal{Z})  \) such that \(\mathcal{W}_{d,r}(\tilde{\mathbb{P}},\mathbb{P}_N) \leq \delta \). By the definition of $\mathcal{W}_{d,r}(\cdot,\cdot)$, there exists \(\tilde{\pi}\in\Pi(\tilde{\mathbb{P}},\mathbb{P}_N ) \)  such that 
\begin{equation*}
	\left(\int_{\mathcal{Z}\times\mathcal{Z}} d^{r}(\tilde{z},z) \mathrm{d}\tilde{\pi}(\tilde{z},z)\right)^{\frac{1}{r}} \leq \delta + \frac{\epsilon}{L_{\beta}^{\mathcal{Z}_N}}.
\end{equation*}
Besides, by the definition of the loss function $\ell(\cdot,\cdot)$, we have 
\begin{equation*}
	\begin{array}{ll}
		\left(\mathrm{E}_{\tilde{\mathbb{P}}} [\ell(Z;\beta)]\right)^{\frac{1}{r}} 
		&= \left(\int_{\mathcal{Z}} \psi_{\beta}^{r}(\tilde{z})\mathrm{d}\tilde{\mathbb{P}}(\tilde{z})\right)^{\frac{1}{r}} =\left(\int_{\mathcal{Z}\times\mathcal{Z}} \psi_{\beta}^{r}(\tilde{z})\mathrm{d}\tilde{\pi}(\tilde{z},z)\right)^{\frac{1}{r}}  \\
		&= \left(\int_{\mathcal{Z}\times\mathcal{Z}} \left(\psi_{\beta}(z)+ \psi_{\beta}(\tilde{z})-\psi_{\beta}(z)\right)^{r}\mathrm{d}\tilde{\pi}(\tilde{z},z)\right)^{\frac{1}{r}} \\
		 &\leq^{(*)} \left(\int_{\mathcal{Z}\times\mathcal{Z}} \psi_{\beta}^r(z)\mathrm{d}\tilde{\pi}(\tilde{z},z)\right)^{\frac{1}{r}}
		  +\left(\int_{\mathcal{Z}\times\mathcal{Z}} \left|\psi_{\beta}(\tilde{z})-\psi_{\beta}(z)\right|^{r}\mathrm{d}\tilde{\pi}(\tilde{z},z)\right)^{\frac{1}{r}}\\
		&= \left(\int_{\mathcal{Z}} \psi_{\beta}^r(z)\mathrm{d}\mathbb{P}_N(z)\right)^{\frac{1}{r}}
		+\left(\int_{\mathcal{Z}\times\mathcal{Z}} \left|\psi_{\beta}(\tilde{z})-\psi_{\beta}(z)\right|^{r}\mathrm{d}\tilde{\pi}(\tilde{z},z)\right)^{\frac{1}{r}}\\
		&=\left(\mathrm{E}_{\mathbb{P}_N}[\ell(Z;\beta)] \right)^{\frac{1}{r}}
		 +\left(\int_{\mathcal{Z}\times\mathcal{Z}} \left|\psi_{\beta}(\tilde{z})-\psi_{\beta}(z)\right|^{r}\mathrm{d}\tilde{\pi}(\tilde{z},z)\right)^{\frac{1}{r}},
	\end{array}
\end{equation*} 
where the inequality \({}^{(*)}\) holds naturally if \(r=1\), and follows from the Minkowski inequality if \(r>1\). Since \(\psi_{\beta}\) is \((L_{\beta}^{\mathcal{Z}_N},d)\)-Lipschitz at \(\mathcal{Z}_N\), it holds that
\begin{equation*}
	\begin{array}{ll}
		\left(\mathrm{E}_{\tilde{\mathbb{P}}} [\ell(Z;\beta)]\right)^{\frac{1}{r}}
		 \leq \left(\mathrm{E}_{\mathbb{P}_N}[\ell(Z;\beta)] \right)^{\frac{1}{r}} + L_{\beta}^{\mathcal{Z}_N} \left(\int_{\mathcal{Z}\times\mathcal{Z}} d^{r}(\tilde{z},z)\mathrm{d}\tilde{\pi}(\tilde{z},z)\right)^{\frac{1}{r}}
		 \leq \left(\mathrm{E}_{\mathbb{P}_N}[\ell(Z;\beta)] \right)^{\frac{1}{r}}+ L_{\beta}^{\mathcal{Z}_N}\delta + \epsilon.
	\end{array}
\end{equation*}
This means that for any \(\epsilon>0\), we have
\begin{equation*}
	\begin{array}{ll}
		\mathcal{S}^{\frac{1}{r}}=\sup_{\mathbb{P}\colon\mathcal{W}_{d,r}(\mathbb{P},\mathbb{P}_N) \leq \delta}\left(\mathrm{E}_{\mathbb{P}} [\ell (Z;\beta) ]\right)^{\frac{1}{r}} 
		\leq  \left(\mathrm{E}_{\mathbb{P}_N}[\ell(Z;\beta)] \right)^{\frac{1}{r}}+ L_{\beta}^{\mathcal{Z}_N}\delta + \epsilon.
	\end{array}
\end{equation*}
By letting \(\epsilon\rightarrow0\), we get the desired inequality.

(c) Since \(\psi_{\beta}\) is \((0,d)\)-Lipschitz at \(\mathcal{Z}_N\), by the convention that \( 0\cdot\infty=0\), one has \(\psi_{\beta}(z') = \psi_{\beta}(z)\) for any \(z'\in\mathcal{Z},z\in\mathcal{Z}_N\). In particular, \(\psi_{\beta}(\bar{z}) =\psi_{\beta}(z)  \) for any \(\bar{z},z\in\mathcal{Z}_N \). Therefore, \(\psi_{\beta}(\cdot) \) is a constant function on $\mathcal{Z}$, and so is \(\ell(\cdot;\beta)\). Thus, we have
\begin{equation*}
	\mathcal{S}=\sup_{\mathbb{P}\colon  \mathcal{W}_{d,r}(\mathbb{P},\mathbb{P}_N) \leq \delta }  \mathrm{E}_{\mathbb{P}} [\ell (Z;\beta) ] = \mathrm{E}_{\mathbb{P}_N}[\ell(Z;\beta)].
\end{equation*}
This completes the proof.
\hfill\(\square\) 

For better understanding, we give an intuitive explanation of the above theorem. The first conclusion shows that the worst-case loss quantity \(\mathcal{S}= \sup_{\mathbb{P}\colon  \mathcal{W}_{d,r}(\mathbb{P},\mathbb{P}_N) \leq \delta }  \mathrm{E}_{\mathbb{P}} [\ell (Z;\beta) ]\) is lower bounded by the weighted average of \(N\) worst-case loss quantities with respect to \(N\) point masses \({\boldsymbol \chi}_{\{Z^{(i)}\}}\), $i=1,\cdots,N$, which is easier to calculate according to Lemma \ref{lemma: singleton}. The second conclusion shows that the worst-case loss quantity can be upper bounded by using the weak Lipschitz property of the kernel function \(\psi_{\beta}\). The third conclusion shows that when the weak Lipschitz constant is zero, then the upper bound is met. In addition, from (a) and (b), we have that if \(\delta=0\), then \(\mathcal{S}\geq\mathcal{L}\geq\mathrm{E}_{\mathbb{P}_N}[\ell(Z;\beta)] = \mathcal{U} \geq \mathcal{S}  \). That is to say, if \(\psi_{\beta}\) is \((L_{\beta}^{\mathcal{Z}_N},d)\)-Lipschitz at \(\mathcal{Z}_N\) and \(\delta=0\), then \(\mathcal{S}=\mathcal{U}=\mathcal{L}=\mathrm{E}_{\mathbb{P}_N}[\ell(Z;\beta)] \). Therefore, in the remaining part of this work, we shall only consider the case when \(L_{\beta}^{\mathcal{Z}_N}\in (0,\infty)\) and \(\delta\in(0,\infty)\).

Note that if one fixes the input data while varying \(\delta\), the worst-case loss quantity \(\mathcal{S}(\cdot)\), the lower bound \(\mathcal{L}(\cdot)\) and the upper bound \(\mathcal{U}(\cdot)\) are all functions of $\delta$ on \([0,\infty) \). A few remarks are in order. First, the lower bound \(\mathcal{L}(\cdot)\) is larger than the trivial lower bound \(\mathrm{E}_{\mathbb{P}_N}[\ell(Z;\beta)]\) given in \citet[Lemma 1]{zhang2022simple}. Second, the upper bound \(\mathcal{U}(\cdot) \) is continuous on \([0,\infty) \) and in particular, right-continuous at \(\delta=0\). This is similar to the continuity of another existing upper bound \(\mathcal{I}(\cdot) \) \cite[Remark 2]{zhang2022simple}. Third, if we have \(\mathcal{S}(\delta) = \mathcal{U}(\delta) \) for each \(\delta\in(0,\infty) \) (for example, when Theorem~\ref{thm:main_r1} or Theorem~\ref{thm:main_r2} holds true), then \(\mathcal{S}(\cdot) \) is continuous on \([0,\infty) \). It is worth noting that another sufficient condition for the continuity of \(\mathcal{S}(\cdot)\) has been studied in  \citet[Remark 2]{zhang2022simple}, where the loss function \(\ell\) is required to be a composition of a non-decreasing concave function and the cost function \(d(\cdot,\cdot)\).

Next, we will analyse various cases on when the lower or upper bound for the worst-case loss quantity provided in Theorem~\ref{thm:main} is achievable. 

\subsection{Equivalence in \eqref{eq:main_equation} when $r=1$}
We first consider the case when \(r=1\). The following theorem provides a sufficient condition on when the lower and upper bounds in Theorem~\ref{thm:main} will coincide.  
\begin{theorem} \label{thm:main_r1}
	Let \(\mathcal{Z}_N\coloneqq\{Z^{(1)},\dots,Z^{(N)}\}\subset\mathcal{Z} \) be a given dataset and \(\mathbb{P}_N\coloneqq\sum_{i=1}^{N}\mu_i{\boldsymbol \chi}_{\{Z^{(i)}\}}\in\mathcal{P}(\mathcal{Z})\) be the corresponding empirical distribution. In addition, let \(d(\cdot,\cdot)\) be a cost function on \(\mathcal{Z}\times\mathcal{Z}\) and  \(\delta\in(0,\infty)\) be a scalar. Suppose the loss function $\ell:\mathcal{Z}\times \mathcal{B}\rightarrow \mathbb{R}$ takes the form as
	\begin{equation*}
		\ell\colon (z;\beta) \mapsto \psi_{\beta}(z), 
	\end{equation*}
	where the function $\psi_{\beta}\colon\mathcal{Z}\rightarrow \mathbb{R}$ satisfies the following assumptions:
	\begin{enumerate}[label=(A\arabic*)]
		\item \(\psi_{\beta}\) is \((L_{\beta}^{\mathcal{Z}_N},d)\)-Lipschitz at \(\mathcal{Z}_N\) with \(L_{\beta}^{\mathcal{Z}_N}\in(0,\infty)\);
		\item for any \(\epsilon\in (0,L_{\beta}^{\mathcal{Z}_N}) \) and each \(Z^{(i)}\in\mathcal{Z}_N\), there exists \(\tilde{Z}^{(i)}_{\epsilon}\in\mathcal{Z}\) such that \(\delta\leq d(\tilde{Z}^{(i)}_{\epsilon},Z^{(i)})<\infty\) and
		\begin{equation*}
			\psi_{\beta}(\tilde{Z}^{(i)}_{\epsilon}) - \psi_{\beta}(Z^{(i)}) \geq   (L_{\beta}^{\mathcal{Z}_N}-\epsilon) d(\tilde{Z}^{(i)}_{\epsilon},Z^{(i)}).
		\end{equation*}
	\end{enumerate}
	Then we have that \(\mathcal{L}=\mathcal{S}=\mathcal{U}\) in Theorem~\ref{thm:main}, that is, 
	\begin{align} 
		\sup_{\mathbb{P}\colon  \mathcal{W}_{d,1}(\mathbb{P},\mathbb{P}_N) \leq \delta }  \mathrm{E}_{\mathbb{P}} [\ell (Z;\beta) ] = \mathrm{E}_{\mathbb{P}_N}[\ell(Z;\beta)] + L_{\beta}^{\mathcal{Z}_N}\delta.\label{eq: r1}
	\end{align}
\end{theorem}
\textbf{Proof.}
Since \(\psi_{\beta}\) is \((L_{\beta}^{\mathcal{Z}_N},d) \)-Lipschitz at \(\mathcal{Z}_N\), by Theorem~\ref{thm:main}, we have that  
\begin{equation*}
	\begin{array}{ll}
		\mathcal{L}\leq\mathcal{S}= \sup_{\mathbb{P}\colon  \mathcal{W}_{d,1}(\mathbb{P},\mathbb{P}_N) \leq \delta }  \mathrm{E}_{\mathbb{P}} [\ell (Z;\beta) ]\\
		\leq \mathcal{U}= \mathrm{E}_{\mathbb{P}_N}[\ell(Z;\beta)] + L_{\beta}^{\mathcal{Z}_N}\delta.
	\end{array}
\end{equation*}
Hence, in order to prove \eqref{eq: r1}, it suffices to show that \(\mathcal{L}\geq\mathcal{U}\).

Let \(\epsilon \in \left(0,\min \{L_{\beta}^{\mathcal{Z}_N}, \delta L_{\beta}^{\mathcal{Z}_N} \}  \right) \) be an arbitrary scalar. By Assumption (A2), for any \(Z^{(i)}\in\mathcal{Z}_N\), there exists \(\tilde{Z}^{(i)}\in\mathcal{Z}\) such that $\delta\leq d(\tilde{Z}^{(i)},Z^{(i)})<\infty$ and
\begin{equation*}
	\psi_{\beta}(\tilde{Z}^{(i)}) - \psi_{\beta}(Z^{(i)}) \geq \left(L_{\beta}^{\mathcal{Z}_N}-\frac{\epsilon}{\delta}\right)d(\tilde{Z}^{(i)},Z^{(i)}).
\end{equation*}
Let \(\eta^{(i)} \coloneqq \delta/d(\tilde{Z}^{(i)},Z^{(i)})\in(0,1] \) and choose 
\begin{equation*}
	\tilde{\mathbb{P}}^{(i)} \coloneqq \eta^{(i)}{\boldsymbol \chi}_{\{\tilde{Z}^{(i)}\}} + (1-\eta^{(i)}) {\boldsymbol \chi}_{\{Z^{(i)}\}}\in\mathcal{P}(\mathcal{Z}).
\end{equation*}
Then one has 
\begin{equation*}
	\begin{array}{ll}
		\mathcal{W}_{d,1}\left(\tilde{\mathbb{P}}^{(i)},{\boldsymbol \chi}_{\{Z^{(i)}\}}\right)
		  =\eta^{(i)} d(\tilde{Z}^{(i)},Z^{(i)}) + (1-\eta^{(i)})d(Z^{(i)},Z^{(i)})
		    = \eta^{(i)}d(\tilde{Z}^{(i)},Z^{(i)}) = \delta,
	\end{array}
\end{equation*}
and
\begin{equation*}
	\begin{array}{ll}
		\mathrm{E}_{\tilde{\mathbb{P}}^{(i)}}[\ell(Z;\beta)] &=  \eta^{(i)} \psi_{\beta}(\tilde{Z}^{(i)}) + (1-\eta^{(i)} )\psi_{\beta}(Z^{(i)})
		= \psi_{\beta}(Z^{(i)}) + \eta^{(i)} \left[\psi_{\beta}(\tilde{Z}^{(i)}) -\psi_{\beta}(Z^{(i)})\right]\\
		&\geq \psi_{\beta}(Z^{(i)}) + \eta^{(i)} \left(L_{\beta}^{\mathcal{Z}_N}-\frac{\epsilon}{\delta}\right)d(\tilde{Z}^{(i)},Z^{(i)}) 
		= \ell(Z^{(i)};\beta)  + L_{\beta}^{\mathcal{Z}_N}\delta - \epsilon.
	\end{array}
\end{equation*}
Let $\epsilon\rightarrow 0$, we have that for any $i=1,\cdots,N$, 
\begin{equation*}
	\begin{array}{ll}
		\mathcal{L}_i=\sup\limits_{\mathbb{P}\in \mathcal{P}(\mathcal{Z}) } \left\{\mathrm{E}_{\mathbb{P}}[\ell(Z;\beta)] \middle\vert  \mathcal{W}_{d,1}\left(\mathbb{P},{\boldsymbol \chi}_{\{Z^{(i)}\}}\right) \leq \delta \right\}
		 \geq  \ell(Z^{(i)};\beta) + L_{\beta}^{\mathcal{Z}_N}\delta.
	\end{array}
\end{equation*} 
Therefore, it holds that
\begin{equation*}
	\begin{array}{ll}
		\mathcal{L} = \sum_{i=1}^N \mu_i \mathcal{L}_i \geq \sum_{i=1}^N \mu_i \left( \ell(Z^{(i)};\beta) + L_{\beta}^{\mathcal{Z}_N}\delta\right) 
		=\mathrm{E}_{\mathbb{P}_N}[\ell(Z;\beta)] + L_{\beta}^{\mathcal{Z}_N}\delta = \mathcal{U}.
	\end{array}
\end{equation*}
This completes the proof.
\hfill\(\square\) 

\begin{figure}[H]
	\centering
	\includegraphics[width=0.7\textwidth]{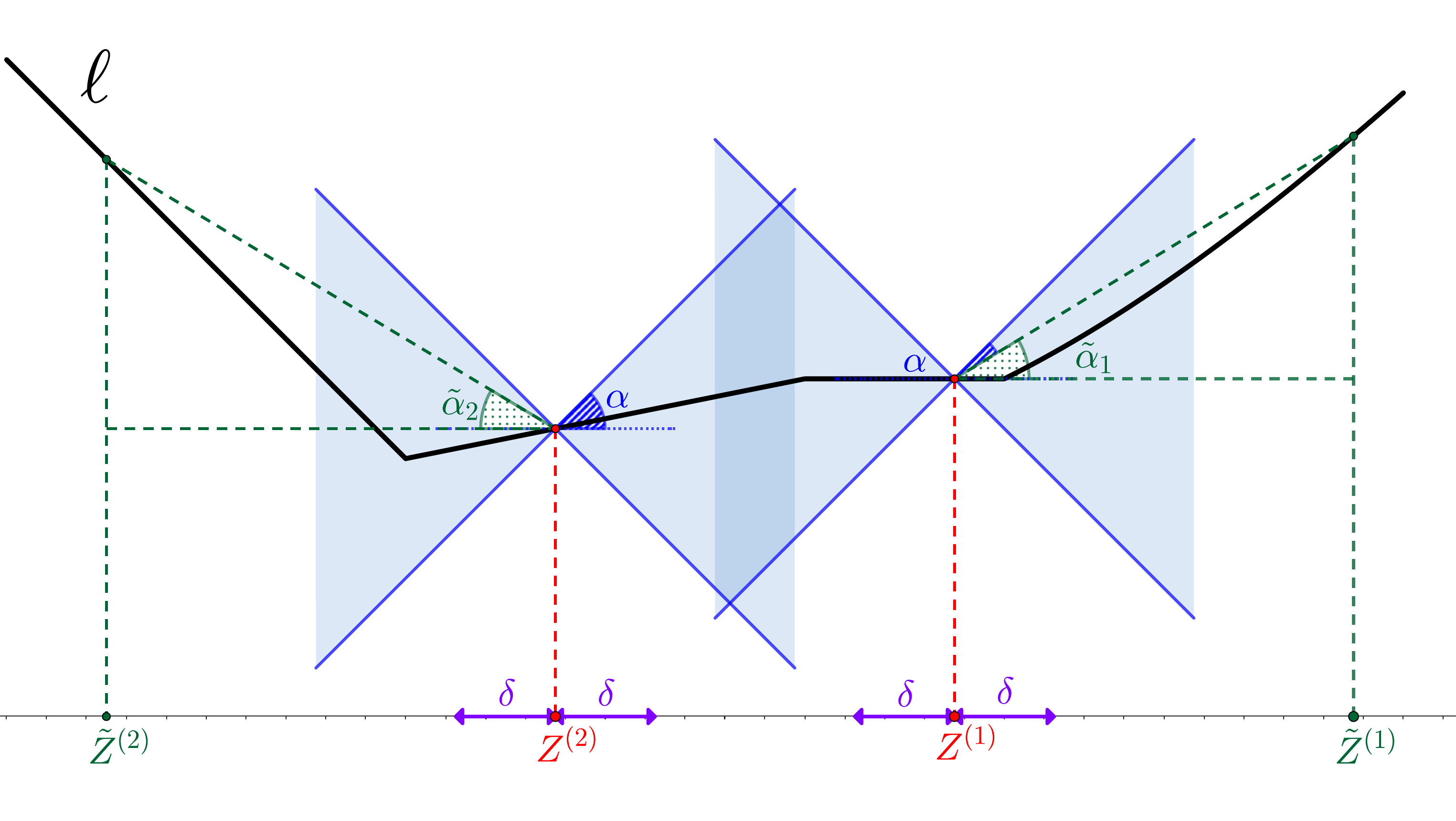}
	\caption{An illustration of Assumptions (A1-A2) (best viewed in color) when \(\mathcal{Z}=\mathbb{R}, \mathcal{Z}_{N}=\{Z^{(1)},Z^{(2)} \}  \) and \(d(z',z) = \left|z'-z\right| \).}
	\label{fig:thm3_2}
\end{figure}

For better illustration, we give a visualization of Assumptions (A1-A2) in Figure~\ref{fig:thm3_2} under the setting where \(\mathcal{Z}=\mathbb{R}, \mathcal{Z}_{N}=\{Z^{(1)},Z^{(2)} \} \) and \(d(z',z) = \left|z'-z\right| \). Assumption (A1) says that for any  \(z\in\mathbb{R}\), we have \(\left|\ell(z;\beta) -  \ell(Z^{(i)};\beta)\right|\leq L_{\beta}^{\mathcal{Z}_N}\left|z-Z^{(i)}\right|  \) for \(i=1,2\), which is equivalent to the condition that the graph of \(\ell(\cdot;\beta)\) must stay inside two blue double cones with the opening angle \(\alpha = \operatorname{arctan}\left(L_{\beta}^{\mathcal{Z}_N}\right) \). Assumption (A2) says that for any \(\epsilon\in (0,L_{\beta}^{\mathcal{Z}_N}) \) and $i = 1,2$, there exists \(\tilde{Z}^{(i)}\) depending on \(\epsilon\) such that \(\left|\tilde{Z}^{(i)}-Z^{(i)}\right|\geq\delta\) and \( \ell(\tilde{Z}^{(i)};\beta) - \ell (Z^{(i)};\beta) \geq ( L_{\beta}^{\mathcal{Z}_N} - \epsilon)\left|\tilde{Z}^{(i)}-Z^{(i)}\right|  \). This is equivalent to requiring the existence of a point  $\tilde{Z}^{(i)}$ such that the green slope \(\tan(\tilde{\alpha}_i)\)  of the line passing through \(\left( \tilde{Z}^{(i)}, \ell( \tilde{Z}^{(i)};\beta) \right) \) and \(\left( {Z}^{(i)}, \ell({Z}^{(i)};\beta) \right) \) satisfies \( \tan(\tilde{\alpha}_i) \geq \tan(\alpha)-\epsilon\), while the distance between $\tilde{Z}^{(i)}$ and ${Z}^{(i)}$ is at least \(\delta\).

\begin{remark}\label{remark}
	In Theorem~\ref{thm:main_r1}, it is required that the condition (A2) holds for any \(i=1,\dots,N\), with respect to the same Lipschitz constant \(L_{\beta}^{\mathcal{Z}_N} \). To relax this condition, one might assume that Assumptions (A1-A2) hold at each \( Z^{(i)}\) with a Lipschitz constant \(L_{\beta}^{\{Z^{(i)}\}}\), for $i = 1,\cdots,N$. Even though it might not guarantee that the lower bound and upper bound for \(\mathcal{S}\) coincide as in Theorem~\ref{thm:main_r1}, we show in Appendix \ref{sec: weaker_main_r1} that one still has closed forms for the lower and upper bounds given by
\begin{equation*} 
\begin{array}{ll}
\widehat{\mathcal{L}}=\mathrm{E}_{\mathbb{P}_N}[\ell(Z;\beta)] + \sum_{i=1}^{N} \mu_i L_{\beta}^{\{Z^{(i)}\}}\delta,\\ 
\widehat{\mathcal{U}}=  \mathrm{E}_{\mathbb{P}_N}[\ell(Z;\beta)] + \max_{i=1,\dots,N} L_{\beta}^{\{Z^{(i)}\}}\delta.
\end{array}
\end{equation*}
\end{remark}

It is worth mentioning that our Assumptions (A1) and (A2) are weaker than those made in the literature. First, in the existing works such as \citet[Theorem 4]{shafieezadeh2015distributionally}, \citet[Theorem 9, Theorem 14]{shafieezadeh2019regularization}, \citet[Assumption 1]{an2021generalization} and \citet[Assumption 1(I)]{gao2022finite}, the Lipschitz assumption on the function \(\psi_{\beta}\) needs to hold globally, while our Assumption (A1) only requires it to hold when the second argument of the cost function $d(\cdot,\cdot)$ is one of the empirical points. Later, Example~\ref{exam:non-Lipschitz} will provide an instance wherein the loss function lacks the global Lipschitz continuity property, but it satisfies our weak Lipschitz property in certain scenarios. Second, \citet[Assumption 1(II)]{gao2022finite} requires that the Lipschitz constant is attained at infinity\footnotemark\footnotetext{This means that for any \(i = 1,\cdots,N\), there exists a sequence \(\left\{ \tilde{Z}^{(i)}_{k}\right\} \) such that \(\lim\limits_{k\rightarrow\infty}d(\tilde{Z}^{(i)}_{k},Z^{(i)}) = \infty \) and  \(\lim\limits_{k\rightarrow\infty} \frac{\psi_{\beta}(\tilde{Z}^{(i)}_{k}) - \psi_{\beta}(Z^{(i)})}{d(\tilde{Z}^{(i)}_{k},Z^{(i)})} =  L_{\beta}^{\mathcal{Z}_N} \).}, and \citet[Assumption 10, Assumption 21]{shafieezadeh2019regularization} requires that the Lipschitz constant is attained exactly at a certain point where the derivative also exists, while our Assumption (A2) only requires it to be (approximately) attained at points far enough from the empirical points. Third, we do not assume any convexity of \(\psi_{\beta}\) as in \citet{wu2022generalization}. In the subsequent contexts, Example~\ref{exam:remove_delta}  will provide an instance where the Lipschitz constant is not attained at infinity and the loss function is not convex but the equivalence \eqref{eq: r1} holds conditionally according to our Theorem \ref{thm:main_r1}. Fourth, we do not require \(d(\cdot,\cdot)\) to be positive definite (that is, \(d(z',z)=0\) if and only if \(z'=z\)) as in \citet[Assumption (A1)]{blanchet2019quantifying} and \citet[Assumption 1]{zhang2022simple}. 

Another important observation of Assumption (A2) is that for any $i=1,\cdots,N$, \(d(\tilde{Z}^{(i)}_{\epsilon},Z^{(i)})\) is required to be at least \(\delta\). That is to say, one requires the knowledge of \(\delta\) to tell whether Assumption (A2), and further the equivalence \eqref{eq: r1}, hold or not in practice. At first glance, it seems restrictive, compared with the existing results where the equivalence \eqref{eq: r1} has been studied with arbitrary \(\delta>0\). Fortunately, as it will be shown later in Section~\ref{sec:app}, our Assumption (A2) indeed holds for most of the commonly-used loss functions for any \(\delta>0\). {In particular,  Propositions~\ref{prop:lin_loss},  \ref{prop:nonlin_loss}, \ref{prop:lin_loss2}, and \ref{prop:nonlin_loss2} serve as guidance on finding \(L_{\beta}^{\mathcal{Z}_N} \) and checking the validity of Assumption (A2) (and Assumption (B) later in Theorem~\ref{thm:main_r2}). }

We give the following two examples to show that the two assumptions in Theorem \ref{thm:main_r1} are not removable. Specifically, we will see that with the same loss function $\ell(\cdot,\cdot)$, the equivalence \eqref{eq: r1} holds for some values of $\delta$ and fails for some other values, which indicates that the condition \(d(\tilde{Z}^{(i)},Z^{(i)})\geq\delta\) in Assumption (A2) is not removable. In addition, Example~\ref{exam:remove_delta} also shows that our weak Lipschitz property in Assumption (A1) needs to depend on the empirical distribution \(\mathbb{P}_N\), which provides the evidence that our generalization is essential compared with the existing results which rely on the global Lipschitz property of $\psi_{\beta}$. {Specifically, these two examples illustrate that the classical Lipschitz terminology fails to capture the exact reformulation of certain WDRO problems. Similar to \cite[Remark 3]{kuhn2019wasserstein}, one can see in Example~\ref{exam:non-Lipschitz} and Example~\ref{exam:remove_delta} that the value of  \(L_{\beta}^{\mathcal{Z}_N} \) is not simple in general cases.  Nevertheless, this notion of the weak Lipschitz constant is better (i.e., lower) than the classical Lipschitz constant, and  applicable for more generic class of loss functions. Efficient scheme to compute  \(L_{\beta}^{\mathcal{Z}_N} \) is {an 
interesting topic to explore, and we leave it} for future research.   } The detailed proofs corresponding to these two examples are given in Appendix \ref{sec: proof_ex31} and Appendix \ref{sec: proof_ex32}.

\begin{example}[Binary cross-entropy \citep{yi2004automated,scott2012calibrated,hurtik2022binary}] \label{exam:non-Lipschitz}
	Let the univariate function \(h\colon (0,1)\rightarrow \mathbb{R} \) be defined as 
	\begin{equation*}
		h(t) =	t\log(t) + (1-t)\log(1-t) .
	\end{equation*} 
	Define the loss function $\ell : (0,1)\times(0,1) \rightarrow \mathbb{R}$ as \( \ell(z;\beta) =\psi_{\beta}(z):= h(\beta z)\).  Consider the cost function \(d\colon(0,1)\times(0,1)\rightarrow [0,1) \) defined as \(d(z',z) = \left|z'-z\right| \) for any $z',z\in (0,1)$. Then the following statements hold true.
	\begin{enumerate}[label=(\alph*)]
		\item \(h\) is convex, continuously differentiable, but not globally Lipschitz on \((0,1)\).
		\item Given any \(\beta\in(0,1)\) and \(\hat{z}\in(0,\frac{1}{2}]\). We have that $\psi_{\beta}$ is \((L_{\beta}^{\{\hat{z}\}},d)\)-Lipschitz at \(\{\hat{z}\}\), where \(L_{\beta}^{\{\hat{z}\}} = -\beta\log(\beta\hat{z}) - (1/\hat{z}-\beta)\log(1-\beta\hat{z}) \). Moreover, for \(\mathbb{P}_N = {\boldsymbol \chi}_{\{\hat{z}\}}\), we have the following two results: 
		\begin{enumerate}[label=(b\arabic*)]
			\item if \(0<\delta < \hat{z} \), then  \(\sup_{\mathbb{P}\colon  \mathcal{W}_{d,1}(\mathbb{P},\mathbb{P}_N) \leq \delta }  \mathrm{E}_{\mathbb{P}} [\ell (Z;\beta) ]=\ell(\hat{z};\beta) + L_{\beta}^{\{\hat{z}\}} \delta\);
			\item if \(\delta \geq \hat{z} \), then \(\sup_{\mathbb{P}\colon  \mathcal{W}_{d,1}(\mathbb{P},\mathbb{P}_N) \leq \delta }  \mathrm{E}_{\mathbb{P}} [\ell (Z;\beta) ]=0\).
		\end{enumerate}
	\end{enumerate}
\end{example}

\begin{example}[Hard sigmoid \citep{howard2019searching} / HardTanh \citep{collobert2004large}]
\label{exam:remove_delta}  
Let the univariate function $h\colon\mathbb{R} \rightarrow \mathbb{R} $ be defined as 
\begin{equation*}
	h(t) = \max \left\{ 0, \min \left\{1, \frac{t+1}{2} \right\} \right\}.
\end{equation*}
Define the loss function $\ell:\mathbb{R}^n\times \mathbb{R}^n\rightarrow \mathbb{R}$ as \( \ell(z;\beta) =\psi_{\beta}(z):= h(\langle \beta,z\rangle )\). Consider the cost function \(d\colon\mathbb{R}^n\times\mathbb{R}^n\rightarrow [0,\infty) \) as \(d(z',z) = \|z'-z\|_{{\mathbb{R}^{n}}} \). For any $\beta\in \mathbb{R}^n$, we denote $\alpha_{\beta}$ to be a vector in $\mathbb{R}^n$ satisfying $\|\alpha_{\beta}\|_{{\mathbb{R}^{n}}}=1$ and $\langle\alpha_{\beta},\beta\rangle = \|\beta\|_{\mathbb{R}^{n},*}$.
		
		\begin{enumerate}[label=(\alph*)]
			\item Given scalars \(0<\vartheta_1\leq  \vartheta_2<\infty \) and any vector $\beta\in \mathbb{R}^n$ satisfying \(\vartheta_1\leq \|\beta\|_{\mathbb{R}^{n},*}\leq \vartheta_2 \). Suppose \(\mathbb{P}_N = {\boldsymbol \chi}_{\{\hat{z}\}}\) with \(\hat{z}={\boldsymbol 0}_{\mathbb{R}^n}\), then $\psi_{\beta}$ is \(\left(\frac{\|\beta\|_{\mathbb{R}^{n},*}}{2},d\right)\)-Lipschitz at \(\{\hat{z}\}\). Moreover, 
			\begin{enumerate}[label=(a\arabic*)]
				\item  if \( 0 <\delta\leq \frac{1}{\vartheta_2} \), then \(\sup_{\mathbb{P}\colon  \mathcal{W}_{d,1}(\mathbb{P},\mathbb{P}_N) \leq \delta }  \mathrm{E}_{\mathbb{P}} [\ell (Z;\beta) ] = \ell(\hat{z};\beta) +  \frac{\|\beta\|_{\mathbb{R}^{n},*}}{2}\delta\);
				\item if \(\delta \geq \frac{1}{\vartheta_1} \), then \(\sup_{\mathbb{P}\colon  \mathcal{W}_{d,1}(\mathbb{P},\mathbb{P}_N) \leq \delta }  \mathrm{E}_{\mathbb{P}} [\ell (Z;\beta) ] = \ell(\hat{z};\beta) +  \frac{1}{2} \). 
			\end{enumerate} 
			\item Given any $\beta\in \mathbb{R}^n$ such that $\|\beta\|_{\mathbb{R}^{n},*}=\vartheta>0$. Suppose \(\mathbb{P}_N = {\boldsymbol \chi}_{\{\bar{z}\}}\) with \(\bar{z}=-\frac{3}{\vartheta}\alpha_{\beta} \), then $\psi_{\beta}$ is \((\frac{\vartheta}{4},d)\)-Lipschitz at \(\{\bar{z}\}\). Moreover, 
			\begin{enumerate}[label=(b\arabic*)]
				\item if \( 0 <\delta\leq \frac{4}{\vartheta} \), then \(\sup_{\mathbb{P}\colon  \mathcal{W}_{d,1}(\mathbb{P},\mathbb{P}_N) \leq \delta }  \mathrm{E}_{\mathbb{P}} [\ell (Z;\beta) ] = \ell(\bar{z};\beta) +  \frac{\vartheta}{4}\delta\);
				\item if \(\delta \geq \frac{4}{\vartheta} \), then \(\sup_{\mathbb{P}\colon  \mathcal{W}_{d,1}(\mathbb{P},\mathbb{P}_N) \leq \delta }  \mathrm{E}_{\mathbb{P}} [\ell (Z;\beta) ] = \ell(\bar{z};\beta) +  1 \). 
			\end{enumerate}
		\end{enumerate} 
	\end{example}
	
\subsection{Equivalence in \eqref{eq:main_equation} when $r>1$}
Next, we derive a sufficient condition on when the equivalence \eqref{eq:main_equation} holds for \(r>1\). 
\begin{theorem}\label{thm:main_r2}
Let \(\mathcal{Z}_N\coloneqq\{Z^{(1)},\dots,Z^{(N)}\}\subset\mathcal{Z} \) be a given dataset and \(\mathbb{P}_N\coloneqq\sum_{i=1}^{N}\mu_i{\boldsymbol \chi}_{\{Z^{(i)}\}}\in\mathcal{P}(\mathcal{Z})\) be the corresponding empirical distribution. In addition, let \(d(\cdot,\cdot)\) be a cost function on \(\mathcal{Z}\times\mathcal{Z}\) and  \(\delta\in(0,\infty)\) be a scalar. Suppose the loss function $\ell:\mathcal{Z}\times \mathcal{B}\rightarrow \mathbb{R}$ takes the form as
\begin{equation*}
\ell\colon (z;\beta) \mapsto \psi_{\beta}^{r}(z), \mbox{ with } r\in(1,\infty),
\end{equation*}
where $\psi_{\beta}\colon\mathcal{Z}\rightarrow \mathbb{R}_{+}$ satisfies Assumption (A1) in Theorem~\ref{thm:main_r1} with \(L_{\beta}^{\mathcal{Z}_N}\in(0,\infty)\), and also satisfies the following assumption:
\begin{enumerate}[label=(B)]
	\item for any \(\epsilon\in (0,L_{\beta}^{\mathcal{Z}_N}) \) and \(Z^{(i)}\in\mathcal{Z}_N\), there exists \(\tilde{Z}^{(i)}_{\epsilon}\in\mathcal{Z}\) such that \(d(\tilde{Z}^{(i)}_{\epsilon},Z^{(i)}) \in  \mathcal{D}(Z^{(i)})\) and
	\begin{equation*}
	\psi_{\beta}(\tilde{Z}^{(i)}_{\epsilon}) - \psi_{\beta}(Z^{(i)}) \geq (L_{\beta}^{\mathcal{Z}_N}-\epsilon) d(\tilde{Z}^{(i)}_{\epsilon},Z^{(i)}),
	\end{equation*}
	where the set \(\mathcal{D}(Z^{(i)})\subset \mathbb{R}\) is defined as
	\begin{equation*}
	\begin{cases}
	\left\{ \frac{\psi_{\beta}(Z^{(i)})}{\left(\mathrm{E}_{\mathbb{P}_N}[\ell(Z;\beta)]\right)^{\frac{1}{r}}}\delta \right\} & \text{if } \mathrm{E}_{\mathbb{P}_N}[\ell(Z;\beta)]\ne 0, \\
	[\delta,\infty) & \text{if } \mathrm{E}_{\mathbb{P}_N}[\ell(Z;\beta)] =  0.
	\end{cases}
	\end{equation*} 
\end{enumerate}
Then we have \(\mathcal{S}=\mathcal{U} \) in Theorem~\ref{thm:main}, that is,
\begin{equation*}
	\sup_{\mathbb{P}  \colon   \mathcal{W}_{d,r}(\mathbb{P},\mathbb{P}_N) \leq \delta }     \mathrm{E}_{\mathbb{P}} [\ell (Z;\beta) ]  =  \left( \left(\mathrm{E}_{\mathbb{P}_N}[\ell(Z;\beta)]\right)^{\frac{1}{r}}   +   L_{\beta}^{\mathcal{Z}_N}\delta\right)^{r}.
\end{equation*}
\end{theorem}
\textbf{Proof.}
Let \(\epsilon \in \left(0,\min \{L_{\beta}^{\mathcal{Z}_N}, \delta L_{\beta}^{\mathcal{Z}_N} \}  \right) \) be a scalar. According to Assumption (B), for any \(Z^{(i)}\in\mathcal{Z}_N\), there exists \(\tilde{Z}^{(i)}\) such that \(d(\tilde{Z}^{(i)},Z^{(i)}) \in\mathcal{D}(Z^{(i)}) \)  and \(\psi_{\beta}(\tilde{Z}^{(i)}) - \psi_{\beta}(Z^{(i)}) \geq  \left(L_{\beta}^{\mathcal{Z}_N}-\frac{\epsilon}{\delta}\right)d(\tilde{Z}^{(i)},Z^{(i)}) \). We consider two cases.
	
\paragraph{Case 1.} When $\mathrm{E}_{\mathbb{P}_N}[\ell(Z;\beta)]\ne 0$, we know that for each \(i = 1,\cdots,N\),
\begin{equation*}
	d(\tilde{Z}^{(i)},Z^{(i)}) = \frac{\psi_{\beta}(Z^{(i)})}{\left(\mathrm{E}_{\mathbb{P}_N}[\ell(Z;\beta)]\right)^{\frac{1}{r}}}\delta.
\end{equation*}
Let \(\tilde{\mathbb{P}} \coloneqq \sum_{i=1}^{N}\mu_i{\boldsymbol \chi}_{\{\tilde{Z}^{(i)}\}} \)and \(\tilde{\pi}\coloneqq \sum_{i=1}^{N}\mu_i{\boldsymbol \chi}_{\{\tilde{Z}^{(i)}\}}\otimes{\boldsymbol \chi}_{\{Z^{(i)}\}} \). Then it can be seen that \(\tilde{\mathbb{P}}\in\mathcal{P}(\mathcal{Z})\), \(\tilde{\pi}\in\Pi(\tilde{\mathbb{P}},\mathbb{P}_N) \), since for any measurable sets \(A,B\subset\mathcal{Z}\),
\begin{equation*}
	\begin{array}{llll}
	\tilde{\pi}(A\times\mathcal{Z})&= \sum_{i=1}^{N}\mu_i{\boldsymbol \chi}_{\{\tilde{Z}^{(i)}\}}(A){\boldsymbol \chi}_{\{Z^{(i)}\}}(\mathcal{Z})
	&=\sum_{i=1}^{N}\mu_i{\boldsymbol \chi}_{\{\tilde{Z}^{(i)}\}}(A)=\tilde{\mathbb{P}}(A),\\
	\tilde{\pi}(\mathcal{Z}\times B)& =\sum_{i=1}^{N}\mu_i{\boldsymbol \chi}_{\{\tilde{Z}^{(i)}\}}(\mathcal{Z}){\boldsymbol \chi}_{\{Z^{(i)}\}}(B) 
	&= \sum_{i=1}^{N}\mu_i{\boldsymbol \chi}_{\{Z^{(i)}\}}(B) = \mathbb{P}_N(B).
\end{array}
\end{equation*}
In addition, it can be seen that
\begin{equation*}
	\begin{array}{ll}
	\mathcal{W}_{d,r}\left(\tilde{\mathbb{P}},\mathbb{P}_N\right)  \leq \left(\int_{\mathcal{Z}\times\mathcal{Z}}d^{r}(\tilde{z},z)\mathrm{d}\tilde{\pi}(\tilde{z},z)\right)^{\frac{1}{r}}
	= \left( \sum_{i=1}^{N}  \mu_i d^{r}(\tilde{Z}^{(i)},Z^{(i)})  \right)^{\frac{1}{r}} =\delta. 
	\end{array}
\end{equation*}
Moreover, we have that
\begin{equation*}
	\begin{array}{ll}
	\left(\mathrm{E}_{\tilde{\mathbb{P}}} [\ell(Z;\beta)]\right)^{\frac{1}{r}} &= \left(\sum_{i=1}^N \mu_i \psi_{\beta}^r(\tilde{Z}^{(i)})\right)^{\frac{1}{r}}
	= \left(\sum_{i=1}^N \mu_i \left(\psi_{\beta}(Z^{(i)})+ \psi_{\beta}(\tilde{Z}^{(i)})-\psi_{\beta}(Z^{(i)})\right)^{r}\right)^{\frac{1}{r}} \\
	&\geq  \left(\sum\limits_{i=1}^N \mu_i \left(\psi_{\beta}(Z^{(i)})+  \left(L_{\beta}^{\mathcal{Z}_N}-\frac{\epsilon}{\delta}\right)d(\tilde{Z}^{(i)},Z^{(i)}) \right)^{r}\right)^{\frac{1}{r}} \\
	&=^{(\Delta)} \left(\sum_{i=1}^N \mu_i \psi_{\beta}^r(Z^{(i)})\right)^{\frac{1}{r}}
	 + \left(L_{\beta}^{\mathcal{Z}_N}-\frac{\epsilon}{\delta}\right)\left(\sum_{i=1}^N \mu_i d^r(\tilde{Z}^{(i)},Z^{(i)})\right)^{\frac{1}{r}}\\
	&= \left(\mathrm{E}_{\mathbb{P}_N}[\ell(Z;\beta)]\right)^{\frac{1}{r}} + (L_{\beta}^{\mathcal{Z}_N}\delta-\epsilon),
	\end{array}
\end{equation*}	
where the equality $^{(\Delta)}$ follows from the fact that for any \(i\in \{1,\cdots,N\}\), we always have $\psi_{\beta}(Z^{(i)}) = \left[ \left(\mathrm{E}_{\mathbb{P}_N}[\ell(Z;\beta)]\right)^{\frac{1}{r}} /\delta\right]d(\tilde{Z}^{(i)},Z^{(i)})$.
	
\paragraph{Case 2.} When \(\mathrm{E}_{\mathbb{P}_N}[\ell(Z;\beta)] =  0\), we have $d(\tilde{Z}^{(i)},Z^{(i)})\in [\delta,\infty)$ for any $i = 1,\cdots,N$. Set \(\eta^{(i)} \coloneqq \delta^{r}/d^r(\tilde{Z}^{(i)},Z^{(i)})\in(0,1]\) for each \(i=1,\dots,N\), and define
\begin{equation*}
	\begin{array}{lll}
	\tilde{\mathbb{P}} &\coloneqq \sum_{i=1}^{N}  \mu_i \eta^{(i)} {\boldsymbol \chi}_{\{\tilde{Z}^{(i)}\}}   +\mu_i(1-\eta^{(i)} ){\boldsymbol \chi}_{\{Z^{(i)}\}},\\
	\tilde{\pi}&\coloneqq\sum_{i=1}^{N} \mu_i\eta^{(i)} {\boldsymbol \chi}_{\{\tilde{Z}^{(i)}\}} \otimes {\boldsymbol \chi}_{\{Z^{(i)}\}}+\mu_i(1-\eta^{(i)} ){\boldsymbol \chi}_{\{Z^{(i)}\}}\otimes{\boldsymbol \chi}_{\{Z^{(i)}\}}.
	\end{array}
\end{equation*}
Then we can see that \(\tilde{\mathbb{P}}\in\mathcal{P}(\mathcal{Z})\) and \(\tilde{\pi}\in\Pi(\tilde{\mathbb{P}},\mathbb{P}_N) \), as for any measurable sets \(A,B\subset\mathcal{Z}\),
\begin{equation*}
	\begin{array}{lll}
	\tilde{\pi}(A\times\mathcal{Z})&=  \sum_{i=1}^{N}   \mu_i\eta^{(i)} {\boldsymbol \chi}_{\{\tilde{Z}^{(i)}\}} (A)+\mu_i(1-\eta^{(i)} ){\boldsymbol \chi}_{\{Z^{(i)}\}}(A) &= \tilde{\mathbb{P}} (A),\\
	\tilde{\pi}(\mathcal{Z}\times B) &= \sum_{i=1}^{N}   \mu_i\eta^{(i)} {\boldsymbol \chi}_{\{Z^{(i)}\}} (B)+\mu_i(1-\eta^{(i)} ){\boldsymbol \chi}_{\{Z^{(i)}\}} (B)
	&=\sum_{i=1}^{N}  \mu_i {\boldsymbol \chi}_{\{Z^{(i)}\}} (B) = \mathbb{P}_N(B).
	\end{array}
\end{equation*}
Moreover, we have \(\int\limits_{\mathcal{Z}\times\mathcal{Z}}d^{r}(\tilde{z},z)\mathrm{d}\tilde{\pi}(\tilde{z},z) =  \sum\limits_{i=1}^{N} \mu_i\eta^{(i)}d^{r}(\tilde{Z}^{(i)},Z^{(i)}) +  \mu_i(1-\eta^{(i)}) d^{r}(Z^{(i)},Z^{(i)}) = \sum\limits_{i=1}^{N} \mu_i\eta^{(i)}d^{r}(\tilde{Z}^{(i)},Z^{(i)})\). Thus, it holds that
\begin{equation*}
	\begin{array}{ll}
	\mathcal{W}_{d,r}\left(\tilde{\mathbb{P}},\mathbb{P}_N\right)  \leq \left(\int_{\mathcal{Z}\times\mathcal{Z}}d^{r}(\tilde{z},z)\mathrm{d}\tilde{\pi}(\tilde{z},z)\right)^{\frac{1}{r}}
	= \left( \sum_{i=1}^{N}  \mu_i\eta^{(i)} d^{r}(\tilde{Z}^{(i)},Z^{(i)})   \right)^{\frac{1}{r}}=\delta. 
	\end{array}
\end{equation*}
The fact that \(\mathrm{E}_{\mathbb{P}_N}[\ell(Z;\beta)] =  0\) together with the nonnegativity of the function $\psi_{\beta}$ implies that $\psi_{\beta}(Z^{(i)})=0$ for any $i=1,\cdots,N$, which further indicates that
\begin{equation*}
	\begin{array}{ll}
	\left(\mathrm{E}_{\tilde{\mathbb{P}}} [\ell(Z;\beta)]\right)^{\frac{1}{r}} &
	= \left(\sum_{i=1}^N \mu_i \eta^{(i)} \psi_{\beta}^r (\tilde{Z}^{(i)}) + \mu_i (1-\eta^{(i)}) \psi_{\beta}^r (Z^{(i)}) \right)^{\frac{1}{r}} \\
	&= \left(\sum_{i=1}^N \mu_i \eta^{(i)} \left( \psi_{\beta} (\tilde{Z}^{(i)})-\psi_{\beta}(Z^{(i)})\right)^{r} \right)^{\frac{1}{r}} \\
	&\geq  \left(L_{\beta}^{\mathcal{Z}_N}-\frac{\epsilon}{\delta}\right)\left(\sum_{i=1}^N \mu_i \eta^{(i)}  d^r(\tilde{Z}^{(i)},Z^{(i)}) \right)^{\frac{1}{r}}
	= L_{\beta}^{\mathcal{Z}_N}\delta-\epsilon
	= \left(\mathrm{E}_{\mathbb{P}_N}[\ell(Z;\beta)]\right)^{\frac{1}{r}} + L_{\beta}^{\mathcal{Z}_N}\delta-\epsilon.
	\end{array}
\end{equation*}	
	
Therefore, combining the above two cases, for any given \(0<\epsilon<\min \{L_{\beta}^{\mathcal{Z}_N}, \delta L_{\beta}^{\mathcal{Z}_N} \} \), we can construct \(\tilde{\mathbb{P}} \in\mathcal{P}(\mathcal{Z})\) such that \(\mathcal{W}_{d,r}\left(\tilde{\mathbb{P}},\mathbb{P}_N\right)\leq \delta \) and 
\begin{equation*}
	\mathrm{E}_{\tilde{\mathbb{P}}}[\ell(Z;\beta)] \geq \left(\mathrm{E}_{\mathbb{P}_N}[\ell(Z;\beta)]^{\frac{1}{r}} + L_{\beta}^{\mathcal{Z}_N}\delta-\epsilon\right)^{r}.
\end{equation*}
By letting $\epsilon\rightarrow 0$ and recalling that  \(\mathcal{U}=\left(\mathrm{E}_{\mathbb{P}_N}[\ell(Z;\beta)]^{\frac{1}{r}} + L_{\beta}^{\mathcal{Z}_N}\delta\right)^{r} \), we can conclude that  
\begin{equation*}
	\sup_{\mathbb{P}\colon  \mathcal{W}_{d,r}(\mathbb{P},\mathbb{P}_N) \leq \delta }  \mathrm{E}_{\mathbb{P}} [\ell (Z;\beta) ] \geq \mathcal{U}.
\end{equation*}
The fact that $\psi_{\beta}$ satisfies Assumption (A1) together with Theorem~\ref{thm:main} also implies that  
\begin{equation*}
	\sup_{\mathbb{P}\colon  \mathcal{W}_{d,r}(\mathbb{P},\mathbb{P}_N) \leq \delta }  \mathrm{E}_{\mathbb{P}} [\ell (Z;\beta) ] \leq \mathcal{U},
\end{equation*}
which completes the proof.
\hfill\(\square\) 

Note that when the loss quantity at the empirical distribution \(\mathrm{E}_{\mathbb{P}_N}[\ell(Z;\beta)]  \) vanishes, one can see that Assumption (A2) and Assumption (B) are the same. On the other hand, when \(\mathrm{E}_{\mathbb{P}_N}[\ell(Z;\beta)] \ne 0 \), if \(\psi_{\beta}(Z^{(i)})=0\) for some $i=1,\cdots,N$, one can always choose \(\tilde{Z}^{(i)}_{\epsilon}= Z^{(i)}\) for this \(i\). 
	
\section{Applications to different function classes} 
\label{sec:app}
In this section, we present the applications of Theorem~\ref{thm:main_r1} and Theorem~\ref{thm:main_r2} in various regression and classification problems. We first give the definition of an absolutely homogeneous function, which will be used in the remaining part of this section. 
\begin{definition} \label{def:abs-homo}
	An extended-valued function  \(\Upsilon \colon\mathcal{Z}\rightarrow [0,\infty] \) on a real vector space \(\mathcal{Z}\) is called absolutely homogeneous if one has \(\Upsilon(tz) = \left|t\right|\Upsilon(z) \) for any \(t\in\mathbb{R}\) and \(z\in\mathcal{Z}\). In addition, \(\Upsilon\) is called proper if there exists \(z_0\in\mathcal{Z}\) such that \(\Upsilon(z_0) =1 \).
\end{definition}

From the above definition, we can see that if \(\Upsilon\) is absolutely homogeneous, then \(\Upsilon(0_{\mathcal{Z}})=0\) and \(\Upsilon^{-1}(0) = \left\{ z \in \mathcal{Z} \mid \Upsilon(z)=0 \right\} \) is a cone on \(\mathcal{Z}\). Besides, the following two functions
\begin{itemize}
	\item \(z \mapsto +\infty \) for all \(0\ne z\in\mathcal{Z}\) and \(z \mapsto 0 \) for \(z=0\);
	\item \(z \mapsto 0 \) for all \( z\in\mathcal{Z}\),
\end{itemize}
are absolutely homogeneous, but not proper. 
	
\subsection{Applications to simple piecewise linear regression loss functions}
We start from the applications of our results to simple piecewise linear regression loss functions as follows. 
\begin{proposition}[Linear loss] \label{prop:lin_loss} 
Let \(\mathcal{Z}\) be a (finite or infinite dimensional) real vector space. Suppose that  \(\left\llbracket\cdot\right\rrbracket\colon\mathcal{Z}\rightarrow[0,\infty]\) is absolutely homogeneous and proper, \(\phi\colon \mathcal{Z}\rightarrow\mathbb{R} \) is linear, \(\left\llbracket \cdot\right\rrbracket^{-1}(0) \subseteq \phi^{-1}(0)\) and
\begin{equation*}
L_{\phi} \coloneqq \sup_{z\in\mathcal{Z}} \left\{ \left|\phi(z)\right| \mid \left\llbracket z\right\rrbracket = 1  \right\} \in [0,+\infty).
\end{equation*}
Let the cost function $d:\mathcal{Z}\times \mathcal{Z}\rightarrow [0,\infty]$ be defined as \(d(z',z)\coloneqq \left\llbracket z'-z\right\rrbracket \) for any $z',z\in \mathcal{Z}$. Given any scalar $\tau\in \mathbb{R}$, then the functions \(\left|\phi\right|\), \(\max\{0,\phi-\tau\}\) and \(\max\{0,\left|\phi\right|-\tau\}\) are $(L_{\phi},d)$-Lipschitz at $\mathcal{Z}$.  Furthermore, they also satisfy Assumptions (A1), (A2) and (B) at $\mathcal{Z}$ for any $\delta>0$ if \(L_{\phi}>0\).
\end{proposition}
\textbf{Proof.}
In order to prove that \(|\phi|\), \(\max\{0,\phi-\tau\}\) and \(\max\{0,|\phi|-\tau\}\) are $(L_{\phi},d)$-Lipschitz at $\mathcal{Z}$, by noting the fact that for any \(z',z\in\mathcal{Z}\),
%\begin{equation*}
%\begin{array}{llll}
%	&\left||\phi(z')| - |\phi(z)| \right| \leq \left|\phi(z') - \phi(z) \right|,\\ 
%	&\left|\max\{0,\phi(z')  -  \tau\}  -  \max\{0,\phi(z)  -  \tau\} \right|  \leq  \left|\phi(z') - \phi(z)\right|,\\
%	&\left|\max\{0,|\phi(z')|-\tau\} - \max\{0,|\phi(z)|-\tau\} \right| \\
%	&\qquad \leq \left||\phi(z')| - |\phi(z)| \right|  \leq \left|\phi(z')-\phi(z)\right|,
%\end{array}
%\end{equation*}
\begin{equation*}
	\begin{array}{lll}
		\left||\phi(z')| - |\phi(z)| \right| &\leq \left|\phi(z') - \phi(z) \right|;\\ 
		\multicolumn{2}{l}{\left|\max\{0,\phi(z')-\tau\} - \max\{0,\phi(z)-\tau\} \right|} 
		&\leq \left|\phi(z')-\phi(z)\right|;\\
		\multicolumn{2}{l}{\left|\max\{0,|\phi(z')|-\tau\} - \max\{0,|\phi(z)|-\tau\} \right|} 
		&\leq   \left||\phi(z')|  -  |\phi(z)| \right|  \leq \left|\phi(z') - \phi(z)\right|;
	\end{array}
\end{equation*}
we only need to prove that
\begin{equation*}
	\left|\phi(z') - \phi(z) \right| = \left|\phi(z'-z) \right|\leq L_{\phi}d(z',z).
\end{equation*}
This can be seen as follows:
\begin{itemize} 
	\item  if  \(d(z',z)=\left\llbracket z'-z\right\rrbracket=\infty \), then it holds true immediately;      
	\item if \(d(z',z) =\left\llbracket z'-z\right\rrbracket = 0 \), since \( \left\llbracket \cdot\right\rrbracket^{-1}(0) \subseteq \phi^{-1}(0)\), one also has \( \left|\phi(z'-z) \right|=0\);
	\item if \(0<d(z',z)=\left\llbracket z'-z\right\rrbracket <\infty \), then 
	\begin{equation*}
	\begin{array}{ll}
		\left|\phi(z'-z)\right| = \left\llbracket z'-z\right\rrbracket \left| \phi\left(\frac{z'-z}{\left\llbracket z'-z\right\rrbracket }\right)\right|  
		\leq \left\llbracket z'-z\right\rrbracket  L_{\phi} = L_{\phi}d(z',z).
	\end{array}
	\end{equation*}
\end{itemize}
Thus, \(|\phi|\), \(\max\{0,\phi-\tau\}\) and \(\max\{0,|\phi|-\tau\}\) are all \((L_{\phi},d)\)-Lipschitz at \(\mathcal{Z}\).

Suppose that \(L_{\phi}>0\). From the previous discussion, we have that \(|\phi|\), \(\max\{0,\phi-\tau\}\) and \(\max\{0,|\phi|-\tau\}\) satisfy Assumption (A1) at $\mathcal{Z}$ with any $\delta>0$. By the definition of \(L_{\phi}\), for any $0<\epsilon <L_{\phi}$, there exists \(\tilde{v}\in\mathcal{Z}\) such that \(\left\llbracket \tilde{v} \right\rrbracket  = 1\) and \(\left|\phi(\tilde{v})\right| \geq L_{\phi}-\epsilon/2\). Since \(\phi\) is linear and \(\left\llbracket \cdot\right\rrbracket\) is absolutely homogeneous, for any $v\in \mathcal{Z}$, we have \(\left\llbracket -v\right\rrbracket=\left\llbracket v\right\rrbracket\) and \(\phi(-v)=-\phi(v) \). Hence, one can always choose \(\tilde{v}\in\mathcal{Z}\) such that 
\begin{equation*}
	\left\llbracket \tilde{v}\right\rrbracket = 1 \quad \text{and}\quad  \phi(\tilde{v}) \geq L_{\phi}-\epsilon/2> L_{\phi}-\epsilon>0.
\end{equation*}
Next we prove that the functions \(|\phi|\), \(\max\{0,\phi-\tau\}\) and \(\max\{0,|\phi|-\tau\}\) satisfy Assumptions (A2) and (B) at $\mathcal{Z}$ for any $\delta>0$.
\begin{itemize}
	\item For the function $|\phi|$, given any $z\in \mathcal{Z}$ and any $\sigma>0$, by letting $\tilde{z} = z +{\rm sgn}(\phi(z))\sigma \tilde{v}$, we have \( d(\tilde{z},z) = \left\llbracket {\rm sgn}(\phi(z))\sigma \tilde{v}\right\rrbracket = \sigma  \) and 
	\begin{equation*}
	\begin{array}{ll}
		\left| \phi(\tilde{z})\right|- \left|\phi(z)\right| = \left|\phi(z) + {\rm sgn}(\phi(z))\sigma \phi(\tilde{v})\right|- \left|\phi(z)\right|
		=\sigma \phi(\tilde{v})  \geq  (L_{\phi}-\epsilon) d(\tilde{z},z).
	\end{array}
	\end{equation*}
	Therefore, \(|\phi|\) satisfies Assumption (A2) at \(\mathcal{Z}\) for any \(\delta>0\). Using some similar analysis as above, we can also show that \(|\phi|\) satisfies Assumption (B) at \(\mathcal{Z}\). 
	\item For the function \(\max\{0,\phi-\tau\} \), we first prove that \(\max\{0,\phi-\tau\} \) satisfies Assumption (A2) at $\mathcal{Z}$. For any $z\in \mathcal{Z}$ and \(\delta>0 \), let
	\begin{equation*}
		\tilde{z} = \left\{\begin{array}{ll}
			z + \delta \tilde{v} &\text{if } \phi(z) \geq \tau\\
			z + \left( 2(\tau-\phi(z))/\epsilon +\delta\right) \tilde{v} &\text{otherwise}
		\end{array}\right. .
	\end{equation*}
	Then if $ \phi(z) \geq \tau$, we have \( d(\tilde{z},z) = \left\llbracket \delta \tilde{v}\right\rrbracket = \delta  \) and 
	\begin{equation*}
	\begin{array}{ll}
		\max\{0,\phi(\tilde{z})-\tau \} - \max\{0,\phi(z)-\tau \} 
		\geq  \phi(z+\delta\tilde{v}) - \phi(z) = \delta \phi(\tilde{v})\geq  (L_{\phi}-\epsilon) d(\tilde{z},z);
	\end{array}
	\end{equation*}
	if $ \phi(z) < \tau$, we have \( d(\tilde{z},z) = \left\llbracket \left( 2(\tau-\phi(z))/\epsilon +\delta\right) \tilde{v}\right\rrbracket =2(\tau-\phi(z))/\epsilon +\delta \geq \delta  \) and
	\begin{equation*}
	\begin{array}{lll}
		\max\{0,\phi(\tilde{z})-\tau \} - \max\{0,\phi(z)-\tau \} \\
		 \geq \phi(z)-\tau + d(\tilde{z},z)  \phi(\tilde{v}) \geq -\frac{\epsilon}{2}d(\tilde{z},z) + d(\tilde{z},z)  \phi(\tilde{v})
		 \geq (L_{\phi}-\epsilon) d(\tilde{z},z).
	\end{array}
	\end{equation*}
	Therefore, \(\max\{0,\phi-\tau\}\) satisfies Assumptions (A2) at $\mathcal{Z}$ for any $\delta>0$.
	
	Next, we turn to Assumption (B). Fix $r> 1$, $\delta>0$, a dataset \(\mathcal{Z}_N\subset\mathcal{Z} \) and the corresponding empirical distribution \(\mathbb{P}_N \). For any $\hat{z}\in \mathcal{Z}_N$, we consider the following cases.
	\begin{itemize}
		\item If \(\mathrm{E}_{\mathbb{P}_N}[\max\left\{0, \phi(Z)-\tau \right\}^r ]=0 \), then (A2) and (B) are equivalent. 
		\item If \(\mathrm{E}_{\mathbb{P}_N}[\max\left\{0, \phi(Z)-c \right\}^r ]\ne0 \) and \( \phi(\hat{z})>\tau  \), for any \(\sigma\geq \delta\), we can set \(\tilde{z} = \hat{z} + \sigma\tilde{v}\). Then we have \(\phi(\tilde{z}) = \phi(\hat{z}) + \sigma\phi(\tilde{v}) > \tau + \sigma(L_{\phi}-\epsilon/2) > \tau \). Therefore, \(d(\tilde{z},\hat{z}) = \left\llbracket \sigma\tilde{v}\right\rrbracket = \sigma \) and
		\begin{equation*}
		\begin{array}{ll}
			\max\left\{0, \phi(\tilde{z})-\tau \right\}-\max\left\{0, \phi(\hat{z})-\tau \right\} 
			= \phi(\tilde{z}) - \phi(\hat{z}) = \sigma\phi(\tilde{v}) \geq (L_{\phi}-\epsilon)d(\tilde{z},\hat{z}).
		\end{array}
		\end{equation*}
		\item If \(\mathrm{E}_{\mathbb{P}_N}[\max\left\{0, \phi(Z)-\tau \right\}^r ]\ne0 \) and \( \phi(\hat{z})\leq \tau  \), then \(\max\{0, \phi(\hat{z})-\tau\}=0 \) and one can choose \(\tilde{z}=\hat{z}\) such that (B) holds.
	\end{itemize}
	This means that \(\max\{0,\phi-\tau\}\) satisfies Assumptions (B) at $\mathcal{Z}$ for any $\delta>0$.
	\item Finally, we consider the function \(\max\{0,|\phi|-\tau\} \). For any $z\in \mathcal{Z}$ and \(\delta>0 \), let
%	\begin{equation*}
%		\tilde{z}  =   \left\{
%		\begin{array}{ll}
%			z  +  {\rm sgn}(\phi(z)) \delta \tilde{v} &\text{if } |\phi(z)|   \geq   \tau,\\
%			z  +  {\rm sgn}(\phi(z)) \left( 2(\tau  -  |\phi(z)|)/\epsilon   +  \delta\right) \tilde{v} &\text{otherwise}.
%		\end{array}\right.
%	\end{equation*}
	\begin{equation*}
		\tilde{z}  =   \left\{
		\begin{array}{ll}
			z + {\rm sgn}(\phi(z)) \delta \tilde{v} &\text{if } |\phi(z)| \geq  \tau,\\
			z + {\rm sgn}(\phi(z)) \left( \frac{2(\tau-|\phi(z)|)}{\epsilon}+\delta\right) \tilde{v} &\text{otherwise}.
		\end{array}\right.
	\end{equation*}
	Then if $ |\phi(z)| \geq \tau$, we have \( d(\tilde{z},z) = \left\llbracket {\rm sgn}(\phi(z)) \delta \tilde{v}\right\rrbracket = \delta  \) and 
	\begin{equation*}
	\begin{array}{lll}
		\max\{0,\left|\phi(\tilde{z})\right|-\tau \} - \max\{0,|\phi(z)|-\tau \} 
		\geq  \left|\phi(z + {\rm sgn}(\phi(z)) \delta \tilde{v})\right| - |\phi(z)| 
		= \delta \phi(\tilde{v})\geq  (L_{\phi}-\epsilon) d(\tilde{z},z);
	\end{array}
	\end{equation*}
	if $ |\phi(z)| < \tau$, we have \( d(\tilde{z},z) = \left\llbracket {\rm sgn}(\phi(z)) \left( 2(\tau-|\phi(z)|)/\epsilon +\delta\right) \tilde{v}\right\rrbracket =2(\tau-|\phi(z)|)/\epsilon +\delta \geq \delta  \) and
	\begin{equation*}
	\begin{array}{llll}
		&\max\{0,\left|\phi(\tilde{z})\right|-\tau \} - \max\{0,|\phi(z)|-\tau \}  \geq \left|\phi(\tilde{z})\right| - \tau \\
		&= \left|\phi(z) + {\rm sgn}(\phi(z)) d(\tilde{z},z) \phi(\tilde{v})\right|-\tau
		= |\phi(z)|-\tau + d(\tilde{z},z)  \phi(\tilde{v}) \\
		&\geq  -\frac{\epsilon}{2}d(\tilde{z},z) + d(\tilde{z},z)  \phi(\tilde{v}) \geq (L_{\phi}-\epsilon) d(\tilde{z},z).
	\end{array}
	\end{equation*}
	Thus \(\max\{0,|\phi|-\tau\}\) satisfies Assumptions (A2) at $\mathcal{Z}$ for any $\delta>0$. 
	
	Fix $r> 1$, $\delta>0$, a dataset \(\mathcal{Z}_N\subset\mathcal{Z} \) and the corresponding empirical distribution \(\mathbb{P}_N \). For each \(\hat{z}\in\mathcal{Z}_N \), since the cases when (1) \(\mathrm{E}_{\mathbb{P}_N}[\max\left\{0, |\phi(Z)|-\tau \right\}^r ]=0 \) or (2) \(\mathrm{E}_{\mathbb{P}_N}[\max\left\{0, |\phi(Z)|-\tau \right\}^r ]\ne0 \) and $\left|\phi(\hat{z})\right|\leq \tau$, are easy to check, we only need to consider the case when \(\mathrm{E}_{\mathbb{P}_N}[\max\left\{0, |\phi(Z)|-\tau \right\}^r ]\ne0 \) and $|\phi(\hat{z})|>\tau$. For any $\sigma\geq \delta$, let \(\tilde{z} = \hat{z} + {\rm sgn}(\phi(z)) \sigma\tilde{v}\). Then we can see that 
	\begin{equation*}
	\begin{array}{ll}
		|\phi(\tilde{z})|&=\left|\phi(\hat{z}) + {\rm sgn}(\phi(z)) \sigma\phi(\tilde{v})\right| 
		= |\phi(\hat{z})| + \sigma\phi(\tilde{v})>\tau.
	\end{array}
	\end{equation*}
	Moreover, we have \(d(\tilde{z},\hat{z}) = \left\llbracket {\rm sgn}(\phi(z)) \sigma\tilde{v}\right\rrbracket = \sigma \) and 
	\begin{equation*}
	\begin{array}{lll}
		\max\left\{0, |\phi(\tilde{z})|-\tau \right\}-\max\left\{0, |\phi(\hat{z})|-\tau \right\} 
		= \left|\phi(\hat{z}) + {\rm sgn}(\phi(z)) \sigma\phi(\tilde{v})\right|-|\phi(\hat{z})|  
		= \sigma\phi(\tilde{v})\\
		\geq (L_{\phi}-\epsilon)d(\tilde{z},\hat{z}).
	\end{array}
	\end{equation*}
	Therefore, \(\max\{0,|\phi|-\tau\}\) satisfies Assumptions (B) at $\mathcal{Z}_N$ for any $\delta>0$. 
\end{itemize}
This completes the proof.
\hfill\(\square\) 

\begin{remark} 
Proposition~\ref{prop:lin_loss} can be viewed as a guidance to find the weak Lipschitz constant \(L_{\phi}\) when certain properties of the loss and cost functions are given. 
\end{remark}
	
Based on Proposition~\ref{prop:lin_loss}, Theorem~\ref{thm:main}(c), Theorem~\ref{thm:main_r1} and Theorem~\ref{thm:main_r2}, we obtain the following corollary stating the equivalence \eqref{eq:main_equation} for simple piecewise linear regression loss functions.

\begin{corollary}\label{coro: linear}
Under the setting of Proposition \ref{prop:lin_loss},  given any scalar \(\delta>0 \) and any empirical distribution \(\mathbb{P}_N\) on \( \mathcal{Z} \), it holds that for any $r\geq 1,\tau\in\mathbb{R}$, 
\begin{equation*}
\sup_{\mathbb{P}\colon  \mathcal{W}_{d,r}(\mathbb{P},\mathbb{P}_N) \leq \delta }  \mathrm{E}_{\mathbb{P}} [\ell(Z)]  = \left(\left(\mathrm{E}_{\mathbb{P}_N}[\ell(Z)]\right)^{\frac{1}{r}} + L_{\phi}\delta\right)^{r},
\end{equation*}
where \(\ell: \mathcal{Z}\rightarrow \mathbb{R}\) takes one of the following forms: \(\ell(\cdot)= \left|\phi(\cdot)\right|^{r}\), \(\ell(\cdot)= \max\{0,\phi(\cdot)-\tau\}^{r}\) or \(\ell(\cdot)= \max\{0,\left|\phi(\cdot)\right|-\tau\}^{r}\).
\end{corollary}

Here are some specific examples of the above corollary. Note that some of the following results have been studied in the literature, whereas we have provided a unified framework to study the equivalence between the worst-case loss quantity in the WDRO problem and the regularization scheme for simple piecewise linear regression and scalar-on-function linear regression functions.
	
\begin{example}[Linear-type regression] \label{ex:linear_regression} 
Let \(\mathcal{Z}=\mathbb{R}^n\times\mathbb{R} \).  Given any $\delta>0$, \(\beta\in\mathbb{R}^n\) and any empirical distribution \(\mathbb{P}_N\) on \( \mathbb{R}^n\times\mathbb{R}\). For any $r\geq 1$ and $\tau\in \mathbb{R}$, we have
\begin{equation*}
\sup_{\mathbb{P}\in\mathfrak{M}_r }  \mathrm{E}_{\mathbb{P}} [\ell(Z;\beta)]  = \left(\left(\mathrm{E}_{\mathbb{P}_N}[\ell(Z;\beta)]\right)^{\frac{1}{r}} + L_{\phi}(\beta)\delta\right)^{r},
\end{equation*}
where \(Z=(X,Y)\), \(\mathfrak{M}_r\coloneqq \left\{\mathbb{P}\in \mathcal{P}(\mathcal{Z}) \mid  \mathcal{W}_{d,r}(\mathbb{P},\mathbb{P}_N) \leq \delta \right\}\) and  \(\ell(Z;\beta)\) takes one of the following forms:
\begin{enumerate}[label=(\alph*)]
	\item \(\ell(Z;\beta)= \left|Y - \langle\beta,X\rangle\right|^r \);
	\item \(\ell(Z;\beta)= \left( Y - \langle\beta,X\rangle-\tau\right)_{+}^r \);
	\item \(\ell(Z;\beta)=  \left( \left| Y - \langle\beta,X\rangle\right|-\tau\right)_{+}^r \);
\end{enumerate}
and \(d(\cdot,\cdot)\) and \(L_{\phi}(\beta)\) take one of the following forms:
\begin{enumerate}[label=(\roman*)]
	\item $d((x',y'),(x,y)) = \|[x'-x;y'-y]\|_{\mathbb{R}^{n+1}}$ and  \(L_{\phi}(\beta)=\|[-\beta;1]\|_{\mathbb{R}^{n+1},*}\);
	\item $d((x',y'),(x,y)) = \|x'-x\|_{\mathbb{R}^n}+{\boldsymbol \delta}_{\{0\}}(y'-y)$ and \(L_{\phi}(\beta)=\|\beta\|_{\mathbb{R}^n,*}\);
	\item $d((x',y'),(x,y)) = \|x'_{\mathcal{I}}-x_{\mathcal{I}}\|_{\mathbb{R}^{|\mathcal{I}|}}+{\boldsymbol \delta}_{\{{\bf 0}^{\mathbb{R}^{|\mathcal{I}^c|+1}}\}}([x'_{\mathcal{I}^c}-x_{\mathcal{I}^c};y'-y])$ and \(L_{\phi}(\beta)=\|\beta_{\mathcal{I}}\|_{\mathbb{R}^{|\mathcal{I}|},*}\) where \(\mathcal{I}\subset\{1,2,\dots,n\}\) and \(\mathcal{I}^c=\{1,2,\dots,n\}\setminus\mathcal{I}\);
	\item \(d((x',y'),(x,y)) =\inf_{\bar{x}\in\mathbb{R}^s}\left\{\|\bar{x}\|_{\mathbb{R}^s} \mid B^T\bar{x}=x'-x \right\} +{\boldsymbol \delta}_{\{0\}}(y'-y) \) and \(L_{\phi}(\beta)= \|B\beta\|_{\mathbb{R}^{s},*} \) where \(B\in\mathbb{R}^{s\times n} \) is a given matrix.
\end{enumerate}
\end{example}

The proof associated with the above example is given in Appendix \ref{sec: proof_ex41}. Note that this example covers many commonly-used regression problems. Specifically, we list some of them here.
\begin{itemize}
	\item Higher-order regression: $\mathrm{E}_{\mathbb{P}} [|Y - \langle\beta,X\rangle|^r ]$, $r\geq 1$.  When \(r=2\), the regularized problem $\sqrt{\mathrm{E}_{\mathbb{P}_N} [|Y - \langle\beta,X\rangle|^2 ]}  + L_{\phi}(\beta)\delta  $ is often referred as a variant of the square-root Lasso model, where \(L_{\phi}(\beta)\) is a function of \(\beta\) which promotes specific structures in \(\beta\), such as smoothness, sparsity, and clustering of coordinates \citep{belloni2011square,bunea2013group,stucky2017sharp,jiang2021simultaneous}.
	\item Lower partial moments \citep{bawa1975optimal,fishburn1977mean,chen2011tight}: $\mathrm{E}_{\mathbb{P}} [\left(\langle\beta,X\rangle-\tau \right)_{+}^r ]$, $r\geq 1$, $\tau\in \mathbb{R}$.
	\item Higher-order $\tau$-insensitive regression: $\mathrm{E}_{\mathbb{P}} [\left(|Y - \langle\beta,X\rangle|-\tau \right)_{+}^r ]$, $r\geq 1$, $\tau\geq 0$. When $r=1$, it is the $\tau$-insensitive support vector regression \citep{drucker1996support,wang2020comprehensive}; when $r=2$, it is the $\tau$-smooth support vector regression \citep{lee2005epsilon}.
\end{itemize}

In the existing literature, the equivalence between the worst-case loss quantity in the WDRO problem and the regularization scheme for the above problems has been studied for some specific settings. For instances, \citet[Proposition 2, Theorem 1]{blanchet2019robust} studied variants of the square-root Lasso model when \(d(\cdot,\cdot)\) is an extended norm and \citet[Theorem 1]{Chu2021OnRS} covered the cases when \(d(\cdot,\cdot)\) is an extended semi-norm; {\cite{kuhn2019wasserstein} studied tractable reformulation of the WDRO problem where \(\ell\) satisfies certain convex/concave properties;} \citet[Corollary 2]{gao2022wasserstein} investigated more general linear models when \(d(\cdot,\cdot)\) is a metric and \citet{wu2022generalization} explored the piecewise linear convex loss function when \(d(\cdot,\cdot)\) is an extended norm. Moreover, when \(r\geq1\), \citet{montiel2022generalization} considered a similar model but with the max-sliced Wasserstein ball.
	
As a natural extension of the ordinary linear regression problem, the scalar-on-function linear regression problem has received increasing attention nowadays, where the feature space \(\mathcal{X}\) is a functional space \(\mathfrak{L}^2[0,1] \) instead of \(\mathbb{R}^n \), endowed with the inner product \( \langle x',x\rangle=\int_{0}^{1}x'(t)x(t)\mathrm{d}t \). We refer the readers to \citet{ramsay1991some,ramsay_functional_2005,wang2016functional} for more details and discussions about the functional linear regression problem. To the best of our knowledge, the equivalence between the worst-case loss quantity in the WDRO problem and the regularization scheme for this class of problems has not been established in the literature. Fortunately, based on our results, we can give the following equivalence for the scalar-on-function linear regression problems, whose proof can be found in Appendix \ref{sec: proof_ex42}. It is worth mentioning that the nonparametric model (a) has been studied in \citet{cardot1999functional,cai2012minimax} and the regularizer involving \(\int_{0}^{1}|\beta(t)|^2\mathrm{d}t\) has been considered in \citet[(2)]{tong2018analysis}. In addition, the parametric model (b) has been introduced to reduce the degree of freedoms. 

\begin{example}[Scalar-on-function linear regression] \label{eg:functional-LR}
Denote the set of real-valued, square-integrable functions on \([0,1]\) as \( \mathfrak{L}^2[0,1]\).  Given any empirical distribution \(\mathbb{P}_N\) on \( \mathcal{Z} :=\mathfrak{L}^2[0,1]\times\mathbb{R} \) and any scalar $\delta>0$,  the following equality holds for any $r\geq 1$ and $\tau\in \mathbb{R}$:
\begin{equation*}
\sup_{\mathbb{P}\colon  \mathcal{W}_{d,r}(\mathbb{P},\mathbb{P}_N) \leq \delta }  \mathrm{E}_{\mathbb{P}} [\ell(Z)]  = \left(\left(\mathrm{E}_{\mathbb{P}_N}[\ell(Z)]\right)^{\frac{1}{r}} + L_{\phi}\delta\right)^{r},
\end{equation*}
where \(d((x',y'),(x,y)) =  \left(\int_{0}^{1}|x'(t)-x(t)|^2\mathrm{d}t\right)^{1/2} + {\boldsymbol \delta}_{\{0\}}(y'-y) \) for any $(x',y'),(x,y)\in \mathfrak{L}^2[0,1]\times\mathbb{R}$;  \(\ell(Z)= |\phi(Z)|^{r}\), \(\ell(Z)= \max\{0,\phi(Z)-\tau\}^{r}\) or \(\ell(Z)= \max\{0,|\phi(Z)|-\tau\}^{r}\); and $\phi$, $L_{\phi}$ take one of the following forms. 
\begin{enumerate}[label=(\alph*)]
	\item (Nonparametric) Given \(\beta\in\mathfrak{L}^2[0,1]\) and let  $\phi:\mathfrak{L}^2[0,1]\times\mathbb{R}\rightarrow \mathbb{R}$ defined as
	\begin{equation*}
	\phi\colon(x,y)\mapsto y - \int\limits_{0}^{1} x(t) \beta(t) \mathrm{d}t,
	\end{equation*}
	for any \((x,y)\in \mathfrak{L}^2[0,1]\times\mathbb{R}\). Let \(L_{\phi}=\left(\int_{0}^{1}|\beta(t)|^2\mathrm{d}t\right)^{1/2} \).
	\item (Parametric) Let \(\beta\in\mathbb{R}^n\), \(\{{\boldsymbol g}_1,\cdots,{\boldsymbol g}_n \} \subset \mathfrak{L}^2[0,1]\). Define $\phi:\mathfrak{L}^2[0,1]\times\mathbb{R}\rightarrow \mathbb{R}$ as 
	\begin{equation*}
	\phi\colon (x,y)\mapsto y - \int\limits_{0}^{1} x(t) \sum_{j=1}^{n} \beta_j {\boldsymbol g}_j(t)\mathrm{d}t,
	\end{equation*}    
	for any \((x,y)\in \mathfrak{L}^2[0,1]\times\mathbb{R}\). Let \(L_{\phi}=\left(\int_{0}^{1}|\sum_{j=1}^{n} \beta_j {\boldsymbol g}_j(t)|^2\mathrm{d}t\right)^{1/2}\).
\end{enumerate}
\end{example}

\subsection{Applications to nonlinear regression loss functions}
Next, we move on to study nonlinear regression problems. 
\begin{proposition}[Nonlinear regression loss] \label{prop:nonlin_loss} 
Let \(\mathcal{Z}\) be a (finite or infinite dimensional)  real vector space. Suppose that  \(\left\llbracket \cdot\right\rrbracket \colon\mathcal{Z}\rightarrow[0,\infty]\) is absolutely homogeneous and proper, \(\phi\colon \mathcal{Z}\rightarrow\mathbb{R} \) is linear, \(\left\llbracket \cdot\right\rrbracket^{-1}(0) \subseteq \phi^{-1}(0)\) and
\begin{equation*} 
	L_{\phi} \coloneqq \sup_{z\in\mathcal{Z}} \left\{ |\phi(z)| \mid \left\llbracket z\right\rrbracket  = 1  \right\} \in [0,+\infty).
\end{equation*}
In addition, given \(\delta>0\), suppose a univariate function \(h\colon\mathbb{R}\rightarrow\mathbb{R}\) satisfies the following assumptions: 	
\begin{enumerate}[label=(H\arabic*)]
	\item \(h\) is globally \(L_{h}\)-Lipschitz on $\mathbb{R}$ with $L_h>0$;
	\item for any \(t_0\in\mathbb{R}\), there exists \(\{t_k\}_{k=1}^{\infty}\) such that \(|t_k| \in [L_{\phi}\delta,\infty)\setminus \{0\} \) and 
	\[\lim_{k\rightarrow\infty} \frac{h(t_k+t_0)-h(t_0)}{|t_k|} = L_{h}.\]
\end{enumerate}
Then the function \(\psi \colon\mathcal{Z} \rightarrow\mathbb{R} \) defined by $\psi(z)=h(\phi(z))$ is $(L_{h}L_{\phi},d)$-Lipschitz at $\mathcal{Z}$, where the cost function $d:\mathcal{Z}\times \mathcal{Z}\rightarrow [0,\infty]$ is defined as \(d(z',z)\coloneqq \left\llbracket z'-z\right\rrbracket\) for any $z',z\in \mathcal{Z}$. Moreover, $\psi$ also satisfies Assumptions (A1) and (A2) at $\mathcal{Z}$ with $\delta$  if \(L_{\phi}>0\).
\end{proposition}
\textbf{Proof.}
By Proposition~\ref{prop:lin_loss}, one has that \(\phi\) is \((L_{\phi},d)\)-Lipschitz at \(\mathcal{Z}\). According to (H1), for any $z',z\in \mathcal{Z}$, we have
\begin{equation*}
\begin{array}{ll}
	\left|\psi(z')-\psi(z)\right| = \left| h(\phi(z')) - h(\phi(z)) \right| 
	\leq L_h \left| \phi(z')-\phi(z)\right|\leq L_{h}L_{\phi} d(z',z).
\end{array}
\end{equation*}
Hence, \(\psi\) is \((L_{h}L_{\phi},d) \)-Lipschitz at $\mathcal{Z}$. 

Next, suppose that \(L_{\phi}>0\). Then we have \(\psi\) satisfies Assumption (A1) at $\mathcal{Z}$ with $\delta$. For any \(\hat{z}\in\mathcal{Z}\), denote \(t_0 =\phi\left(\hat{z}\right) \), then we have \(\psi(\hat{z})=h(t_0)\). Let $0<\epsilon<L_{h}L_{\phi}$ be any scalar. As in the proof of Proposition~\ref{prop:lin_loss}, there exists \(\tilde{v}\in\mathcal{Z}\) such that \(\left\llbracket \tilde{v}\right\rrbracket  = 1 \) and $0 < L_{\phi} - \frac{\epsilon}{2L_{h}}< \phi(\tilde{v})\leq L_{\phi} $. We know from (H2) that there exists $ |\tilde{t}|\geq L_{\phi}\delta$ such that
\begin{equation*} 
	\left| \frac{h(\tilde{t}+t_0)-h(t_0)}{|\tilde{t}|} -L_{h}\right| \leq   \frac{\epsilon}{2\phi(\tilde{v})}.
\end{equation*}
Since \(h\) is globally \(L_{h}\)-Lipschitz continuous, we have \( \left| h(\tilde{t}+t_0)-h(t_0)\right|\leq L_{h} |\tilde{t}| \). Thus, we have
\begin{equation*}
\begin{array}{ll}
	\left| \frac{h(\tilde{t}+t_0)-h(t_0)}{|\tilde{t}|} -L_{h} \right| = L_{h}-\frac{h(\tilde{t}+t_0)-h(t_0)}{|\tilde{t}|}
	\leq   \frac{\epsilon}{2\phi(\tilde{v})},
\end{array}
\end{equation*}
which implies that
\begin{equation*}
	h(\tilde{t}+t_0) - h(t_0) \geq \left(L_{h}-\frac{\epsilon}{2\phi(\tilde{v})}\right)|\tilde{t}|.
\end{equation*}
Let \(\tilde{z}\coloneqq \hat{z}+ \tilde{t}\tilde{v}/\phi(\tilde{v})\). Then we have
\begin{equation*}
	d\left(\tilde{z},\hat{z}\right) = \left\llbracket \frac{\tilde{t}}
		{\phi(\tilde{v}) }\tilde{v}\right\rrbracket = \frac{|\tilde{t}|}{\phi(\tilde{v})} \geq \frac{|\tilde{t}|}{L_{\phi}} \geq \delta,
\end{equation*}
and
\begin{equation*}
\begin{array}{lllll}
	&\psi(\tilde{z}) - \psi(\hat{z}) =h \left( \phi\left(\hat{z}+ \frac{\tilde{t}\tilde{v}}{\phi(\tilde{v}) }  \right)\right) - h(\phi(\hat{z})) 
	= h(t_0+\tilde{t}) - h(t_0) \\ 
	&\geq \left(L_{h}-\frac{\epsilon}{2\phi(\tilde{v}) }\right)|\tilde{t}| = \left(L_{h} \phi(\tilde{v}) - \frac{\epsilon}{2} \right) \frac{|\tilde{t}|}{\phi(\tilde{v}) } \\
	&\geq \left(L_{h}\left(L_{\phi}-\frac{\epsilon}{2L_{h}} \right)-\frac{\epsilon}{2} \right) \frac{|\tilde{t}|}{\phi(\tilde{v})} 
	=   ( L_{\phi}L_h-\epsilon)d(\tilde{z},\hat{z}).
\end{array}
\end{equation*}
Therefore, \(\psi\) satisfies Assumption (A2) at $\mathcal{Z}$ with the given $\delta$. This completes the proof.
\hfill\(\square\) 

\begin{remark}
Proposition~\ref{prop:nonlin_loss} characterizes a certain case where \(\ell\) is nonlinear but Theorem~\ref{thm:main_r1} still holds true. In this case, the  loss function \(\ell\) can be written as a composition of a linear map \(\phi\) and a univariate \(L_{h}\)-Lipschitz function \(h\), where the Lipschitz constant \(L_{h}\) is approximately attainable on \([L_{\phi}\delta,\infty)\setminus \{0\}\). In particular, if \(L_{h}\) is approximately attainable at infinity, then Assumption (H2) holds regardless of the choice of \(\phi\), as long as \(L_{\phi} \) is finite. 
\end{remark}
	
Thanks to Proposition~\ref{prop:nonlin_loss}, Theorem~\ref{thm:main}(c) and Theorem~\ref{thm:main_r1}, we now can cover a broader class of nonlinear regression loss functions in the following corollary.

\begin{corollary}\label{coro: nonlinear}
Under the setting of Proposition~\ref{prop:nonlin_loss}, for any empirical distribution \(\mathbb{P}_N\) on \( \mathcal{Z} \), the following equivalence holds for any given $\delta>0$:
\begin{equation*}
\sup_{\mathbb{P}\in\mathfrak{M}_1 }  \mathrm{E}_{\mathbb{P}} [h(\phi(Z)) ] = \mathrm{E}_{\mathbb{P}_N}[h(\phi(Z)) ] + L_h L_{\phi}\delta,
\end{equation*}
where \(\mathfrak{M}_1\coloneqq \left\{\mathbb{P} \in \mathcal{P}(\mathcal{Z})\mid  \mathcal{W}_{d,1}(\mathbb{P},\mathbb{P}_N) \leq \delta \right\}\).
\end{corollary}
	
We give the following example as an application of Corollary \ref{coro: nonlinear} to show the equivalence of the worst case loss quantity in the WDRO problem and the regularization scheme for nonlinear regression loss functions, including the log-cosh loss, the quantile loss and Huber loss \citep{huber1973robust,koenker1978regression,koenker2001quantile,wang2020comprehensive}, whose proof can be found in Appendix \ref{sec: proof_ex44}.

\begin{example}\label{ex: logcosh}
Given any $\delta>0$, \(\beta\in\mathbb{R}^n\) and any empirical distribution \(\mathbb{P}_N\) on \( \mathbb{R}^n\times\mathbb{R}\), we have
\begin{equation*}
	\sup_{\mathbb{P}\in\mathfrak{M}_1 }  \mathrm{E}_{\mathbb{P}} [ h(Y - \langle\beta,X\rangle) ] = \mathrm{E}_{\mathbb{P}_N}[h(Y - \langle\beta,X\rangle)] +L_{\phi} \delta,
\end{equation*}
where \(\mathfrak{M}_1\coloneqq \left\{\mathbb{P}\in \mathcal{P}(\mathcal{Z}\mid  \mathcal{W}_{d,1}(\mathbb{P},\mathbb{P}_N) \leq \delta \right\}\) and \(h\) takes one of the following forms:
\begin{enumerate}[label=(\alph*)]
	\item log-cosh loss: \(h \colon t\mapsto \log(\cosh(t)) \);
	\item Huber loss: \(h\colon t\mapsto \begin{cases}
		\frac{1}{2}t^2 & \text{if } |t| \leq 1, \\
		|t| - \frac{1}{2} & \text{otherwise};
	    \end{cases}\)   
	\item quantile loss: \(h\colon t\mapsto \begin{cases}
		\gamma t & \text{if } t \geq 0,\\
		-t & \text{otherwise},
		\end{cases}\) \quad with $\gamma\in (0,1)$; 
\end{enumerate}
and \(d(\cdot,\cdot)\) and \(L_{\phi} \) take one of the following forms
\begin{enumerate}[label=(\roman*)]
	\item $d((x',y'),(x,y)) = \|[x'-x;y'-y]\|_{\mathbb{R}^{n+1}}$ and \(L_{\phi}=\|[-\beta;1]\|_{\mathbb{R}^{n+1},*} \);
	\item $d((x',y'),(x,y)) = \|x'-x\|_{\mathbb{R}^n}+{\boldsymbol \delta}_{\{0\}}(y'-y)$ and \(L_{\phi}=\|\beta\|_{\mathbb{R}^{n},*} \).
\end{enumerate}
\end{example}	

\subsection{A special regression model}
In this subsection, we introduce an interesting example wherein our results can be applied to the cases when the cost function \(d(\cdot,\cdot)\) is nonconvex, not positive definite and the weak Lipschitz constant is not in the popular form of the norm of the regression vector \(\beta\). In \citet[Remark 19]{shafieezadeh2019regularization}, it has been pointed out that the Tikhonov-regularized problem with respect to a Lipschitz loss function can not be explained by a distributionally robust learning problem. {In addition, it has been shown in \cite{li2022tikhonov} that the Tikhonov-regularized problem with respect to the squared loss has an equivalence interpretation as a martingale DRO problem with respect to the quadratic cost \(d(z',z)= \| z'-z\|^2 \).} Nevertheless, according to Example~\ref{ex:linear-square}, whose proof is in Appendix \ref{appen:proof_special_reg}, we show that with the notion of the weak Lipschitz property in Definition~\ref{def:Lip}, the Tikhonov-regularized squared loss problem is equivalent to a {WDRO problem} swith a specifically designed cost function \(d(\cdot,\cdot)\) on the sample space. {In particular, unlike \(d(z',z)= \| z'-z\|^2 \), our designed cost function is nonconvex and not positive definite.}
	
\begin{example}[Ridge linear ordinary regression \citep{horel1962application,hoerl1970ridge}] \label{ex:linear-square}
For any $z'=(x',y'),z=(x,y)\in \mathbb{R}^n\times\mathbb{R}$, define $d(z',z) = \|z'-z\|_2\|z'+z\|_2$. Given any $\delta>0$, \(\beta\in\mathbb{R}^n\) and any empirical distribution \(\mathbb{P}_N\) on \( \mathbb{R}^n\times\mathbb{R}\). For \(\mathfrak{M}_1\coloneqq \left\{\mathbb{P}\in \mathcal{P}(\mathcal{Z})\mid  \mathcal{W}_{d,1}(\mathbb{P},\mathbb{P}_N) \leq \delta \right\}\), we have
\begin{equation*}
\sup_{\mathbb{P}\in\mathfrak{M}_1}    \mathrm{E}_{\mathbb{P}} [(Y   +   \langle\beta,   X\rangle)^2 ]  =   \mathrm{E}_{\mathbb{P}_N}[(Y   +   \langle\beta,  X\rangle)^2] + \|\beta\|_{2}^2\delta + \delta.
\end{equation*}
\end{example}

\subsection{Applications to classification loss functions}
Next we turn our attention to  linear-type classification loss functions. For the detailed proof of the following proposition, see Appendix~\ref{sec:proof_prop_lin_loss2}.
\begin{proposition}[Linear classification loss]
\label{prop:lin_loss2}
Let \(\mathcal{Z}=\mathcal{X}\times \{-1,1\}\), where \(\mathcal{X}\) is a (finite or infinite dimensional) real vector space. Suppose that  \(\left\llbracket\cdot\right\rrbracket\colon\mathcal{X}\rightarrow[0,\infty]\) is absolutely homogeneous and proper, \(\phi\colon \mathcal{X}\rightarrow\mathbb{R} \) is linear, \(\left\llbracket \cdot\right\rrbracket ^{-1}(0) \subseteq \phi^{-1}(0)\), and
\begin{equation*} 
L_{\phi} \coloneqq \sup_{x\in\mathcal{X}} \left\{ |\phi(x)| \mid \left\llbracket x\right\rrbracket = 1  \right\} \in [0,+\infty).
\end{equation*}
Let the cost function $d:\mathcal{Z}\times \mathcal{Z}\rightarrow [0,\infty]$ be defined as \(d(z',z)\coloneqq \left\llbracket x'-x\right\rrbracket   + {\boldsymbol \delta}_{\{0\}}(y'-y) \) and the function $\psi:\mathcal{Z}\rightarrow \mathbb{R}$ be defined  as $\psi(z)=y\cdot \phi(x)$  for any $z'=(x',y'),z=(x,y)\in \mathcal{Z}$. Then for any $\tau\in \mathbb{R}$, the functions \(\left|\tau-\psi\right|\) and \(\max\{0,\tau-\psi\}\) are $(L_{\phi},d)$-Lipschitz at $\mathcal{Z}$. Furthermore, they also satisfy Assumptions (A1), (A2) and (B) at $\mathcal{Z}$ for any $\delta>0$  if \(L_{\phi}>0\).
\end{proposition}

Together with Theorem~\ref{thm:main}(c), Theorem~\ref{thm:main_r1} and Theorem~\ref{thm:main_r2}, we have the following corollary on the equivalence  \eqref{eq:main_equation} for linear-type classification loss functions.
	
\begin{corollary}\label{coro: linear2}
Under the setting of Proposition \ref{prop:lin_loss2},  given any scalars $r\geq 1$, \(\delta>0 \) and any empirical distribution \(\mathbb{P}_N\) on \( \mathcal{Z} \). We have that for any $\tau\in \mathbb{R}$,
\begin{equation*}
\sup_{\mathbb{P}\colon  \mathcal{W}_{d,r}(\mathbb{P},\mathbb{P}_N) \leq \delta }  \mathrm{E}_{\mathbb{P}} [\ell(Z)]  = \left(\left(\mathrm{E}_{\mathbb{P}_N}[\ell(Z)]\right)^{\frac{1}{r}} + L_{\phi}\delta\right)^{r},
\end{equation*}
where \(\ell(Z)=\left|\tau-\phi (Z)\right|^r \) or \(\ell(Z)= \max\{0,\tau-\phi (Z)\}^r \).
\end{corollary}

Now we give an example of the higher-order hinge loss binary classification and the higher-order support vector machine classification, which is an application of the above corollary. The detailed proof can be found in Appendix \ref{sec: proof_ex43}. Note that the second conclusion in \citet[Theorem 2]{blanchet2019robust} can be viewed as a special case of this example with \(r=1\).
	
\begin{example}[Higher-order hinge loss /Higher-order support vector machine] \label{eg:hinge-class}
For any $(x',y'),(x,y)\in \mathbb{R}^n\times \{-1,1\}$, define the cost function $d((x',y'),(x,y))=\|x'-x\|_{\mathbb{R}^n}+ {\boldsymbol \delta}_{\{0\}}(y'-y)$. Given any $\delta>0$, $\beta\in \mathbb{R}^n$ and any empirical distribution $\mathbb{P}_N$ on $\mathbb{R}^n\times \mathbb{R}$, we have
\begin{equation*}
\sup_{\mathbb{P}\in\mathfrak{M}_r }  \mathrm{E}_{\mathbb{P}} [\ell(Z;\beta)]  = \left(\left(\mathrm{E}_{\mathbb{P}_N}[\ell(Z;\beta)]\right)^{\frac{1}{r}} + \|\beta\|_{\mathbb{R}^n,*}\delta\right)^{r},
\end{equation*}
where \(\mathfrak{M}_r\coloneqq \left\{\mathbb{P}\in \mathcal{P}(\mathcal{Z})\mid  \mathcal{W}_{d,r}(\mathbb{P},\mathbb{P}_N) \leq \delta \right\}\) and \(\ell(Z;\beta)\) takes one of the following forms:
\begin{enumerate}[label=(\alph*)]
	\item \(\ell(Z;\beta)= \left| 1-Y\cdot\langle \beta,X\rangle \right|^r \);
	\item \(\ell(Z;\beta)= \left( 1-Y\cdot\langle\beta,X\rangle \right)_{+}^r\).
\end{enumerate}
\end{example}

Then, we consider the applications of our results to nonlinear classification loss functions in the following proposition, whose proof can be found in Appendix \ref{sec: prof_nonlin_loss2}.
\begin{proposition}[Nonlinear classification loss] \label{prop:nonlin_loss2} 
Let \(\mathcal{Z}=\mathcal{X}\times \{-1,1\}\), where \(\mathcal{X}\) is a (finite or infinite dimensional) real vector space. Suppose that  \(\left\llbracket\cdot\right\rrbracket \colon\mathcal{X}\rightarrow[0,\infty]\) is absolutely homogeneous and proper, \(\phi\colon \mathcal{X}\rightarrow\mathbb{R} \) is linear, \(\left\llbracket \cdot\right\rrbracket ^{-1}(0) \subseteq \phi^{-1}(0)\), and 
\begin{equation*} 
L_{\phi} \coloneqq \sup_{x\in\mathcal{X}} \left\{ |\phi(x)| \mid \left\llbracket x\right\rrbracket  = 1  \right\} \in [0,+\infty).
\end{equation*}
Let the cost function $d:\mathcal{Z}\times \mathcal{Z}\rightarrow [0,\infty]$ be defined as \(d(z',z)\coloneqq \left\llbracket x'-x\right\rrbracket   + {\boldsymbol \delta}_{\{0\}}(y'-y) \) and the function $\psi:\mathcal{Z}\rightarrow \mathbb{R}$ be defined  as $\psi(z)=h(y\cdot \phi(x))$  for any $z'=(x',y'),z=(x,y)\in \mathcal{Z}$. Given $\delta>0$ and suppose \(h\) satisfies Assumptions (H1-H2) in Proposition~\ref{prop:nonlin_loss}, then $\psi$ is $(L_h L_{\phi},d)$-Lipschitz at $\mathcal{Z}$. Moreover, $\psi$ also satisfies Assumptions (A1) and (A2) at $\mathcal{Z}$ with $\delta$ if \(L_{\phi}>0\).
\end{proposition}

The following corollary allows us to establish the equivalence of the worst-case loss quantity in the WDRO problem and the regularization scheme for nonlinear classification loss functions, based on Proposition~\ref{prop:nonlin_loss2}, Theorem~\ref{thm:main}(c) and Theorem~\ref{thm:main_r1}.
	
\begin{corollary}\label{coro: nonlinear2}
Under the setting of Proposition~\ref{prop:nonlin_loss2}, for any empirical distribution \(\mathbb{P}_N\) on \( \mathcal{Z} \), the following equivalence holds for any given $\delta>0$:
\begin{equation*}
\sup_{\mathbb{P}\in\mathfrak{M}_1 }  \mathrm{E}_{\mathbb{P}} [h(Y\cdot \phi(X)) ] = \mathrm{E}_{\mathbb{P}_N}[h(Y\cdot \phi(X)) ] + L_h L_{\phi}\delta.
\end{equation*}
where \(\mathfrak{M}_1\coloneqq \left\{\mathbb{P}\in \mathcal{P}(\mathcal{Z})\mid  \mathcal{W}_{d,1}(\mathbb{P},\mathbb{P}_N) \leq \delta \right\}\).
\end{corollary}

Nonlinear classification loss functions have gained widespread popularity due to their efficacy in handling real-world datasets where the relationships between the variables are intricate and nonlinear. We give the following example to show that our results are applicable in many popular instances, whose proof can be found in Appendix \ref{sec: proof_ex45}. Here, the log-exponential loss example is equivalent to the first conclusion in \citet[Theorem 2]{blanchet2019robust} and a tractable reformulation of the worst-case loss quantity with respect to the smooth hinge loss function has also been studied in \citet[Corollary 16]{shafieezadeh2019regularization}. 
	
\begin{example} \label{exam:nonlin-class}
For any $(x',y'),(x,y)\in \mathbb{R}^n\times\{-1,1\} $, define the cost function $d((x',y'),(x,y))=\|x'-x\|_{\mathbb{R}^n}+ {\boldsymbol \delta}_{\{0\}}(y'-y)$. Given any $\delta>0$, $\beta\in \mathbb{R}^n$ and any empirical distribution $\mathbb{P}_N$ on $\mathbb{R}^n\times \{-1,1\}$. We have
\begin{equation*}
\sup_{\mathbb{P}\in\mathfrak{M}_1 }  \mathrm{E}_{\mathbb{P}} [h(Y \cdot \langle \beta,X\rangle )]  =  \mathrm{E}_{\mathbb{P}_N}[h(Y \cdot \langle \beta,X\rangle )] + \|\beta\|_{\mathbb{R}^n,*} \delta,
\end{equation*}
where \(\mathfrak{M}_1\coloneqq \left\{\mathbb{P}\in \mathcal{P}(\mathcal{Z})\mid  \mathcal{W}_{d,1}(\mathbb{P},\mathbb{P}_N) \leq \delta \right\}\) and \(h\) takes one of the following forms.
\begin{enumerate}[label=(\alph*)]
	\item Log-exponential loss: \(h\colon t\mapsto \log (1+\exp(-t)) \);
	\item Smooth hinge loss: \(h \colon t\mapsto \begin{cases}
		0 & \text{if } t \geq 1,\\
		\frac{1}{2}(1-t)^2 & \text{if } 0<t<1,\\
		\frac{1}{2}-t & \text{otherwise};
		\end{cases}\)
	\item Truncated pinball loss \citep{shen2017support}: \(h\colon t\mapsto  \begin{cases}
		1-t & \text{if } t \leq 1,\\
		\tau_1 (t-1) & \text{if } 1 < t < \tau_2+1,\\
		\tau_1 \tau_2 & \text{otherwise},
	   \end{cases} \)\\ where \(\tau_1 \in [0,1],\tau_2 \geq 0 \) are two given constants.
\end{enumerate}
\end{example}

\section{Generalization to risk measure}
\label{sec: generalrisk}
In this section, we shall generalize the expectation function in the equivalence \eqref{eq:main_equation} to a risk measure, which might be nonlinear in distribution. 
	
Given any cost function $d(\cdot,\cdot)$ and any scalar $r\geq 1$, it can be easily seen that \(\mathcal{W}_{d,r}(\mathbb{P},\mathbb{P})=0 \) for any \(\mathbb{P}\in\mathcal{P}(\mathcal{Z})\) by choosing \(\pi = \mathbb{P}\otimes\mathbb{P} \). In addition, we use the convention that \(\mathcal{W}_{d,r}(\mathbb{P},\mathbb{Q}) = \infty \) when \(\Pi(\mathbb{P},\mathbb{Q}) = \emptyset \). Then we can see that \(\mathcal{W}_{d,r}(\cdot,\cdot)\) is a cost function on \(\mathcal{P}(\mathcal{Z})\). Moreover, one can define the weak Lipschitz property on \(\mathcal{P}(\mathcal{Z})\) with respect to the cost function \(\mathcal{W}_{d,r}(\cdot,\cdot)\). Inspired by \citet{wu2022generalization}, we propose the following sufficient condition under which a supremum on the Wasserstein neighborhood of \(\mathbb{P}_N\) and an infimum on \(\mathbb{R}\) are interchangeable. 

\begin{theorem} \label{thm:exchange}
Given any empirical distribution \(\mathbb{P}_N\) on \( \mathcal{P}(\mathcal{Z})\) and any given $r\geq 1$. Suppose that a function  \(\mathcal{F}\colon \mathcal{P}(\mathcal{Z})\times \mathbb{R} \rightarrow \mathbb{R} \) satisfies the following conditions:
\begin{itemize}
	\item \(\mathcal{F}\) on $\mathcal{P}(\mathcal{Z}) \times \mathbb{R}$ is concavelike in \(\mathcal{P}(\mathcal{Z})\), i.e., for any $\mathbb{P}_1,\mathbb{P}_2\in\mathcal{P}(\mathcal{Z}) $ and $0\leq \nu \leq 1$, there exists $\mathbb{P}\in \mathcal{P}(\mathcal{Z})$ such that for all \(t\in \mathbb{R}\),
	\begin{equation*}
	\nu \mathcal{F}(\mathbb{P}_1,t)+(1-\nu ) \mathcal{F}(\mathbb{P}_2,t)\leq \mathcal{F}(\mathbb{P},t);
	\end{equation*}
	\item \(\mathcal{F}(\mathbb{P},\cdot)\) is convex, coercive  and lower semi-continuous for each $\mathbb{P}\in \mathcal{P}(\mathcal{Z})$;
	\item \(\mathcal{F}(\cdot,t) \) is \( \left(L_{\mathcal{F}},\mathcal{W}_{d,r}\right) \)-Lipschitz at \(\{\mathbb{P}_N\} \) for some \(L_{\mathcal{F}}\in(0,\infty) \) which does not depend on \(t\).
\end{itemize}
Then we have that for any \(\delta>0\),
\begin{equation*}
\sup_{\mathbb{P}\colon \mathcal{W}_{d,r}(\mathbb{P},\mathbb{P}_N) \leq \delta  } \inf_{t\in\mathbb{R}} \mathcal{F}(\mathbb{P},t) =   \inf_{t\in\mathbb{R}}\sup_{\mathbb{P}\colon \mathcal{W}_{d,r}(\mathbb{P},\mathbb{P}_N) \leq \delta  } \mathcal{F}(\mathbb{P},t).
\end{equation*}
\end{theorem}
\textbf{Proof.}
Since \(\mathcal{F}(\cdot,t) \) is \( \left(L_{\mathcal{F}},\mathcal{W}_{d,r}\right) \)-Lipschitz at \(\{\mathbb{P}_N\} \), for any \(\tilde{\mathbb{P}}\in\mathcal{P}(\mathcal{Z}) \) such that \(\mathcal{W}_{d,r}(\tilde{\mathbb{P}},\mathbb{P}_N)\leq\delta \), we have for any \(t\in\mathbb{R}\),
\begin{align}
\left|\mathcal{F}(\tilde{\mathbb{P}},t) - \mathcal{F}(\mathbb{P}_N,t) \right| \leq L_{\mathcal{F}}\mathcal{W}_{d,r}(\tilde{\mathbb{P}},\mathbb{P}_N) \leq L_{\mathcal{F}}\delta.\label{eq:proof_duo}
\end{align}
As \(\mathcal{F}(\mathbb{P}_N,\cdot)\) is convex and coercive, it admits at least one minimizer. Let \(t_N \) be a minimizer of \(\mathcal{F}(\mathbb{P}_N,\cdot) \), i.e., \(t_N \in \arg\min_{t\in\mathbb{R}} \mathcal{F}(\mathbb{P}_N,t) \). Since  \(\mathcal{F}(\mathbb{P}_N,\cdot)\) is coercive, there exists \(\Delta_N > 0 \) such that for any \(t\notin [t_N-\Delta_N, t_N+\Delta_N] \), 
\begin{equation*}
\mathcal{F}(\mathbb{P}_N,t) \geq \mathcal{F}(\mathbb{P}_N,t_N) + 3L_{\mathcal{F}}\delta.
\end{equation*}
This together with \eqref{eq:proof_duo} implies that for any \(\tilde{\mathbb{P}}\in\mathcal{P}(\mathcal{Z}) \) such that \(\mathcal{W}_{d,r}(\tilde{\mathbb{P}},\mathbb{P}_N)\leq\delta \) and any scalar \(t\notin [t_N-\Delta_N, t_N+\Delta_N] \), it holds that
\begin{align}
\mathcal{F}(\tilde{\mathbb{P}},t) \geq \mathcal{F}(\mathbb{P}_N,t) - L_{\mathcal{F}}\delta \geq \mathcal{F}(\mathbb{P}_N,t_N) + 2L_{\mathcal{F}}\delta.\label{eq: Fpt_upper_bound}
\end{align}
On the other hand, \eqref{eq:proof_duo} also implies that
\begin{align}
\mathcal{F}(\tilde{\mathbb{P}},t_N) \leq \mathcal{F}(\mathbb{P}_N,t_N) + L_{\mathcal{F}}\delta.\label{eq: Fpt_lower_bound}
\end{align}
Thus, according to \eqref{eq: Fpt_upper_bound} and \eqref{eq: Fpt_lower_bound}, we have for any \(\tilde{\mathbb{P}}\in\mathcal{P}(\mathcal{Z}) \) such that \(\mathcal{W}_{d,r}(\tilde{\mathbb{P}},\mathbb{P}_N)\leq\delta \),
\begin{equation*}
\inf_{t\in\mathbb{R}} \mathcal{F}(\tilde{\mathbb{P}},t) = \inf_{t\in[t_N-\Delta_N, t_N+\Delta_N]} \mathcal{F}(\tilde{\mathbb{P}},t).
\end{equation*}
This further means that
\begin{equation*}
\begin{array}{ll}
\sup\limits_{\mathbb{P}\colon \mathcal{W}_{d,r}(\mathbb{P},\mathbb{P}_N) \leq \delta  } \inf\limits_{t\in\mathbb{R}} \mathcal{F}(\mathbb{P},t) 
= \sup\limits_{\mathbb{P}\colon \mathcal{W}_{d,r}(\mathbb{P},\mathbb{P}_N) \leq \delta  }\inf\limits_{t\in[t_N-\Delta_N, t_N+\Delta_N]} \mathcal{F}(\mathbb{P},t).
\end{array}
\end{equation*}
By \citet[Theorem 4.2']{sion1958general}, since $\mathcal{F}(\cdot,\cdot)$ is a concave-convexlike on $\mathcal{P}(\mathcal{Z})\times [t_N-\Delta_N, t_N+\Delta_N]$ and \(\mathcal{F}(\mathbb{P},\cdot)\) is lower semi-continuous for each $\mathbb{P}\in \mathcal{P}(\mathcal{Z})$ {on the compact set $[t_N-\Delta_N, t_N+\Delta_N]$}, we have
\begin{equation*}
\begin{array}{ll}
	\sup\limits_{\mathbb{P}\colon \mathcal{W}_{d,r}(\mathbb{P},\mathbb{P}_N) \leq \delta  } \inf\limits_{t\in[t_N-\Delta_N, t_N+\Delta_N]} \mathcal{F}(\mathbb{P},t)
	=   \inf\limits_{t\in[t_N-\Delta_N, t_N+\Delta_N]}\sup\limits_{\mathbb{P}\colon \mathcal{W}_{d,r}(\mathbb{P},\mathbb{P}_N) \leq \delta  } \mathcal{F}(\mathbb{P},t).
\end{array}
\end{equation*}
Therefore, it holds that
\begin{equation*}
\begin{array}{lll}
	\sup\limits_{\mathbb{P}\colon \mathcal{W}_{d,r}(\mathbb{P},\mathbb{P}_N) \leq \delta  } \inf\limits_{t\in\mathbb{R}} \mathcal{F}(\mathbb{P},t)  
	=   \inf\limits_{t\in[t_N-\Delta_N, t_N+\Delta_N]}\sup\limits_{\mathbb{P}\colon \mathcal{W}_{d,r}(\mathbb{P},\mathbb{P}_N) \leq \delta  } \mathcal{F}(\mathbb{P},t) 
	\geq \inf\limits_{t\in\mathbb{R}}\sup\limits_{\mathbb{P}\colon \mathcal{W}_{d,r}(\mathbb{P},\mathbb{P}_N) \leq \delta  } \mathcal{F}(\mathbb{P},t) .
\end{array}
\end{equation*}
As it is obvious that
\begin{equation*}
\sup_{\mathbb{P}\colon \mathcal{W}_{d,r}(\mathbb{P},\mathbb{P}_N) \leq \delta  } \inf_{t\in\mathbb{R}} \mathcal{F}(\mathbb{P},t) \leq \inf_{t\in\mathbb{R}}\sup_{\mathbb{P}\colon \mathcal{W}_{d,r}(\mathbb{P},\mathbb{P}_N) \leq \delta  } \mathcal{F}(\mathbb{P},t),
\end{equation*}
we complete the proof.
\hfill\(\square\)

Based on the above theorem, we give the following corollary, which generalizes the expectation function in the equivalence \eqref{eq:main_equation} to a risk measure. 
	
\begin{corollary}\label{coro: r1}
Given any scalar $\alpha\in (0,1)$. Let \(d(\cdot,\cdot)\) be a cost function on \(\mathcal{Z}\times\mathcal{Z}\). Suppose the function $G:\mathcal{Z}\rightarrow \mathbb{R}$ is $(L_G,d)$-Lipschitz at $\mathcal{Z}$ with $L_G\in (0,\infty)$. Let \(\mathcal{Z}_N\coloneqq\{Z^{(1)},\dots,Z^{(N)}\}\subset\mathcal{Z} \) be a given dataset and \(\mathbb{P}_N\coloneqq\sum_{i=1}^{N}\mu_i{\boldsymbol \chi}_{\{Z^{(i)}\}}\in\mathcal{P}(\mathcal{Z})\) be the corresponding empirical distribution.  Then we have the following conclusions.
		
\begin{enumerate}[label=(\alph*)] 
	\item Denote \(\mathfrak{M}_1\coloneqq \left\{\mathbb{P}\in \mathcal{P}(\mathcal{Z})\mid  \mathcal{W}_{d,1}(\mathbb{P},\mathbb{P}_N) \leq \delta \right\}\). Given $\delta>0$, suppose that for all \(t\in \mathbb{R}\),
	\begin{equation*}
	\sup_{\mathbb{P}\in\mathfrak{M}_1 }  \mathrm{E}_{\mathbb{P}} [\left( G(Z) -t\right)_{+} ] = \mathrm{E}_{\mathbb{P}_N}[\left( G(Z) -t\right)_{+} ] + L_{G}\delta.
	\end{equation*}
	Then it holds that
	\begin{equation*}
	\sup_{\mathbb{P}\in\mathfrak{M}_1 }\mathrm{CVaR}_{\alpha}^{\mathbb{P}}( G(Z)) = \mathrm{CVaR}_{\alpha}^{\mathbb{P}_N}( G(Z)) + \frac{1}{1-\alpha} L_{G}\delta.
	\end{equation*}
	\item Given \(r\geq 1\), let \(\mathfrak{M}_r\coloneqq \left\{\mathbb{P}\in \mathcal{P}(\mathcal{Z})\mid  \mathcal{W}_{d,r}(\mathbb{P},\mathbb{P}_N) \leq \delta \right\}\). Given $\delta>0$, suppose that for any $t\in \mathbb{R}$,
	\begin{equation*}
	\sup_{\mathbb{P}\in\mathfrak{M}_r }   \mathrm{E}_{\mathbb{P}} [\left( G(Z)  - t\right)_{+}^r ]  =   \left[ \left( \mathrm{E}_{\mathbb{P}_N} [\left( G(Z)  -  t\right)_{+}^r ]\right)^{\frac{1}{r}}   +  L_{G}\delta\right]^{r}.
	\end{equation*}
	Then it holds that
	\begin{equation*}
	\begin{array}{ll}
	\sup\limits_{\mathbb{P}\in\mathfrak{M}_r }
	\inf\limits_{t\in \mathbb{R}} \Big\{ t+\frac{1}{1-\alpha}
	\left(\mathrm{E}_{\mathbb{P}} [\left( G(Z) -t\right)_{+}^r ]\right)^{\frac{1}{r}}\Big\}
	= \inf\limits_{t\in\mathbb{R}} \Big\{ t+\frac{1}{1-\alpha}
	\left( \mathrm{E}_{\mathbb{P}_N}[\left( G(Z) -t\right)_{+}^r ]\right)^{\frac{1}{r}}\Big\} + \frac{1}{1-\alpha} L_{G}\delta.
	\end{array}
	\end{equation*}
\end{enumerate}
\end{corollary}
\textbf{Proof.}
(a) Define the function  \(\mathcal{F}\colon \mathcal{P}(\mathcal{Z})\times \mathbb{R} \rightarrow \mathbb{R} \) by
\begin{equation*}
	\mathcal{F}(\mathbb{P},t) \coloneqq  t  + \frac{1}{1-\alpha}\mathrm{E}_{\mathbb{P}}[\left( G(Z) -t\right)_{+}].
\end{equation*} 
By the definition of $\mathrm{CVaR}_{\alpha}^{\mathbb{P}}(\cdot)$, we have that
\begin{equation*}
	\mathrm{CVaR}_{\alpha}^{\mathbb{P}}( G(Z)) = \inf_{t\in\mathbb{R}} \mathcal{F}(\mathbb{P},t).
\end{equation*}

Now we verify that the function $\mathcal{F}(\cdot,\cdot)$ satisfies the three conditions in Theorem \ref{thm:exchange}.  It is easy to verify that the first two conditions hold. Next we turn to the third condition. Fix $t\in \mathbb{R}$ and any \(\tilde{\mathbb{P}} \in\mathcal{P}(\mathcal{Z}) \) such that  \(\mathcal{W}_{d,1}(\tilde{\mathbb{P}},\mathbb{P}_N)<\infty \). For any $\tilde{\pi}\in \Pi(\tilde{\mathbb{P}} ,\mathbb{P}_N)$, we have 
%\begin{equation*}
%\begin{array}{lllllll}
%	&\left|\mathcal{F}(\tilde{\mathbb{P}},t) - \mathcal{F}(\mathbb{P}_N,t) \right| \\
%	&=\frac{1}{1-\alpha}\left| \mathrm{E}_{\tilde{\mathbb{P}}}[\left( G(Z) -t\right)_{+}]-\mathrm{E}_{\mathbb{P}_N}[\left( G(Z) -t\right)_{+}]\right|
%	\\
%	&=  \frac{1}{1-\alpha} \Big| \int_{\mathcal{Z}\times \mathcal{Z}} \left( G(z') -t\right)_{+} \mathrm{d} \tilde{\pi}(z',z) \\
%	&\qquad - \int_{\mathcal{Z}\times \mathcal{Z}} \left( G(z) -t\right)_{+} \mathrm{d} \tilde{\pi}(z',z)\Big|\\
%	&\leq \frac{1}{1-\alpha} \int_{\mathcal{Z}\times \mathcal{Z}} \left| \left( G(z') -t\right)_{+} - \left( G(z) -t\right)_{+}\right| \mathrm{d} \tilde{\pi}(z',z)\\
%	&\leq \frac{1}{1-\alpha} \int_{\mathcal{Z}\times \mathcal{Z}} \left|  G(z') - G(z) \right| \mathrm{d} \tilde{\pi}(z',z)\\
%	& \leq \frac{L_G }{1-\alpha} \int_{\mathcal{Z}\times \mathcal{Z}} d(z,z') \mathrm{d} \tilde{\pi}(z',z).
%\end{array}
%\end{equation*}
\begin{equation*}
	\begin{array}{lllllll}
		\left|\mathcal{F}(\tilde{\mathbb{P}},t) - \mathcal{F}(\mathbb{P}_N,t) \right| 
		&=\frac{1}{1-\alpha}\left| \mathrm{E}_{\tilde{\mathbb{P}}}[\left( G(Z) -t\right)_{+}]-\mathrm{E}_{\mathbb{P}_N}[\left( G(Z) -t\right)_{+}]\right|\\
		&=  \frac{1}{1-\alpha} \Big| \int\limits_{\mathcal{Z}\times \mathcal{Z}} \left( G(z') -t\right)_{+} \mathrm{d} \tilde{\pi}(z',z) 
		- \int\limits_{\mathcal{Z}\times \mathcal{Z}} \left( G(z) -t\right)_{+} \mathrm{d} \tilde{\pi}(z',z)\Big|\\
		&\leq \frac{1}{1-\alpha} \int\limits_{\mathcal{Z}\times \mathcal{Z}} \left| \left( G(z') -t\right)_{+} - \left( G(z) -t\right)_{+}\right| \mathrm{d} \tilde{\pi}(z',z)\\
		&\leq \frac{1}{1-\alpha} \int\limits_{\mathcal{Z}\times \mathcal{Z}} \left|  G(z') - G(z) \right| \mathrm{d} \tilde{\pi}(z',z)\\
		& \leq \frac{L_G }{1-\alpha} \int\limits_{\mathcal{Z}\times \mathcal{Z}} d(z,z') \mathrm{d} \tilde{\pi}(z',z).
	\end{array}
\end{equation*}
By taking infimum over all $\tilde{\pi}\in \Pi(\tilde{\mathbb{P}} ,\mathbb{P}_N)$, we can see that
\begin{equation*}
	\left|\mathcal{F}(\tilde{\mathbb{P}},t) - \mathcal{F}(\mathbb{P}_N,t) \right| \leq \frac{L_G }{1-\alpha} \mathcal{W}_{d,1}(\tilde{\mathbb{P}},\mathbb{P}_N),
\end{equation*}
which means that \(\mathcal{F}(\cdot,t) \) is \( \left(\frac{L_G }{1-\alpha},\mathcal{W}_{d,1}\right) \)-Lipschitz at \(\{\mathbb{P}_N\} \). Then according to Theorem \ref{thm:exchange}, for any $\delta>0$, we have that
\begin{equation*}
\begin{array}{ll}
	\sup\limits_{\mathbb{P}\in \mathfrak{M}_1} \mathrm{CVaR}_{\alpha}^{\mathbb{P}}( G(Z)) &= \sup\limits_{\mathbb{P}\in \mathfrak{M}_1 } \inf\limits_{t\in\mathbb{R}} \mathcal{F}(\mathbb{P},t) 
	 = \inf\limits_{t\in\mathbb{R}} \sup\limits_{\mathbb{P}\in \mathfrak{M}_1 }  \mathcal{F}(\mathbb{P},t).
\end{array}
\end{equation*}
Thus, it holds that for any $\delta>0$,
\begin{equation*}
\begin{array}{lllll}
	\sup\limits_{\mathbb{P}\in \mathfrak{M}_1 } \mathrm{CVaR}_{\alpha}^{\mathbb{P}}( G(Z)) 
	&= \inf\limits_{t\in\mathbb{R}} \sup\limits_{\mathbb{P}\in \mathfrak{M}_1 } \Big\{  t  + \frac{1}{1-\alpha}\mathrm{E}_{\mathbb{P}}[\left( G(Z) -t\right)_{+}]
	\Big\}\\
	&= \inf\limits_{t\in\mathbb{R}} \Big\{  t  + \frac{1}{1-\alpha} \sup\limits_{\mathbb{P}\in \mathfrak{M}_1 }  \mathrm{E}_{\mathbb{P}}[\left( G(Z) -t\right)_{+}]
	\Big\}\\
	&=\inf\limits_{t\in\mathbb{R}} \Big\{  t  + \frac{1}{1-\alpha} \mathrm{E}_{\mathbb{P}_N}[\left( G(Z) -t\right)_{+} ] +  \frac{1}{1-\alpha}L_{G}\delta
	\Big\}\\
	&= \mathrm{CVaR}_{\alpha}^{\mathbb{P}_N}( G(Z)) + \frac{1}{1-\alpha} L_{G}\delta,
\end{array}
\end{equation*}
where the third equality following from the fact that for all \(t\in \mathbb{R}\),
\begin{equation*}
\sup_{\mathbb{P}\in\mathfrak{M}_1 }  \mathrm{E}_{\mathbb{P}} [\left( G(Z) -t\right)_{+} ] = \mathrm{E}_{\mathbb{P}_N}[\left( G(Z) -t\right)_{+} ] + L_{G}\delta.
\end{equation*}

(b) The proof of this part is similar to the one for part (a). We omit the details here.
\hfill\(\square\)

As applications of the above corollary, we give the following three examples, stating the equivalence between the worst case loss quantity in the WDRO problem and the regularization scheme for $\nu$-support vector regression, $\nu$-support vector machine, and higher moment coherent risk measures. The proof of these examples can be found in Appendices \ref{sec: proof_ex51}-\ref{sec: proof_ex53}.
	
\begin{example}[$\nu$-support vector regression \citep{scholkopf1998support}]  
\label{ex:nusvr}
For any $(x',y'),(x,y)\in \mathbb{R}^n\times\mathbb{R}$, define the cost function $d((x',y'),(x,y)) = \|(x',y')-(x,y)\|_{\mathbb{R}^{n+1}}$. Given $\alpha\in (0,1)$, $\delta>0$, \(\beta\in\mathbb{R}^n\) and any empirical distribution \(\mathbb{P}_N\) on \( \mathbb{R}^n\times\mathbb{R}\), we have 
\begin{equation*}
	\begin{array}{ll}
	\sup_{\mathbb{P}\colon \mathcal{W}_{d,1}(\mathbb{P},\mathbb{P}_N) \leq \delta } \mathrm{CVaR}_{\alpha}^{\mathbb{P}}(\left|Y - \langle \beta,X\rangle \right|) 
	= \mathrm{CVaR}_{\alpha}^{\mathbb{P}_N}(\left|Y - \langle\beta,X\rangle \right|)  + \frac{1}{1-\alpha} \|[-\beta;1]\|_{\mathbb{R}^{n+1},*}\delta.
	\end{array}
\end{equation*}
\end{example}
	
\begin{example}[$\nu$-support vector machine \citep{scholkopf2000new}] 
\label{ex: nusvm}
For any $(x',y'),(x,y)\in \mathbb{R}^n\times \mathbb{R}$, define the cost function $d((x',y'),(x,y))=\|x'-x\|_{\mathbb{R}^n}+ {\boldsymbol \delta}_{\{0\}}(y'-y)$. Given any $\alpha\in (0,1)$, $\delta>0$, $\beta\in \mathbb{R}^n$ and any empirical distribution $\mathbb{P}_N$ on $\mathbb{R}^n\times \mathbb{R}$, we have
\begin{equation*}
\begin{array}{ll}
	\sup_{\mathbb{P}\colon \mathcal{W}_{d,1}(\mathbb{P},\mathbb{P}_N) \leq \delta } \mathrm{CVaR}_{\alpha}^{\mathbb{P}}(-Y\cdot\langle\beta,X\rangle ) 
	= \mathrm{CVaR}_{\alpha}^{\mathbb{P}_N}(-Y\cdot\langle\beta,X\rangle)  + \frac{1}{1-\alpha} \|\beta\|_{\mathbb{R}^{n},*}\delta.
\end{array}
\end{equation*}
\end{example}

\begin{example}[Higher moment coherent risk measures \citep{krokhmal2007higher}] 
\label{ex: highmoment_crm}
For any $z',z\in \mathbb{R}^n$, define the cost function $d(z',z) = \|z'-z\|_{\mathbb{R}^{n}}$. Given $\alpha\in (0,1)$, $\delta>0$, $r\geq 1$, \(\beta\in\mathbb{R}^n\) and any empirical distribution \(\mathbb{P}_N\) on \( \mathbb{R}^n\), we have 
\begin{equation*}
\begin{array}{ll}
	\sup\limits_{\mathbb{P}\colon \mathcal{W}_{d,r}(\mathbb{P},\mathbb{P}_N) \leq \delta }
	\inf\limits_{t\in \mathbb{R}} \Big\{ t+\frac{1}{1-\alpha}
	\left(\mathrm{E}_{\mathbb{P}} [\left( \langle \beta,Z\rangle -t\right)_{+}^r ]\right)^{\frac{1}{r}}\Big\}
	= \inf\limits_{t\in\mathbb{R}} \Big\{ t  + \frac{1}{1-\alpha}
\left( \mathrm{E}_{\mathbb{P}_N}[\left( \langle\beta,Z\rangle -t\right)_{+}^r ]\right)^{\frac{1}{r}}\Big\} + \frac{1}{1-\alpha} \|\beta\|_{\mathbb{R}^{n},*}\delta.
	\end{array}
\end{equation*}
\end{example}

\section{Conclusion}
\label{sec: conclusion}
In this paper, we studied a variety of the Wasserstein distributionally robust optimization problems and proposed certain conditions to quantify the corresponding worst-case loss quantity. Specifically, we drew connections and established the equivalence between the worst-case loss quantity and its associated regularization scheme. Our proposed results generalized the existing results from various perspectives, particularly by relaxing the required assumptions on the loss function and the cost function. Moreover, our constructive approaches and elementary proofs directly characterized the closed forms of the approximate worst-case distributions. Extensive examples demonstrated that our theoretical results can be applied to various problems, including regression, classification and risk measure problems.
	
Following the presented results, there are some possible topics for future studies on the WDRO problems. For example, by similar arguments as in our proposed weak Lipschitz property, the notion of the growth rate in \citet[Lemma 2]{gao2023distributionally}  can be readily extended to be dependent on the cost function and its variables. On the other hand, recent works on the WDRO problems such as  \citet{blanchet2019quantifying,zhang2022simple,gao2023distributionally}, suggest that further assumptions might be required when the empirical distribution is not discrete. In summary, we hope our results can inspire more fruitful studies on the behavior of the worst-case loss quantity and the applications of the associated regularization scheme in the machine learning and operations research.

\section*{Acknowledgement}
{Meixia Lin is supported by The Singapore University of Technology and Design under MOE Tier 1 Grant SKI 2021{\_}02{\_}08. Kim-Chuan Toh is supported by the Ministry of Education, Singapore, under its Academic Research Fund Tier 3 grant call (MOE-2019-T3-1-010).}

\appendix

\section{Proof of auxiliary results}
\subsection{Proof of $\mathcal{S}\leq \mathcal{I}$} 
\label{append:upper-bounds}
We can see that
\begin{equation*}
\begin{array}{lllllllllllll}
	& \mathcal{S}  =  \sup_{Z'\sim\mathbb{P}\colon  \mathcal{W}_{d,r}(\mathbb{P},\mathbb{P}_N) \leq \delta }  \mathrm{E}_{\mathbb{P}} [\ell (Z';\beta) ] \\
	&=   \sup_{Z'\sim\mathbb{P}\in\mathcal{P}(\mathcal{Z})} \ \inf_{\rho\geq 0} \Big\{\mathrm{E}_{\mathbb{P}} [\ell (Z';\beta) ]  \\
	&\quad  + \rho \left(\delta^r -  \mathcal{W}^{r}_{d,r}(\mathbb{P},\mathbb{P}_N) \right)\Big\}\\
	&= \sup_{Z'\sim\mathbb{P}\in\mathcal{P}(\mathcal{Z})}\ \inf_{\rho\geq 0}  \Big\{\mathrm{E}_{\mathbb{P}} [\ell (Z';\beta) ]   \\
	&\quad  + \rho \left(\delta^r -  \inf_{\pi \in \Pi(\mathbb{P},\mathbb{P}_N)} \left\{\mathrm{E}_{(Z',Z)\sim\pi}\left[d^r(Z',Z)\right] \right\} \right)\Big\}\\
	& =  \sup_{\mathbb{P}\in\mathcal{P}(\mathcal{Z})} \ \inf_{\rho\geq 0} \ \sup_{\pi \in \Pi(\mathbb{P},\mathbb{P}_N)} \Big\{ \rho\delta^{r} \\
	&\quad  + \mathrm{E}_{(Z',Z)\sim\pi}\left[\ell (Z';\beta) - \rho d^r(Z',Z) \right] \Big\} \\
	& \leq \inf_{\rho\geq 0}  \ \sup_{\mathbb{P}\in\mathcal{P}(\mathcal{Z})}\ \sup_{\pi \in \Pi(\mathbb{P},\mathbb{P}_N)} \Big\{ \rho\delta^{r}  \\
	&\quad  + \mathrm{E}_{(Z',Z)\sim\pi}\left[\ell (Z';\beta) - \rho d^r(Z',Z) \right] \Big\} \\
	& \leq \inf_{\rho\geq 0}  \ \sup_{\mathbb{P}\in\mathcal{P}(\mathcal{Z}),\pi \in \Pi(\mathbb{P},\mathbb{P}_N)} \Big\{ \rho\delta^{r}  \\
	&\quad  + \mathrm{E}_{(Z',Z)\sim\pi}\left[ \sup_{z'\in\mathcal{Z}} \left\{\ell (z';\beta) - \rho d^r(z',Z)\right\} \right]\Big \} \\
	& = \inf_{\rho\geq 0}   \left\{ \rho\delta^{r}  +  \mathrm{E}_{\mathbb{P}_N}\left[ \sup_{z'\in\mathcal{Z}} \left\{\ell (z';\beta)  -  \rho d^r(z',Z)\right\} \right] \right\} \\
	& =\mathcal{I},
\end{array}
\end{equation*}
which completes the proof.

\subsection{Proof of Lemma \ref{lemma: singleton}}  
\label{append:proof_lemma1}
For any \(\pi\in\Pi(\mathbb{P},{\boldsymbol \chi}_{\{\hat{z}\}}) \), we have \(\pi(A\times\mathcal{Z}) = \mathbb{P}(A), \pi(\mathcal{Z}\times B) = {\boldsymbol \chi}_{\{\hat{z}\}}(B)  \) for any measurable sets \(A,B\subset\mathcal{Z}\). In particular, it holds that
\begin{align}
	\pi(\mathcal{Z}\times \left(\mathcal{Z}\setminus \{\hat{z} \} \right) ) = {\boldsymbol \chi}_{\{\hat{z}\}}(\mathcal{Z}\setminus \{\hat{z} \} ) =0.\label{eq:pi_on_single}
\end{align}
This implies that for any measurable set \(A\subset\mathcal{Z}\), \(\pi(A\times\left(\mathcal{Z}\setminus \{\hat{z} \} \right))  =0  \) and hence 
\begin{equation*}
\begin{array}{ll}
	\pi(A\times\{\hat{z} \}) &=\pi(A\times\mathcal{Z}) - \pi(A\times\left(\mathcal{Z}\setminus \{\hat{z} \} \right))\\
	&=\pi(A\times\mathcal{Z})  = \mathbb{P}(A).
\end{array}
\end{equation*}
Moreover, \eqref{eq:pi_on_single} also implies that
\begin{equation*}
	\int_{\mathcal{Z}\times\left(\mathcal{Z}\setminus \{\hat{z} \} \right)} d^r(z',z)\mathrm{d}\pi(z',z) = 0.
\end{equation*}
Therefore, one has that
\begin{equation*}
\begin{array}{ll}
	\int_{\mathcal{Z}\times\mathcal{Z}} d^r(z',z)\mathrm{d}\pi(z',z)&=\int_{\mathcal{Z}\times\{\hat{z} \} } d^r(z',z)\mathrm{d}\pi(z',z)\\
	 &= \int_{\mathcal{Z}}d^{r}(z',\hat{z})\mathrm{d}\mathbb{P}(z').
\end{array}
\end{equation*}
This completes the proof.

\subsection{Proof of Proposition~\ref{prop:lin_loss2}} \label{sec:proof_prop_lin_loss2}
For any $z,z'\in \mathcal{Z}$, we have
\begin{equation*}
\begin{array}{lll}
	&\left|\left|\tau-\psi(z')\right|-\left|\tau-\psi(z)\right|\right|\leq \left|\psi(z')-\psi(z)\right|,\\
	&\left|\max\{0,\tau-\psi(z')\}-\max\{0,\tau-\psi(z)\}\right| \\
	&\quad \leq \left|\psi(z')-\psi(z)\right|.
\end{array}
\end{equation*}
Moreover, we can see that
\begin{equation*}
\begin{array}{ll}
	\left|\psi(z')-\psi(z)\right| &= \left|y'\cdot \phi(x')-y\cdot \phi(x)\right| \\
	&= \left|\phi(y' x'-y x)\right|\leq L_{\phi} d(z',z),
\end{array}
\end{equation*}
where the last inequality holds as follows.
\begin{itemize} 
	\item If $d(z',z)=\infty$, then it holds true;
	\item If \(d(z',z)=0\), then \(y'=y\) and \(\left\llbracket x'-x\right\rrbracket=0 \). Since \(\left\llbracket \cdot\right\rrbracket^{-1}(0) \subseteq \phi^{-1}(0) \), this implies that \(|\psi(z')-\psi(z)| = \left|\phi(y' x'-y x)\right| = \left|\phi(x'-x)\right|=0 \);
	\item If \(0<d(z',z)<\infty\), then we have $y'=y$, \(\left\llbracket x'-x\right\rrbracket\ne0 \) and
	\begin{equation*}
	\begin{array}{lll}
		&\left|\phi(y' x'-y x)\right| =\left|\phi( x'-x)\right| \\
		&= \left\llbracket x'-x\right\rrbracket \left|\phi\left(\frac{x'-x}{\left\llbracket x'-x\right\rrbracket }\right)\right|  \\
		&\leq \left\llbracket x'-x\right\rrbracket L_{\phi} = L_{\phi}d(z',z).
	\end{array}
	\end{equation*}
\end{itemize}
Therefore, \(|\tau-\psi|\) and \(\max\{0,\tau-\psi\}\) are $(L_{\phi},d)$-Lipschitz at $\mathcal{Z}$.

Suppose that \(L_{\phi}>0\). Then it can be seen that \(\left|\tau-\psi\right|\) and \(\max\{0,\tau-\psi\}\) satisfy Assumption (A1) at $\mathcal{Z}$ for any $\delta >0$ thanks to previous discussions. By the definition of \(L_{\phi}\), for any $0<\epsilon <L_{\phi}$, there exists \(\tilde{v}\in\mathcal{X}\) such that \(\left\llbracket \tilde{v}\right\rrbracket  = 1 \text{ and } \phi(\tilde{v}) \geq L_{\phi}-\epsilon/2>0\). For any $z=(x,y)\in \mathcal{Z}$ and $\sigma>0$, define $\tilde{z} = (\tilde{x},\tilde{y})$ as
\begin{equation*}
	\tilde{x} = x - {\rm sgn}(\tau-y\cdot \phi(x)) y \sigma \tilde{v} ,\quad \tilde{y} = y.
\end{equation*}
Then we can see that $d(\tilde{z},z) = \left\llbracket {\rm sgn}(\tau-y\cdot \phi(x)) y \sigma \tilde{v} \right\rrbracket =\sigma $ and 
\begin{equation*}
\begin{array}{llll}
	& \left| \tau-\psi(\tilde{z})\right| - \left| \tau-\psi(z)\right| \\
	&  =  \left|\tau    -   y  \cdot   \phi(x   -    {\rm sgn}(\tau  -  y\cdot \phi(x)) y \sigma \tilde{v})\right|  -  \left|\tau  -  y\cdot \phi(x)\right|\\
	& =  \left|\tau  -  y \cdot  \phi(x)   +   {\rm sgn}(\tau  -  y\cdot \phi(x)) \sigma \phi(\tilde{v})\right|   -   \left|\tau  -  y\cdot \phi(x)\right|\\
	& \geq \sigma \phi(\tilde{v})\geq (L_{\phi}-\epsilon) d(\tilde{z},z).
\end{array}
\end{equation*}
%\begin{equation*}
%	\begin{array}{ll}
%		&\left| \tau-\psi(\tilde{z})\right| - \left| \tau-\psi(z)\right| \\
%		& = \left|\tau  - y \cdot  \phi(x -  {\rm sgn}(\tau-y\cdot \phi(x)) y \sigma \tilde{v})\right| - \left|\tau-y\cdot \phi(x)\right|\\
%		&= \left|\tau-y\cdot \phi(x) + {\rm sgn}(\tau-y\cdot \phi(x)) \sigma \phi(\tilde{v})\right| - \left|\tau-y\cdot \phi(x)\right|\\
%		& \geq \sigma \phi(\tilde{v})\geq (L_{\phi}-\epsilon) d(\tilde{z},z).
%	\end{array}
%\end{equation*}
Therefore, for any $\delta>0$ and $z\in \mathcal{Z}$, by setting $\sigma\in [\delta,\infty)$ or $\sigma \in \mathcal{D}(z)$, we can see that \(|\tau-\psi|\) satisfies both Assumption (A2) and (B) at $\mathcal{Z}$.

Finally, we are going to prove that \(\max\{0,\tau-\psi\}\) satisfy Assumptions (A2) and (B) when \(L_{\phi}>0\). For any $z\in \mathcal{Z}$ and $\delta>0$, define $\tilde{z} = (\tilde{x},\tilde{y})$ as $\tilde{y} = y$ and
\begin{equation*}
	\tilde{x}  =  \left\{
	\begin{array}{ll}
		x  -  y \delta \tilde{v} &\text{if } \tau - y\cdot \phi(x)  \geq  0 \\
		x  -  \left( 2(y\cdot \phi(x) - \tau)/\epsilon  + \delta\right) y \tilde{v} &\text{otherwise.}
	\end{array}\right.
\end{equation*}
Then if $\tau-y\cdot\phi(x) \geq 0$, we have $d(\tilde{z},z) = \left\llbracket y \delta \tilde{v}\right\rrbracket =\delta$ and
\begin{equation*}
\begin{array}{lll}
	&\max\{0,\tau-\psi(\tilde{z})\} - \max\{\tau-\psi(z)\} \\
	&= \max\{0,\tau-y\cdot\phi(\tilde{x})\} - \max\{0,\tau-y\cdot \phi(x)\}\\
	&\geq y \cdot \phi(x-\tilde{x}) = \delta \phi(\tilde{v})\geq (L_{\phi}-\epsilon) d(\tilde{z},z);
\end{array}
\end{equation*}
if $\tau-y\cdot \phi(x) < 0$, we have $d(\tilde{z},z) = \left\llbracket \left( 2(y\cdot \phi(x)-\tau)/\epsilon +\delta\right) y \tilde{v}\right\rrbracket =2(y\cdot \phi(x)-\tau)/\epsilon +\delta \geq \delta $ and
\begin{equation*}
\begin{array}{lllll}
	& \max\{0,\tau-\psi(\tilde{z})\} - \max\{\tau-\psi(z)\} \\
	&= \max\{0,\tau-\tilde{y}\cdot\phi(\tilde{x})\} - \max\{0,\tau-y\cdot \phi(x)\}\\
	&\geq \tau-y\cdot \phi(x) + \phi(\tilde{v})d(\tilde{z},z)\\
	&\geq -\frac{\epsilon}{2}d(\tilde{z},z)+ \phi(\tilde{v})d(\tilde{z},z)\\
	&\geq (L_{\phi}-\epsilon) d(\tilde{z},z).
\end{array}
\end{equation*}
This means that \(\max\{0,\tau-\psi\}\) satisfies Assumptions (A2) at $\mathcal{Z}$ for any $\delta>0$. Next, we turn to Assumption (B). Fix $r> 1$, $\delta>0$, a dataset \(\mathcal{Z}_N\subset\mathcal{Z} \) and the corresponding empirical distribution \(\mathbb{P}_N \). For any $\hat{z}\in \mathcal{Z}_N$, we consider the following cases.
\begin{itemize}
	\item If \(\mathrm{E}_{\mathbb{P}_N}[\max\left\{0, \tau-\phi(Z) \right\}^r ]=0 \), (A2) and (B) are equivalent. 
	\item If \(\mathrm{E}_{\mathbb{P}_N}[\max\left\{0, \tau-\phi(Z) \right\}^r ]\ne0 \) and \( \psi(\hat{z})\geq \tau  \), one can choose \(\tilde{z}=\hat{z}\) such that (B) holds.
	\item If \(\mathrm{E}_{\mathbb{P}_N}[\max\left\{0, \tau-\phi(Z) \right\}^r ]\ne0 \) and \( \psi(\hat{z})<\tau  \). For any \(\sigma\geq \delta\), let $\tilde{z}=(\tilde{x},\tilde{y})\in \mathcal{Z}$ be defined as \(\tilde{x} = \hat{x} -\hat{y} \sigma\tilde{v}\) and $\tilde{y} = \hat{y}$. Then we have
	\begin{equation*}
	\begin{array}{ll}
		\psi(\tilde{z}) = \tilde{y}\cdot \phi(\tilde{x}) &= \hat{y}\cdot\phi(\hat{x}) - \sigma \phi(\tilde{v}) \\
		&= \psi(\hat{z}) - \sigma \phi(\tilde{v}) <\tau.
	\end{array}
	\end{equation*}
	Moreover, one can see that \(d(\tilde{z},\hat{z}) = \left\llbracket \tilde{x}-\hat{x}\right\rrbracket  = \left\llbracket \hat{y} \sigma\tilde{v}\right\rrbracket =\sigma \) and
	\begin{equation*}
    \begin{array}{ll}
		&\max\left\{0, \tau-\psi(\tilde{z}) \right\}-\max\left\{0,\tau- \psi(\hat{z})\right\} \\
		&=  \hat{y} \cdot \phi(\hat{x} - \tilde{x})= \sigma\phi(\tilde{v}) \geq (L_{\phi}-\epsilon)d(\tilde{z},\hat{z}).
	\end{array}
	\end{equation*}
\end{itemize}
This means that \(\max\{0,\tau-\phi\}\) satisfies Assumptions (B) at $\mathcal{Z}$ for any $\delta>0$.

\subsection{Proof of Proposition \ref{prop:nonlin_loss2} }
\label{sec: prof_nonlin_loss2}
As the proof of this proposition is similar to that of Proposition \ref{prop:nonlin_loss}, we only sketch it. For any $z',z\in \mathcal{Z}$, we can see that
\begin{equation*}
\begin{array}{llll}
	&\left|\psi(z')-\psi(z)\right| = \left| h(y'\cdot\phi(x')) - h(y\cdot\phi(x)) \right| \\
	&\leq L_{h}\left|\phi(y'\cdot x' - y\cdot x)\right| \\
	&\leq L_{h}L_{\phi} 
	\left\llbracket y'\cdot x' - y\cdot x\right\rrbracket  \\
	&\leq L_{h}L_{\phi} \left(\left\llbracket x'-x\right\rrbracket  + {\boldsymbol \delta}_{\{0\}}(y'-y)\right) = L_{h}L_{\phi}  d(z',z).
\end{array}
\end{equation*}
Hence, \(\psi\) is \((L_{h}L_{\phi},d) \)-Lipschitz at $\mathcal{Z}$.

Next, suppose that \(L_{\phi}>0\). Then we can see that \(\psi\) satisfies Assumption (A1) at $\mathcal{Z}$ with $\delta$. For any \(\hat{z}=(\hat{x},\hat{y})\in\mathcal{Z}\), denote \(t_0\coloneqq \hat{y}\cdot\phi\left(\hat{x}\right) \), then we have \(\psi(\hat{z})=h(t_0)\). 
Let \(0<\epsilon<L_{h}L_{\phi}\). As in the proof of Proposition~\ref{prop:lin_loss}, there exists \(\tilde{v}\in\mathcal{V}\) such that \(\left\llbracket \tilde{v}\right\rrbracket = 1 \) and $0 < L_{\phi} - \frac{\epsilon}{2L_{h}}< \phi(\tilde{v})\leq L_{\phi} $. By Assumption (H2) and the Lipschitz property of $h$, there exists \(\tilde{t}\in\mathbb{R}\) such that \(|\tilde{t}|\geq L_{\phi}\delta\) and 
\begin{equation*}
	h(\tilde{t}+t_0) - h(t_0) \geq \left(L_{h}-\frac{\epsilon}{2\phi(\tilde{v})}\right)|\tilde{t}|.
\end{equation*}
Define \(\tilde{z}\coloneqq \left(\hat{x}+\tilde{t} \hat{y}\tilde{v}/\phi(\tilde{v}),\hat{y}\right) \), then we have 
\begin{equation*}
	d\left(\tilde{z},\hat{z}\right) = \left\llbracket \frac{\tilde{t}\hat{y}}
		{\phi(\tilde{v})}\tilde{v}\right\rrbracket =\frac{|\tilde{t}|}{\phi(\tilde{v})} \geq \frac{|\tilde{t}|}{L_{\phi}}\geq \delta,
\end{equation*}
and 
\begin{equation*}
\begin{array}{lllll}
	&\psi(\tilde{z}) - \psi(\hat{z}) = h \left(\hat{y}\cdot \phi\left( \hat{x}+ \frac{\tilde{t}\hat{y}}{\phi(\tilde{v})}\tilde{v}\right) \right) - h(\hat{y}\cdot  \phi(\hat{x}))\\
	&=
	h(t_0+\tilde{t}) - h(t_0) \\
	&\geq \left(L_{h}-\frac{\epsilon}{2\phi(\tilde{v})}\right)|\tilde{t}| 
	= \left(L_{h} \phi(\tilde{v}) - \frac{\epsilon}{2} \right) \frac{|\tilde{t}|}{\phi(\tilde{v}) } \\
	&\geq \left(L_{h}\left(L_{\phi}-\frac{\epsilon}{2L_{h}} \right)-\frac{\epsilon}{2} \right) \frac{|\tilde{t}|}{\phi(\tilde{v})} \\
	&=   ( L_{\phi}L_h-\epsilon)d(\tilde{z},\hat{z}).
\end{array}
\end{equation*}
Therefore, \(\psi\) satisfies Assumption (A2) at $\mathcal{Z}$ with the given $\delta$. This completes the proof.

\section{Proof of examples}
We start with the following technical lemma, which will be used later in the proof of Examples \ref{exam:non-Lipschitz} and \ref{exam:remove_delta}.
\begin{lemma} \label{lem:H} 
	Let \(T\in (0,\infty] \) and  \(h\colon(0,T)\rightarrow\mathbb{R} \) be a convex, continuously differentiable function.   Given any \(\hat{t}\in(0,T)\), let \(\mathcal{H}_{h,\hat{t}} \colon (0,T)\rightarrow \mathbb{R} \) defined as
	\begin{equation*}
		\mathcal{H}_{h,\hat{t}}(t) = \begin{cases}
			\frac{h(t)-h(\hat{t})}{t-\hat{t}}   & \text{if } t \ne \hat{t},\\
			\nabla h(\hat{t})  & \text{otherwise.}
		\end{cases}
	\end{equation*}
	Then we have that \(\mathcal{H}_{h,\hat{t}}\) is continuous and  non-decreasing on \((0,T)\). Moreover, 
	\begin{equation*}
		\sup_{t\in(0,T)} \left|\mathcal{H}_{h,\hat{t}}(t)\right| = \max \left\{ \left|\mathcal{H}^+_{h,\hat{t}}(0) \right|, \left|\mathcal{H}^-_{h,\hat{t}}(T)\right|  \right\}.
	\end{equation*}
	where \(\mathcal{H}^+_{h,\hat{t}}(0)\coloneqq \lim_{t\downarrow0}\mathcal{H}_{h,\hat{t}}(t)  \) and \(\mathcal{H}^-_{h,\hat{t}}(T)\coloneqq \lim_{t\uparrow T}\mathcal{H}_{h,\hat{t}}(t)  \).
\end{lemma}
\textbf{Proof.}
	Since $ h$ is continuously differentiable, we can see that \(\mathcal{H}_{h,\hat{t}}\) is continuous on \((0,T)\).  For any \(t\in(0,T)\) such that \(t\ne\hat{t}\), we can see
	\begin{equation*}
		\nabla \mathcal{H}_{h,\hat{t}}(t) = \frac{ (t-\hat{t})\nabla h(t)-h(t) + h(\hat{t}) }{(t-\hat{t})^2}.
	\end{equation*}
	Note that \(h\) is convex, hence \(h(\hat{t}) \geq h(t) + (\hat{t}-t)\nabla h(t) \) for any $t\in (0,T)$, and thus \(\nabla \mathcal{H}_{h,\hat{t}}(t)\geq  0 \) for any \(t\in(0,T)\setminus \{\hat{t}\} \). Therefore, \(\mathcal{H}_{h,\hat{t}}\) is non-decreasing on \((0,T)\), and the remaining conclusion follows.
\hfill\(\square\) 

\subsection{Proof of Example \ref{exam:non-Lipschitz}}
\label{sec: proof_ex31}
(a)	Note that \(\nabla h(t) =\log(t) - \log(1-t) \)  for any \(t\in(0,1)\), we can easily see that \(h\) is convex and continuously differentiable. Moreover, we can see that
\begin{equation*}
	\sup_{t',t\in(0,1), t' \ne t} \frac{\left| h(t')-h(t)\right|}{\left|t'-t\right|} \geq \lim_{t\rightarrow0} \left|\nabla h(t)\right| = \infty.
\end{equation*}
Thus \(h\) is not globally Lipschitz on \((0,1)\). 

(b)  For any \(z'\in(0,1) \) such that \(z'\ne \hat{z}\), we have \(\beta z', \beta \hat{z}\in (0,1) \) and
\begin{equation*}
\begin{array}{lll}
	&\left|\psi_{\beta}(z')-\psi_{\beta}(\hat{z})\right|=\left|h(\beta z') - h(\beta\hat{z})\right| \\
	&= \beta|z'-\hat{z}| \left|\frac{h(\beta z') - h(\beta\hat{z})}{\beta z' - \beta\hat{z}} \right|\\
	&\leq \beta|z'-\hat{z}| \sup_{t\in(0,1) } \left|\mathcal{H}_{h,\beta\hat{z}}(t)\right|,
\end{array}
\end{equation*} 
where \(\mathcal{H}_{h,\beta\hat{z}}\) is defined as in Lemma~\ref{lem:H}. According to Lemma~\ref{lem:H}, we have that \(\mathcal{H}_{h,\beta\hat{t}}\) is continuous and non-decreasing on \((0,1)\). Moreover, we can see that
\begin{equation*}
\begin{array}{llll}
	&\mathcal{H}^+_{h,\beta\hat{z}}(0)=
	\frac{ h(\beta\hat{z})-\lim_{t\downarrow 0}h(t)}{\beta \hat{z}}=\frac{h(\beta\hat{z})}{\beta \hat{z}}\\
	&\qquad =\frac{\beta\hat{z}\log(\beta\hat{z}) + (1-\beta\hat{z})\log(1-\beta\hat{z}) }{\beta\hat{z}}<0,\\
	&\mathcal{H}^-_{h,\beta\hat{z}}(1)=\frac{ \lim_{t\uparrow 1}h(t)-h(\beta\hat{z})}{1-\beta \hat{z}}=\frac{- h(\beta\hat{z})}{1-\beta \hat{z}}\\
	&\qquad =\frac{-\beta\hat{z}\log(\beta\hat{z}) - (1-\beta\hat{z})\log(1-\beta\hat{z}) }{1-\beta\hat{z}}>0.
\end{array}
\end{equation*}
Together with Lemma~\ref{lem:H}, we have that
\begin{equation*}
\begin{array}{ll}
	\beta \sup_{t\in(0,1) } \left|\mathcal{H}_{h,\beta\hat{z}}(t)\right| &=\beta \max \left\{ -\mathcal{H}^+_{h,\beta\hat{t}}(0) , \mathcal{H}^-_{h,\beta\hat{t}}(1) \right\}\\
	&=-\beta\mathcal{H}^+_{h,\beta\hat{z}}(0) = L_{\beta}^{\{\hat{z}\}},
\end{array}
\end{equation*}
where the second equality follows from the fact that $\hat{z}\in (0,\frac{1}{2}]$. Thus, \(\psi_{\beta}\) is \((L_{\beta}^{\{\hat{z}\}},d)\)-Lipschitz at \(\{\hat{z}\}\).

(b1) Suppose \(0 <\delta< \hat{z}\). For any \(\epsilon>0\), since \(\mathcal{H}_{h,\beta\hat{z}}\) is continuous, non-decreasing on \((0,1)\), and also satisfies \(\mathcal{H}^+_{h,\beta\hat{z}}(0) < 0\), there exists \(\tilde{z} \in (0,\hat{z}-\delta) \) such that 
\begin{equation*}
	\mathcal{H}_{h,\beta\hat{z}}(\beta\tilde{z}) < 0 \text{ and } 0 \leq \mathcal{H}_{h,\beta\hat{z}}(\beta\tilde{z}) - \mathcal{H}^+_{h,\beta\hat{z}}(0) < \frac{\epsilon}{\beta}.
\end{equation*}
Then \(d(\tilde{z},\hat{z}) = \left|\tilde{z}-\hat{z}\right| = \hat{z} - \tilde{z} > \delta  \) and
\begin{equation*}
\begin{array}{lll}
	&\psi_{\beta}(\tilde{z}) - \psi_{\beta}(\hat{z}) = -\beta (\hat{z}-\tilde{z}) \frac{h(\beta \tilde{z}) - h(\beta\hat{z})}{\beta \tilde{z} - \beta\hat{z}} \\
	&= -\beta d(\tilde{z},\hat{z})\mathcal{H}_{h,\beta\hat{z}} (\beta\tilde{z}) \\
	&\geq \beta d(\tilde{z},\hat{z}) \left(-\mathcal{H}^+_{h,\beta\hat{z}}(0) - \frac{\epsilon}{\beta} \right) =  (L_{\beta}^{\{\hat{z}\}}-\epsilon) d(\tilde{z},\hat{z}).
\end{array}
\end{equation*}
Hence, Assumption (A2) is satisfied. By Theorem~\ref{thm:main_r1}, we have that 
\begin{equation*}
	\sup_{\mathbb{P}\colon  \mathcal{W}_{d,1}(\mathbb{P},\mathbb{P}_N) \leq \delta }  \mathrm{E}_{\mathbb{P}} [\ell (Z;\beta) ]=\ell(\hat{z};\beta) + L_{\beta}^{\{\hat{z}\}} \delta.
\end{equation*}

(b2) Suppose \(\delta\geq \hat{z} \). Note that \(h\) is upper bounded by 0. Thus, we have that
\begin{equation*}
	\sup_{\mathbb{P}\colon  \mathcal{W}_{d,1}(\mathbb{P},\mathbb{P}_N) \leq \delta }  \mathrm{E}_{\mathbb{P}} [\ell (Z;\beta) ] \leq 0.
\end{equation*}
On the other hand, if we choose \(\tilde{\mathbb{P}}_k \coloneqq {\boldsymbol \chi}_{\{\tilde{z}_k\}}\) with \(\tilde{z}_k = \frac{1}{k}\) for each \(k>1\), then by Lemma~\ref{lemma: singleton}, we have \(\mathcal{W}_{d,1}(\tilde{\mathbb{P}}_k,\mathbb{P}_N) = \mathcal{W}_{d,1}({\boldsymbol\chi}_{\{\tilde{z}_k\}},{\boldsymbol\chi}_{\{\hat{z}\}}) = d(\tilde{z}_k,\hat{z}) = \hat{z}-\frac{1}{k} < \delta\) for sufficiently large $k$. Moreover, it holds that
\begin{equation*}
	\mathrm{E}_{\tilde{\mathbb{P}}_k}[\ell(Z;\beta)] = \ell(\tilde{z}_k;\beta) = h(\beta/k) \rightarrow 0
\end{equation*}
as \(k\rightarrow \infty\). Therefore, we have
\begin{equation*}
	\sup_{\mathbb{P}\colon  \mathcal{W}_{d,1}(\mathbb{P},\mathbb{P}_N) \leq \delta }  \mathrm{E}_{\mathbb{P}} [\ell (Z;\beta) ] = 0.
\end{equation*}
This completes the proof.

\subsection{Proof of Example \ref{exam:remove_delta}}
\label{sec: proof_ex32}
(a) For any $z\in \mathbb{R}$, we have $d(z,\hat{z}) = \|z\|_{{\mathbb{R}^{n}}}$ and
\begin{equation*}
\begin{array}{lll}
	& |\psi_{\beta}(z)-\psi_{\beta}(\hat{z})| = \left|h(\langle \beta,z\rangle)-\frac{1}{2}\right|\\
	&= \begin{cases}
		\frac{1}{2} & \text{if } \left|\langle\beta,z\rangle\right|\geq 1,\\
		\frac{1}{2}\left|\langle \beta,z\rangle \right| & \text{otherwise}
	\end{cases}  \\
	&\leq \frac{1}{2}\left|\langle \beta,z\rangle \right| \leq \frac{\|\beta\|_{\mathbb{R}^{n},*}}{2}\|z\|_{{\mathbb{R}^{n}}}.
\end{array}
\end{equation*}
Therefore, we can see that \(\psi_{\beta}\) is \(\left(\frac{\|\beta\|_{\mathbb{R}^{n},*}}{2},d\right)\)-Lipschitz at $\{\hat{z}\}$. 

(a1) Suppose $0<\delta\leq \frac{1}{\vartheta_2}$. For any \(0< \epsilon<\frac{\|\beta\|_{\mathbb{R}^{n},*}}{2}\), let \(\tilde{z}\coloneqq \frac{1}{\|\beta\|_{{\mathbb{R}^{n},*}}}\alpha_{\beta} \), then we have that \(d(\tilde{z},\hat{z}) = \|\tilde{z}-\hat{z}\|_{{\mathbb{R}^{n}}} = \frac{1}{\|\beta\|_{{\mathbb{R}^{n},*}}} \geq \frac{1}{\vartheta_2} \geq \delta \) and 
\begin{equation*}
\begin{array}{ll}
	&\psi_{\beta}(\tilde{z}) - \psi_{\beta}(\hat{z}) = h(1) - h\left(0\right) =\frac{1}{2} \\
	&>\left(\frac{\|\beta\|_{\mathbb{R}^{n},*}}{2}-\epsilon\right) \frac{1}{\|\beta\|_{\mathbb{R}^{n},*}} = \left(\frac{\|\beta\|_{\mathbb{R}^{n},*}}{2}-\epsilon\right)d(\tilde{z},\hat{z}),
\end{array}
\end{equation*}
which means that $\psi_{\beta}$ satisfies Assumption (A2) at $\{\hat{z}\}$ for $0<\delta\leq\frac{1}{\vartheta_2}$. Therefore, according to Theorem~\ref{thm:main_r1}, we have 
\begin{equation*}
\begin{array}{ll}
	&\sup_{\mathbb{P}\colon  \mathcal{W}_{d,1}(\mathbb{P},\mathbb{P}_N) \leq \delta }   \mathrm{E}_{\mathbb{P}} [\ell (Z;\beta) ]  =  \mathrm{E}_{\mathbb{P}_N}[\ell(Z;\beta)]   +   \frac{\|\beta\|_{\mathbb{R}^{n},*}}{2}\delta \\
	&= \ell(\hat{z};\beta)+\frac{\|\beta\|_{\mathbb{R}^{n},*}}{2}\delta,\quad 0<\delta\leq\frac{1}{\vartheta_2}.
\end{array}
\end{equation*}

(a2) Suppose \(\delta \geq \frac{1}{\vartheta_1
}\). We have that for any \(\mathbb{P}\in \mathcal{P}(\mathcal{Z})\),
\begin{equation*}
	\mathrm{E}_{\mathbb{P}}[\ell(Z;\beta)]  = \int_{\mathbb{R}^n}h(\langle\beta, z\rangle )\mathrm{d}\mathbb{P}(z)   \leq 1 \int_{\mathbb{R}^n}\mathrm{d}\mathbb{P}(z) = 1.
\end{equation*}
Hence, \(\sup_{\mathbb{P}\colon  \mathcal{W}_{d,1}(\mathbb{P},\mathbb{P}_N) \leq \delta }  \mathrm{E}_{\mathbb{P}} [\ell (Z;\beta) ] \leq 1\). On the other hand, if we choose \(\tilde{\mathbb{P}} \coloneqq {\boldsymbol \chi}_{\{\tilde{z}\}}\) with \(\tilde{z} = \frac{1}{\|\beta\|_{\mathbb{R}^{n},*}}\alpha_{\beta} \), then by Lemma~\ref{lemma: singleton}, we have \(\mathcal{W}_{d,1}(\tilde{\mathbb{P}},\mathbb{P}_N) = \mathcal{W}_{d,1}({\boldsymbol\chi}_{\{\tilde{z}\}},{\boldsymbol\chi}_{\{\hat{z}\}}) = d(\tilde{z},\hat{z}) = \frac{1}{\|\beta\|_{\mathbb{R}^{n},*}} \leq \frac{1}{\vartheta_1} \leq \delta\). Moreover, it holds that
\begin{equation*}
	\mathrm{E}_{\tilde{\mathbb{P}}}[\ell(Z;\beta)] = \ell(\tilde{z};\beta) = h(1) = 1 ,
\end{equation*}
which means
\begin{equation*}
	\sup_{\mathbb{P}\colon  \mathcal{W}_{d,1}(\mathbb{P},\mathbb{P}_N) \leq \delta }  \mathrm{E}_{\mathbb{P}} [\ell (Z;\beta) ]= 1=\ell(\hat{z};\beta) +  \frac{1}{2}.
\end{equation*}

(b) For any $z\in \mathbb{R}^n$, we can see that $d(z,\bar{z}) = \left\|z+\frac{3}{\vartheta}\alpha_{\beta} \right\|_{\mathbb{R}^n}$ and
\begin{equation*}
\begin{array}{ll}
	&|\psi_{\beta}(z)-\psi_{\beta}(\bar{z})| = \left|h(\langle \beta,z\rangle)-h(-3)\right|=\left|h(\langle \beta,z\rangle )\right|\\
	&= \begin{cases}
		0 & \text{if } \langle\beta,z\rangle \leq -1,\\
		\frac{\langle \beta,z\rangle +1}{2} & \text{if } -1 < \langle \beta,z\rangle < 1,\\
		1 & \text{otherwise},
	\end{cases} 
\end{array}
\end{equation*}
which further implies that
\begin{equation*}
\begin{array}{lll}
	&\left|\psi_{\beta}(z)-\psi_{\beta}(\bar{z})\right| \leq \frac{1}{4} \left| \langle \beta,z\rangle  +3\right| \\
	&=\frac{1}{4} \left| \langle \beta,z\rangle  +\langle \frac{3}{\vartheta}\alpha_{\beta},\beta \rangle \right|  \\
	&\leq \frac{\|\beta\|_{\mathbb{R}^n,*}}{4}\left\|z+\frac{3}{\vartheta}\alpha_{\beta} \right\|_{\mathbb{R}^n}   =   \frac{\vartheta}{4}\left\| z+\frac{3}{\vartheta}\alpha_{\beta} \right\|_{\mathbb{R}^n}  =  \frac{\vartheta}{4}d(z,\bar{z}).
\end{array}
\end{equation*} 
Therefore, we can see that \(\psi_{\beta}\) is \((\frac{\vartheta}{4},d) \)-Lipschitz at $\{\bar{z}\}$.

(b1) Suppose $0<\delta\leq \frac{4}{\vartheta}$. For any \(0<\epsilon < \frac{\vartheta}{4}\), let \(\tilde{z}\coloneqq \frac{1}{\vartheta}\alpha_{\beta} \), then \(d(\tilde{z},\bar{z}) = \|\tilde{z}-\bar{z}\|_{\mathbb{R}^n} = \frac{4}{\vartheta} \geq \delta \) and 
\begin{equation*}
\begin{array}{ll}
	\psi_{\beta}(\tilde{z}) - \psi_{\beta}(\bar{z}) &= h(1) - h\left(-3\right) =1 \\
	&>\left(\frac{\vartheta}{4}-\epsilon\right) \frac{4}{\vartheta} = \left(\frac{\vartheta}{4}-\epsilon\right)d(\tilde{z},\bar{z}),
\end{array}
\end{equation*}
which means that $\psi_{\beta}$ satisfies Assumption (A2) at $\{\bar{z}\}$ for $0<\delta\leq\frac{4}{\vartheta}$. Therefore, according to Theorem~\ref{thm:main_r1}, we have 
\begin{equation*}
\begin{array}{ll}
	&\sup_{\mathbb{P}\colon  \mathcal{W}_{d,1}(\mathbb{P},\mathbb{P}_N) \leq \delta }  \mathrm{E}_{\mathbb{P}} [\ell (Z;\beta) ] = \mathrm{E}_{\mathbb{P}_N}[\ell(Z;\beta)] + \frac{\vartheta}{4}\delta \\
	&= \ell(\bar{z};\beta)+\frac{\vartheta}{4}\delta,\quad 0<\delta\leq\frac{4}{\vartheta}.
\end{array}
\end{equation*}

(b2) Suppose \(\delta \geq \frac{4}{\vartheta
}\). Similar to (b1), we have \(\sup_{\mathbb{P}\colon  \mathcal{W}_{d,1}(\mathbb{P},\mathbb{P}_N) \leq \delta }  \mathrm{E}_{\mathbb{P}} [\ell (Z;\beta) ] \leq 1\). On the other hand, if we choose \(\tilde{\mathbb{P}} \coloneqq {\boldsymbol \chi}_{\{\tilde{z}\}}\) with \(\tilde{z} = \frac{1}{\vartheta}\alpha_{\beta} \), then by Lemma~\ref{lemma: singleton}, we have \(\mathcal{W}_{d,1}(\tilde{\mathbb{P}},\mathbb{P}_N) = d(\tilde{z},\bar{z}) = \frac{4}{\vartheta} \leq \delta\). Moreover, it holds that \(\mathrm{E}_{\tilde{\mathbb{P}}}[\ell(Z;\beta)]=\ell(\tilde{z};\beta)= h(1)=1\). Therefore,
\begin{equation*}
	\sup_{\mathbb{P}\colon  \mathcal{W}_{d,1}(\mathbb{P},\mathbb{P}_N) \leq \delta }  \mathrm{E}_{\mathbb{P}} [\ell (Z;\beta) ]= 1=\ell(\bar{z};\beta) + 1.
\end{equation*}		
This completes the proof.

\subsection{Proof of Example \ref{ex:linear_regression} }
\label{sec: proof_ex41}
It is obvious that the function \(\phi\colon z= (x,y) \mapsto y - \langle\beta,x\rangle \) is linear on \(\mathcal{Z}\). We are going to apply Proposition~\ref{prop:lin_loss} and Corollary~\ref{coro: linear} to draw the conclusions. Here we will check the conditions in Proposition~\ref{prop:lin_loss} and Corollary~\ref{coro: linear} case by case.

(i) Let \(\left\llbracket (x,y)\right\rrbracket = \|[x;y]\|_{\mathbb{R}^{n+1}} \), then it can be seen that \(\left\llbracket \cdot\right\rrbracket \) is absolutely homogeneous on \(\mathcal{Z}\). Then \(d((x',y'),(x,y)) = \left\llbracket (x'-x,y'-y)\right\rrbracket \) and \(\left\llbracket \cdot\right\rrbracket ^{-1}(0) = \{{\boldsymbol 0}_{\mathbb{R}^{n+1}}\} \subseteq \phi^{-1}(0)  \). In addition, 
\begin{equation*}
\begin{array}{ll}
	L_{\phi} &= \sup_{x\in\mathbb{R}^{n},y\in\mathbb{R} } \left\{ | y-\langle \beta,x\rangle | \mid \|[x;y]\|_{\mathbb{R}^{n+1} } = 1  \right\} \\
	&= \|[-\beta;1]\|_{\mathbb{R}^{n+1},*} < \infty.
\end{array}
\end{equation*}          

(ii) Let \(\left\llbracket (x,y)\right\rrbracket = \|x\|_{\mathbb{R}^{n}}+{\boldsymbol \delta}_{\{0\}}(y) \), then \(\left\llbracket \cdot\right\rrbracket \) is absolutely homogeneous,  \(d((x',y'),(x,y)) = \left\llbracket(x'-x,y'-y)\right\rrbracket \) and \(\left\llbracket \cdot\right\rrbracket ^{-1}(0) = \{{\boldsymbol 0}_{\mathbb{R}^{n+1}}\} \subseteq \phi^{-1}(0)  \). In addition, 
\begin{equation*}
\begin{array}{ll}
	L_{\phi} &= \sup_{x\in\mathbb{R}^{n},y\in\mathbb{R} } \left\{ \left| y-\langle \beta,x\rangle \right| \mid \|x\|_{\mathbb{R}^{n} } = 1, y=0  \right\} \\
	&= \|\beta\|_{\mathbb{R}^{n},*} < \infty.
\end{array}
\end{equation*}

(iii) Let \(\left\llbracket (x,y)\right\rrbracket = \|x_{\mathcal{I}}\|_{\mathbb{R}^{|\mathcal{I}|}}+{\boldsymbol \delta}_{\{{\bf 0}_{\mathbb{R}^{|\mathcal{I}^c|+1}}\}}([x_{\mathcal{I}^c};y]) \), then we have that \(\left\llbracket \cdot\right\rrbracket \) is absolutely homogeneous,  \(d((x',y'),(x,y)) = \left\llbracket (x'-x,y'-y)\right\rrbracket \) and \(\left\llbracket \cdot\right\rrbracket^{-1}(0) = \{{\boldsymbol 0}_{\mathbb{R}^{n+1}}\} \subseteq \phi^{-1}(0)  \). In addition, 
%\begin{equation*}
%\begin{array}{ll}
%	L_{\phi} & =  \sup_{\substack{ x\in\mathbb{R}^{n}\\ y\in\mathbb{R}} } \left\{ \left| y  -  \langle \beta,x\rangle \right| \mid \left\|x_{\mathcal{I}}\right\|_{\mathbb{R}^{|\mathcal{I}|} }   =   1, x_{\mathcal{I}^c}  =  0, y  =  0  \right\} \\
%	&= \left\|\beta_{\mathcal{I}}\right\|_{\mathbb{R}^{|\mathcal{I}|},*} < \infty.
%\end{array}
%\end{equation*}

\begin{equation*}
	\begin{array}{ll}
		L_{\phi} &= \sup\limits_{\substack{ x\in\mathbb{R}^{n}\\ y\in\mathbb{R}} } \left\{ \left| y-\langle \beta,x\rangle \right| \mid \left\|x_{\mathcal{I}}\right\|_{\mathbb{R}^{|\mathcal{I}|} } = 1, x_{\mathcal{I}^c}=0, y=0  \right\} \\
		&= \left\|\beta_{\mathcal{I}}\right\|_{\mathbb{R}^{|\mathcal{I}|},*} < \infty.
	\end{array}
\end{equation*}  

(iv) Let \(\left\llbracket(x,y)\right\rrbracket =\inf_{\bar{x}\in\mathbb{R}^s}\left\{\|\bar{x}\|_{\mathbb{R}^s} \mid B^T\bar{x}=x \right\} +{\boldsymbol \delta}_{\{0\}}(y)  \), then we have that \(d((x',y'),(x,y)) = \left\llbracket (x'-x,y'-y)\right\rrbracket \). It follows from \citet[Proposition 2(a)]{Chu2021OnRS} that \(\left\llbracket \cdot\right\rrbracket\) is a norm on \(\mathtt{Range}(B^T)\times \{0\}\) and infinite otherwise, hence it is absolutely homogeneous. Moreover, we can see that \(\left\llbracket \cdot\right\rrbracket ^{-1}(0) = \{{\boldsymbol 0}_{\mathbb{R}^{n+1}}\} \subseteq \phi^{-1}(0)  \). In addition, it follows from \citet[Proposition 2(c)]{Chu2021OnRS} that
%\begin{equation*}
%\begin{array}{llll}
%	&L_{\phi} \\
%	&=   \sup_{\substack{ x\in\mathbb{R}^{n}\\ y\in\mathbb{R}} }    \left\{  \left| y  -  \langle \beta,  x\rangle \right|   \middle \vert  \inf_{\bar{x}\in\mathbb{R}^s}  \left\{  \|\bar{x}\|_{\mathbb{R}^s}   \mid   B^T\bar{x}  =  x \right\}   + {\boldsymbol \delta}_{\{0\}} (y)   =   1   \right\} \\
%	&=  \sup_{x\in\mathbb{R}^{n}} \left\{ \left|\langle \beta,x\rangle \right| \mid  \inf_{\bar{x}\in\mathbb{R}^s}\left\{\|\bar{x}\|_{\mathbb{R}^s} \mid B^T\bar{x}=x \right\}   = 1  \right\}  \\
%	&=  \|B\beta\|_{\mathbb{R}^{n},*} < \infty.
%\end{array}
%\end{equation*}
\begin{equation*}
	\begin{array}{llll}
		&L_{\phi} \\
		&=  \sup\limits_{\substack{ x\in\mathbb{R}^{n}\\ y\in\mathbb{R}} }   \left\{  \left| y - \langle \beta,x\rangle \right|  \middle \vert \inf\limits_{\bar{x}\in\mathbb{R}^s} \left\{  \|\bar{x}\|_{\mathbb{R}^s}   \mid   B^T\bar{x}  =  x \right\} +{\boldsymbol \delta}_{\{0\}} (y) = 1   \right\} \\
		&=  \sup\limits_{x\in\mathbb{R}^{n}} \left\{ \left|\langle \beta,x\rangle \right| \mid  \inf\limits_{\bar{x}\in\mathbb{R}^s}\left\{\|\bar{x}\|_{\mathbb{R}^s} \mid B^T\bar{x}=x \right\}   = 1  \right\}  \\
		&= \|B\beta\|_{\mathbb{R}^{s},*} < \infty.
	\end{array}
\end{equation*}

The desired conclusions then follow from Proposition~\ref{prop:lin_loss} and Corollary~\ref{coro: linear}.

\subsection{Proof of Example \ref{eg:functional-LR}}
\label{sec: proof_ex42}
It is obvious that the function \(\phi\) defined in (a) or (b) is linear on \(\mathcal{Z} \). In addition, let \(\left\llbracket (x,y)\right\rrbracket \coloneqq \left(\int_{0}^{1}|x(t)|^2\mathrm{d}t\right)^{1/2}  + {\boldsymbol \delta}_{\{0\}}(y) \), then we can see that \(\left\llbracket \cdot\right\rrbracket \) is absolutely homogeneous on \(\mathcal{Z} \), \(d((x',y'),(x,y)) = \left\llbracket (x'-x,y'-y)\right\rrbracket \) and  
%\begin{equation*}
%\begin{array}{lll}
%	& \left\llbracket \cdot\right\rrbracket^{-1}(0) \\
%	&=  \left\{ (x,  y) \in  \mathcal{L}^2[0,1] \times \mathbb{R}  \colon  \left(\int_{0}^{1}  |x(t)|^2\mathrm{d}t\right)^{1/2}    =  0, y  =  0 \right\} \\
%	&\subseteq \phi^{-1}(0).
%\end{array}
%\end{equation*}

\begin{equation*}
\begin{array}{lll}
	& \left\llbracket \cdot\right\rrbracket^{-1}(0) \\
	&=\left\{ (x, y)\in \mathfrak{L}^2[0,1] \times \mathbb{R} :  \left(\int\limits_{0}^{1}  |x(t)|^2\mathrm{d}t\right)^{1/2}    = 0, y = 0 \right\} \\
	&\subseteq \phi^{-1}(0).
\end{array}
\end{equation*}
The desired conclusions follows from Proposition~\ref{prop:lin_loss} and Corollary~\ref{coro: linear} since we have
\begin{equation*}
\begin{array}{lll}
	&\sup\limits_{\substack{ x\in\mathfrak{L}^2[0,1]\\ y\in\mathbb{R}}} 
	   \left\{ \left|y  -  \int\limits_{0}^{1}  x(t) \beta(t) \mathrm{d}t\right|    \colon     \left(  \int\limits_{0}^{1} |x(t)|^2\mathrm{d}t \right)^{1/2}            + {\boldsymbol \delta}_{\{0\}}(y)  =  1 \right\} \\
	&=\sup\limits_{x\in\mathfrak{L}^2[0,1]} \left\{ \left|\int\limits_{0}^{1} x(t) \beta(t) \mathrm{d}t\right| \ \colon  \ \left(\int\limits_{0}^{1}|x(t)|^2\mathrm{d}t\right)^{1/2}  = 1\right\}\\
	&= \left(\int\limits_{0}^{1}\left|\beta(t)\right|^2\mathrm{d}t\right)^{1/2},
\end{array}
\end{equation*}
and
\begin{equation*}
	\begin{array}{ll}
		&\sup\limits_{\substack{ x\in\mathfrak{L}^2[0,1]\\ y\in\mathbb{R}}} 
		\left\{\begin{matrix}
			\left|y -\int_{0}^{1}x(t) \sum_{j=1}^{n} \beta_j {\boldsymbol g}_j(t) \mathrm{d}t\right|  :\\
			  \left( \int_{0}^{1} |x(t)|^2\mathrm{d}t \right)^{1/2}  +{\boldsymbol \delta}_{\{0\}}(y)=1
		\end{matrix} \right\} \\
		&= \left(\int_{0}^{1}\left|\sum_{j=1}^{n} \beta_j {\boldsymbol g}_j(t)\right|^2\mathrm{d}t\right)^{1/2}.
	\end{array}
\end{equation*}
This completes the proof.

\subsection{Proof of Example \ref{ex: logcosh}}
\label{sec: proof_ex44}
(a) We first see that
\begin{equation*}
	h'(t) = \tanh t = \frac{e^{2t}-1}{e^{2t}+1},
\end{equation*}
which means that $h$ is globally $1$-Lipschitz on $\mathbb{R}$ and for any $t_0\in \mathbb{R}$, we have
\begin{equation*}
\begin{array}{lll}
	&\lim_{k\rightarrow\infty} \frac{h(k+t_0)-h(t_0)}{k} \\
	&= \lim_{k\rightarrow\infty} \frac{\log(\cosh(k+t_0))-\log(\cosh(t_0))}{k} \\
	&= \lim_{k\rightarrow\infty} \tanh (k+t_0) = 1.
\end{array}
\end{equation*}
This is to say, \(h\) defined in (a) satisfies Assumptions (H1-H2) with  \(L_h=1\). Then the rest of the proof which involves finding \(L_{\phi}\) is similar to Example~\ref{ex:linear_regression}. We omit the details here. Then our conclusion follows from Corollary \ref{coro: nonlinear}.

(b) We have that \( h'(t) = \min \{1, \max \{-1,t \} \}  \), which means that $h$ is globally $1$-Lipschitz on $\mathbb{R}$. For any $t_0\in \mathbb{R}$, we have
\begin{equation*}
	\lim_{k\rightarrow\infty} \frac{h(k+t_0)-h(t_0)}{k} = \lim_{k\rightarrow\infty} \frac{k-h(t_0)}{k} = 1,
\end{equation*}
which means that \(h\) defined in (b) satisfies Assumptions (H1-H2) with  \(L_h=1\). The rest of the proof is similar to that of (a).

(c) It is easy to see that \(h\) defined in (c) is globally $1$-Lipschitz on $\mathbb{R}$. Moreover, for any $t_0\in \mathbb{R}$, we can see that
\begin{equation*}
\lim_{k\rightarrow\infty} \frac{h(-k+t_0)-h(t_0)}{k} = \lim_{k\rightarrow\infty} \frac{k-t_0-h(t_0)}{k} = 1,
\end{equation*}
which means that \(h\) defined in (c) satisfies Assumptions (H1-H2) with  \(L_h=1\). The rest of the proof is similar to that of (a).

\subsection{Proof of Example \ref{ex:linear-square}}
\label{appen:proof_special_reg}
We first show that \(\psi_{\beta}(z) := \langle [\beta;1],z\rangle^2\) for any $z\in \mathbb{R}^{n+1}$ satisfies Assumption (A1) with \(L_{\beta}\coloneqq \|\beta\|_2^2 +1 \). For any \(z',z\in\mathbb{R}^{n+1}\), we have 
\begin{equation*}
\begin{array}{llll}
	& \left|\psi_{\beta}(z')-\psi_{\beta}(z)\right|  = \left|\langle [\beta;1],z'\rangle^2 - \langle [\beta;1],z\rangle^2\right| \\
	&= \left|\langle [\beta;1],z'-z\rangle\right|\left|\langle [\beta;1],z'+z\rangle \right|\\
	& \leq \|[\beta;1]\|_{2}\|z'-z\|_{2} \cdot \|[\beta;1]\|_{2}\|z'+z\|_{2}\\
	& = \left(\|\beta\|_2^2 + 1 \right) d(z',z).
\end{array}
\end{equation*}
Hence, \(\psi_{\beta} \) is \((L_{\beta},d) \)-Lipschitz at $\mathbb{R}^{n+1}$. 

Next, we show that \(\psi_{\beta}\) satisfies Assumption (A2). For any \(z\in\mathbb{R}^{n+1}\) and \(k>0\), let \(\tilde{z}\coloneqq z + k \Delta  \) with \(\Delta\coloneqq \frac{[\beta;1]}{\|[\beta;1]\|_2}\), then \(\|\tilde{z}-z\|_2  = \|k\Delta\|_2 = k, \|\tilde{z}+z\|_2 = \|2z + k\Delta\|_2 \) and
\begin{equation*}
\begin{array}{ll}
	d(\tilde{z},z)& = \|\tilde{z}-z\|_2\|\tilde{z}+z\|_2 = k\|2z + k\Delta\|_2 \\
	&\geq k |k - 2\|z\|_2| \rightarrow  \infty,
\end{array}
\end{equation*}
as $k\rightarrow\infty$. On the other hand, we have
\begin{equation*}
	\frac{\left|\langle \Delta,\tilde{z}+z\rangle \right|}{\|\tilde{z}+z\|_2} \leq \|\Delta\|_2= 1,
\end{equation*}
and
\begin{equation*}
\begin{array}{ll}
	\frac{\langle \Delta,\tilde{z}+z\rangle }{\|\tilde{z}+z\|_2} &= \frac{\langle \Delta,2z+k\Delta\rangle}{\|2z+k\Delta\|_2} \\
	&= \frac{\sum_{i=1}^{n+1}\Delta_i \left(2z_i/k+\Delta_i\right)}{\sqrt{\sum_{i=1}^{n+1}\left(2z_i/k+\Delta_i\right)^2}} \rightarrow 1,
\end{array}
\end{equation*}
as $k\rightarrow\infty$. Thus, for the given $\delta$ and any \(0<\epsilon<L_{\beta}\), there exists a positive integer \(k_{\epsilon}\) such that for \(\tilde{z}_{\epsilon}=z+k_{\epsilon}\Delta\), one has 
\begin{equation*}
\begin{array}{ll}
	d(\tilde{z}_{\epsilon},z) &= k_{\epsilon}\left\|\tilde{z}_{\epsilon}+z\right\|_2 \geq \delta ,\\
	\frac{\langle\Delta,\tilde{z}_{\epsilon}+z\rangle}{\|\tilde{z}_{\epsilon}+z\|_2} &\geq 1 - \frac{\epsilon}{\|[\beta;1]\|_2^2}.
\end{array}
\end{equation*}
This implies that 
\begin{equation*}
\begin{array}{llll}
	&\psi_{\beta}(\tilde{z}_{\epsilon}) - \psi_{\beta}(z) =\langle [\beta;1],\tilde{z}_{\epsilon}-z\rangle\cdot \langle [\beta;1],\tilde{z}_{\epsilon}+z\rangle \\
	& = \|[\beta;1]\|_2^2\langle \Delta,k_{\epsilon}\Delta\rangle \langle \Delta,\tilde{z}_{\epsilon}+z\rangle \\
	& \geq   \|[\beta;1]\|_2^2 k_{\epsilon} \left(1 - \frac{\epsilon}{\|[\beta;1]\|_2^2}\right)\left\|\tilde{z}_{\epsilon}+z\right\|_2 \\
	& = \left(\|[\beta;1]\|_2^2-\epsilon\right) k_{\epsilon}\left\|\tilde{z}_{\epsilon}+z\right\|_2 = \left(L_{\beta}-\epsilon \right) d(\tilde{z}_{\epsilon},z).
\end{array}
\end{equation*}
Therefore, it satisfies Assumption~(A2). By Theorem~\ref{thm:main_r1}, we have the required result. 

\subsection{Proof of Example \ref{eg:hinge-class}}
\label{sec: proof_ex43}
It is obvious that \(\phi(\cdot) := \langle \beta,\cdot\rangle\) is linear on $\mathbb{R}^n$. Let \(\left\llbracket x\right\rrbracket = \|x\|_{\mathbb{R}^{n}}\), then we can see that $\left\llbracket \cdot\right\rrbracket$ is absolutely homogeneous,  \(d((x',y'),(x,y)) = \left\llbracket x'-x\right\rrbracket + {\boldsymbol \delta}_{\{0\}} (y'-y) \) and \(\left\llbracket \cdot\right\rrbracket ^{-1}(0) = \{{\boldsymbol 0}_{\mathbb{R}^n}\} \subseteq \phi^{-1}(0)  \). In addition, 
\begin{equation*}
\begin{array}{ll}
	&L_{\phi} = \sup_{x\in\mathbb{R}^{n},y\in\mathbb{R} } \left\{ \left|y-\langle \beta,x\rangle \right| \mid \|x\|_{\mathbb{R}^{n} } = 1, y=0  \right\} \\
	&=
	\sup_{x\in\mathbb{R}^{n}} \left\{ |\langle \beta,x\rangle | \mid \|x\|_{\mathbb{R}^{n} } = 1  \right\}
	= \|\beta\|_{\mathbb{R}^{n},*} < \infty.
\end{array}
\end{equation*}
According to Corollary \ref{coro: linear2}, we have the conclusions.

\subsection{Proof of Example \ref{exam:nonlin-class}}
\label{sec: proof_ex45}
According to Corollary \ref{coro: nonlinear2} and the proof to Example~\ref{ex:linear_regression}, it suffices to show that each \(h\) defined in this example satisfies Assumption (H1-H2) with  \(L_h=1\). Next, we discuss them one by one.

(a) For $h(t) = \log (1+\exp(-t))$, it can be seen that 
\begin{equation*}
	h'(t) = -\frac{\exp(-t)}{1+\exp(-t)}=-\frac{1}{1+\exp(t)},
\end{equation*}
which indicates that $h$ is globally $1$-Lipschitz on $\mathbb{R}$. Moreover, given any $t_0\in \mathbb{R}$, we have
\begin{equation*}
\begin{array}{lll}
	&\lim_{k\rightarrow\infty} \frac{h(-k+t_0)-h(t_0)}{k} \\
	&= \lim_{k\rightarrow\infty} \frac{\log (1+\exp(k-t_0))-h(t_0)}{k} \\
	&= \lim_{k\rightarrow\infty} \frac{\exp(k-t_0)}{1+\exp(k-t_0)}= 1.
\end{array}
\end{equation*}

(b) Given 
\begin{equation*}
	h(t)= \left\{ \begin{array}{lll}
		0 & \text{if } t \geq 1\\
		\frac{1}{2}(1-t)^2 & \text{if } 0<t<1\\
		\frac{1}{2}-t & \text{otherwise.}
		\end{array}\right.
\end{equation*}
It can be easily seen that $h$ is globally $1$-Lipschitz on $\mathbb{R}$. In addition, for any $t_0\in \mathbb{R}$, we have
\begin{equation*}
\begin{array}{ll}
	&\lim_{k\rightarrow\infty} \frac{h(-k+t_0)-h(t_0)}{k} \\
	&= \lim_{k\rightarrow\infty} \frac{\frac{1}{2}+k-t_0 - h(t_0)}{k} = 1.
\end{array}
\end{equation*}

(c) Let $h$ be defined as 
\begin{equation*}
	h(t) = \left\{ \begin{array}{lll}
		1-t & \text{if } t \geq 1\\
		-\tau_1 (1-t) & \text{if } - \tau_2 < t < 0\\
		\tau_1 \tau_2 & \text{otherwise.}
	\end{array}\right.
\end{equation*}
The piecewise linear function $h$ is obviously globally $1$-Lipschitz on $\mathbb{R}$. For any $t_0\in \mathbb{R}$, we can see that 
\begin{equation*}
\begin{array}{ll}
	&\lim_{k\rightarrow\infty} \frac{h(-k+t_0)-h(t_0)}{k} \\
	&= \lim_{k\rightarrow\infty} \frac{1+k-t_0 -h(t_0)}{k} = 1.
\end{array}
\end{equation*}
This completes the proof.

\subsection{Proof of Example \ref{ex:nusvr}}
\label{sec: proof_ex51}
According to Example \ref{ex:linear_regression} , for any $t\in \mathbb{R}$, we know that 
\begin{equation*}
\begin{array}{ll}
	&\sup_{\mathbb{P}\colon  \mathcal{W}_{d,1}(\mathbb{P},\mathbb{P}_N) \leq \delta }  \mathrm{E}_{\mathbb{P}} [\left( \left| Y - \langle\beta,X\rangle \right|-t\right)_{+} ] \\
	&= \mathrm{E}_{\mathbb{P}_N}[\left( \left| Y - \langle\beta,X\rangle \right|-t\right)_{+}] + \|[-\beta;1]\|_{\mathbb{R}^{n+1},*}\delta.
\end{array}
\end{equation*}
Define the function $G:\mathbb{R}^n\times\mathbb{R}\rightarrow \mathbb{R}$ as
\begin{equation*}
	G(Z) = \left| Y - \langle \beta,X\rangle \right|,
\end{equation*}
for any \(Z=(X,Y)\in \mathbb{R}^n\times \mathbb{R}\). Then we can see that $G$ is $(\|[-\beta;1]\|_{\mathbb{R}^{n+1},*},d)$-Lipschitz at $\mathbb{R}^n\times \mathbb{R}$. Therefore, the conclusion follows from Corollary \ref{coro: r1}(a).

\subsection{Proof of Example \ref{ex: nusvm}}
\label{sec: proof_ex52}
From Corollary \ref{coro: linear2}, we can see that for any $t\in \mathbb{R}$,
\begin{equation*}
\begin{array}{ll}
	&\sup_{\mathbb{P}\colon  \mathcal{W}_{d,1}(\mathbb{P},\mathbb{P}_N) \leq \delta }  \mathrm{E}_{\mathbb{P}} [\left(-Y\cdot\langle\beta,X\rangle- t \right)_{+} ] \\
	&= \mathrm{E}_{\mathbb{P}_N}[\left(-Y\cdot\langle \beta,X\rangle - t\right)_{+}] + \|\beta\|_{\mathbb{R}^{n},*}\delta.
\end{array}
\end{equation*}
Let $G:\mathbb{R}^n\times \mathbb{R}\rightarrow \mathbb{R}$ be defined as
\begin{equation*}
	G(Z) = -Y\cdot\langle \beta,X\rangle ,\quad \mbox{for any } Z=(X,Y)\in \mathbb{R}^n\times \mathbb{R}.
\end{equation*}
It can be seen that \(G(\cdot)\) is $(\|\beta\|_{\mathbb{R}^{n},*},d)$-Lipschitz at $\mathbb{R}^n\times \mathbb{R}$. The conclusion then  follows from Corollary \ref{coro: r1}(a).

\subsection{Proof of Example \ref{ex: highmoment_crm}}
\label{sec: proof_ex53}
The conclusion follows directly from Corollary \ref{coro: linear} and Corollary \ref{coro: r1}(b), as the function $G:\mathbb{R}^n\rightarrow \mathbb{R}$ defined as
\begin{equation*}
	G(Z) = \langle \beta,X\rangle ,\quad \mbox{for any } Z\in \mathbb{R}^n,
\end{equation*}
is $(\|\beta\|_{\mathbb{R}^{n},*},d)$-Lipschitz at $\mathbb{R}^n$.

\section{A weaker version of Theorem~\ref{thm:main_r1} with relaxed assumptions}
\label{sec: weaker_main_r1}
\begin{theorem} \label{thm:weaker_main_r1}
	Let \(\mathcal{Z}_N\coloneqq\{Z^{(1)},\dots,Z^{(N)}\}\subset\mathcal{Z} \) be a given dataset and \(\mathbb{P}_N\coloneqq\sum_{i=1}^{N}\mu_i{\boldsymbol \chi}_{\{Z^{(i)}\}}\in\mathcal{P}(\mathcal{Z})\) be the corresponding empirical distribution. In addition, let \(d(\cdot,\cdot)\) be a cost function on \(\mathcal{Z}\times\mathcal{Z}\) and  \(\delta\in(0,\infty)\) be a scalar. Suppose the loss function $\ell:\mathcal{Z}\times \mathcal{B}\rightarrow \mathbb{R}$ takes the form as
	\begin{equation*}
		\ell\colon (z;\beta) \mapsto \psi_{\beta}(z), 
	\end{equation*}
	where the function $\psi_{\beta}\colon\mathcal{Z}\rightarrow \mathbb{R}$ satisfies the following assumptions:
	\begin{enumerate}[label=(C\arabic*)]
		\item \(\psi_{\beta}\) is \((L_{\beta}^{\{Z^{(i)}\}},d)\)-Lipschitz at \(\{Z^{(i)}\}\) with \(L_{\beta}^{\{Z^{(i)}\}}\in(0,\infty)\) for each \(1\leq i\leq N\);
		\item for any \(\epsilon\in (0,\min_{i}L_{\beta}^{\{Z^{(i)}\}} )\) and each \(Z^{(i)}\in\mathcal{Z}_N\), there exists \(\tilde{Z}^{(i)}_{\epsilon}\in\mathcal{Z}\) such that \(\delta\leq d(\tilde{Z}^{(i)}_{\epsilon},Z^{(i)})<\infty \) and
		\begin{equation*}
		\psi_{\beta}(\tilde{Z}^{(i)}_{\epsilon}) - \psi_{\beta}(Z^{(i)}) \geq   (L_{\beta}^{\{Z^{(i)}\}}-\epsilon) d(\tilde{Z}^{(i)}_{\epsilon},Z^{(i)}).
		\end{equation*}
	\end{enumerate}
	Then we have that \(\widehat{\mathcal{L}}\leq \mathcal{S}\leq \widehat{\mathcal{U}} \), where
	\begin{equation*} 
	\begin{array}{ll}
		& \widehat{\mathcal{L}}=\mathrm{E}_{\mathbb{P}_N}[\ell(Z;\beta)] + \sum_{i=1}^{N} \mu_i L_{\beta}^{\{Z^{(i)}\}}\delta,\\
		& \widehat{\mathcal{U}}=  \mathrm{E}_{\mathbb{P}_N}[\ell(Z;\beta)] + \max_{i=1,\dots,N} L_{\beta}^{\{Z^{(i)}\}}\delta,
	\end{array}
	\end{equation*}
	which means that
	\begin{equation*} 
	\begin{array}{lll}
		&\mathrm{E}_{\mathbb{P}_N}[\ell(Z;\beta)] + \sum_{i=1}^{N} \mu_i L_{\beta}^{\{Z^{(i)}\}}\delta \\
		&\leq \sup_{\mathbb{P}\colon  \mathcal{W}_{d,1}(\mathbb{P},\mathbb{P}_N) \leq \delta }  \mathrm{E}_{\mathbb{P}} [\ell (Z;\beta) ] \\
		&\leq \mathrm{E}_{\mathbb{P}_N}[\ell(Z;\beta)] + \max_{i=1,\dots,N} L_{\beta}^{\{Z^{(i)}\}}\delta.
	\end{array}
	\end{equation*}
\end{theorem}
\textbf{Proof.} 
Since \(\psi_{\beta}\) is \((L_{\beta}^{\{Z^{(i)}\}},d) \)-Lipschitz at \(\{Z^{(i)}\}\) for each \(1\leq i\leq N\), it implies that \(\psi_{\beta}\) is \((L_{\beta}^{\mathcal{Z}_N},d) \)-Lipschitz at \(\mathcal{Z}_N\) with
\begin{equation*}
	L_{\beta}^{\mathcal{Z}_N} = \max_{i=1,\dots,N} L_{\beta}^{\{Z^{(i)}\}}.
\end{equation*}
By Theorem~\ref{thm:main}, by letting \(\mathcal{L}_i=\sup_{\mathbb{P}\in \mathcal{P}(\mathcal{Z})} \left\{\mathrm{E}_{\mathbb{P}}[\ell(Z;\beta)] \middle\vert  \mathcal{W}_{d,1}(\mathbb{P},{\boldsymbol \chi}_{\{Z^{(i)}\}}) \leq \delta \right\}\), \(i=1,\cdots,N\), we have that 
\begin{equation*}
\begin{array}{ll}
	&\sum_{i=1}^{N}\mu_i \mathcal{L}_i \leq \mathcal{S}=\sup_{\mathbb{P}\colon  \mathcal{W}_{d,1}(\mathbb{P},\mathbb{P}_N) \leq \delta }  \mathrm{E}_{\mathbb{P}} [\ell (Z;\beta) ] \\
	& \leq \mathrm{E}_{\mathbb{P}_N}[\ell(Z;\beta)] + L_{\beta}^{\mathcal{Z}_N}\delta =\widehat{\mathcal{U}}.
\end{array}
\end{equation*}
Then by applying Theorem~\ref{thm:main_r1} for each \(\mathcal{L}_i\), we can see that
\begin{equation*}
	\mathcal{L}_i = \ell(Z^{(i)};\beta) + L_{\beta}^{\{Z^{(i)}\}}\delta,
\end{equation*}
which means that \(\sum_{i=1}^{N}\mu_i \mathcal{L}_i = \sum_{i=1}^{N}\mu_i \ell(Z^{(i)};\beta) + \sum_{i=1}^{N}\mu_i L_{\beta}^{\{Z^{(i)}\}}\delta =\widehat{\mathcal{L}} \). This completes the proof.
\hfill\(\square\)
%\end{APPENDICES}

%\theendnotes

%\ACKNOWLEDGMENT{Meixia Lin is supported by The Singapore University of Technology and Design under MOE Tier 1 Grant SKI 2021{\_}02{\_}08. Kim-Chuan Toh is supported by the Ministry of Education, Singapore, under its Academic Research Fund Tier 3 grant call (MOE-2019-T3-1-010).}	

\end{document}